\documentclass[%
aip,
cha,
amsmath,amssymb,
reprint,%
]{revtex4-2}

\usepackage{graphicx}
\usepackage{dcolumn}
\usepackage{bm}

\usepackage[utf8]{inputenc}
\usepackage[T1]{fontenc}
\usepackage{mathptmx}
\usepackage{etoolbox}

\makeatletter
\def\@email#1#2{%
	\endgroup
	\patchcmd{\titleblock@produce}
	{\frontmatter@RRAPformat}
	{\frontmatter@RRAPformat{\produce@RRAP{*#1\href{mailto:#2}{#2}}}\frontmatter@RRAPformat}
	{}{}
}%
\makeatother

\begin{document}

\title[Reconstruction of phase-amplitude dynamics from electrophysiological signals]{Reconstruction of phase-amplitude dynamics from electrophysiological signals}
	\author{Azamat Yeldesbay}
	\affiliation{University of Cologne, Institute of Zoology, Cologne, Germany}
    \affiliation{Cognitive Neuroscience, Institute of Neuroscience and Medicine (INM-3),  Research Centre J\"{u}lich, J\"{u}lich, Germany}
    \email{a.yeldesbay@fz-juelich.de}
	\author{Gemma Huguet}%
	\affiliation{Universitat Polit\`{e}cnica de Catalunya, Departament de Matem\`{a}tiques, Barcelona, Spain}%
    \affiliation{Centre de Recerca Matem\`{a}tica, Barcelona, Spain}
    \email{gemma.huguet@upc.edu}
 	\author{Silvia Daun}
	\affiliation{Cognitive Neuroscience, Institute of Neuroscience and Medicine (INM-3), Research Centre J\"{u}lich, J\"{u}lich, Germany}%
    \affiliation{University of Cologne, Institute of Zoology, Cologne, Germany}
    \email{s.daun@fz-juelich.de}
    \date{\today}








\begin{abstract}
We present a novel method of reconstructing the phase-amplitude dynamics directly from measured electrophysiological signals to estimate the coupling between brain regions. For this purpose, we use the recent advances in the field of phase-amplitude reduction of oscillatory systems, which allow the representation of an uncoupled oscillatory system as a phase-amplitude oscillator in a unique form using transformations (parameterizations) related to the eigenfunctions of the Koopman operator. By combining the parameterization method and the Fourier-Laplace averaging method for finding the eigenfunctions of the Koopman operator, we developed a method of assessing the transformation functions from the signals of the interacting oscillatory systems. The resulting reconstructed dynamical system is a network of phase-amplitude oscillators with the interactions between them represented as coupling functions in phase and amplitude coordinates.
\end{abstract}

\maketitle

\begin{quotation}
Signals from interacting brain regions display transient synchronization of phases and amplitudes in different frequencies. Commonly, the interaction between regions of the brain is quantitatively described by either analyzing the correlations of amplitudes of the measured signals or by calculating phase-synchronization measures. However, for a complete picture of the interactions it is important to analyze the dynamics of both the amplitude and the phase. 

In this work, we present a new method for finding the interaction between brain regions by reconstructing the phase-amplitude dynamics directly from the measured electrophysiological signals. The method is based on the representation of the oscillatory activities in the brain as phase-amplitude oscillators in a unique form using transformations (parameterizations) related to the eigenfunctions of the Koopman operator. We combined the parameterization method \cite{Cabre2003a, Cabre2003b, Cabre2005, Guillamon2009, Huguet2013} and the Fourier-Laplace averaging method \cite{Mauroy2012, Mauroy2018} to develop a novel method of estimating the transformation functions from the oscillatory activities in the measured electrophysiological signals. Thereby, we reconstruct a dynamical system as a network of phase-amplitude oscillators with the edges representing coupling functions in phase and amplitude coordinates. 

Using synthetic signals generated from several models with known and unknown theoretical phase-amplitude reduced forms, we demonstrate that our method is capable of finding the proper unique dynamic form of these oscillatory systems in the reduced phase-amplitude space. 

Our method can be applied to describe any network of interacting oscillators as a dynamical system using signals of the network elements. In particular, to analyze the coupling between distant brain regions using high time resolution signals, such as electroencephalography (EEG) or magnetoencephalography (MEG). Further simulation and study of the reconstructed dynamical system then enables the construction and investigation of a mathematical model of various neural pathologies and disorders of the brain. 
\end{quotation}

Oscillators, synchronization, experimental data, network, communication between brain regions	


\section{Introduction}

Performing a cognitive task involves large-scale integration of multiple brain regions, which is accompanied by an emergence of transient synchronization over multiple frequency bands \cite{Varela2001, Buzsaki2004, VanEde2018}. The transient synchronization observed during the stimulus of a cognitive task, e.g. the processing of a movement, its preparation and execution stages, is functionally relevant for the activation and coordination of a widespread network of oscillatory neural populations across brain areas \cite{Engel2001, Buzsaki2004}. A detailed study of the emergence mechanisms of the transient synchronization is important to understand the interaction between brain areas involved in such a task. 

An assessment of the coupling between distant brain regions by means of the synchronization of oscillatory activity observed in EEG/MEG signals can be performed in different ways \cite{Greenblatt2012b,LopesdaSilva2013}: by considering the relation between the phase of the signals and calculating phase synchronization measures \cite{Lachaux1999a, Sauseng2008a, Sauseng2008b}, by analyzing the relationship between amplitudes with respect to an event using event-related synchronization (ERS) and desynchronization (ERD)  measures \cite{Pfurtscheller1999}, or by calculating phase-amplitude coupling (PAC) measures \cite{Tort2010,Voytek2013,Hulsemann2019}, which consider the correlation between phase and amplitude of different frequency bands. All these conventional measures provide only an estimation of the interaction magnitude between the oscillatory sources, but do not provide information about the direction of the coupling and can not explain the mechanism of the transient synchronization. The latter is possible by analyzing a 
dynamical system that describes the oscillatory network and reproduces the observed transient synchronization.

The reconstruction of a dynamical system in a form of coupled phase oscillators from observations was first presented in Ref.~\onlinecite{Rosenblum2001} and was developed further in the following works Refs.~\onlinecite{Kralemann2007,Kralemann2008,Kralemann2011,Kralemann2013,Kralemann2014,Penny2009,Stankovski2017,Stankovski2019,Yeldesbay2019}. These methods are based on the phase reduction approach \cite{Hoppensteadt1997, Pikovsky2001, Pikovsky2015}, which is usually limited by weak coupling between oscillators, which, however, is usually not the case for the connection between brain regions. For strong coupling cases consideration of the amplitude of the oscillations is required. In the recent decades several approaches of phase-amplitude reduction were developed \cite{Monga2019}, namely methods based on the parameterization approach \cite{Cabre2003a, Cabre2003b, Cabre2005, Guillamon2009, Huguet2013, Perez-Cervera2019, Perez-Cervera2020}, the Koopman operator \cite{Mezic2004, Mezic2005, Budisic2012, Mauroy2012, Mauroy2018} and  the isostable coordinates approaches \cite{Mauroy2013,Wilson2016,Wilson2018,Wilson2020a,Wilson2021a}, with further development as data-driven phase-amplitude reduction approaches using isostable coordinates \cite{Wilson2020b, Wilson2021b, Ahmed2022, Cestnik2022, Wilson2023}, using dynamic mode decomposition (DMD) \cite{Tu2014, Williams2015,Proctor2016, Proctor2018, Kaiser2021}, and as an extension of dynamic causal modelling \cite{Fagerholm2020}. 

\begin{figure*}
	\centering
	\includegraphics[width=\textwidth]{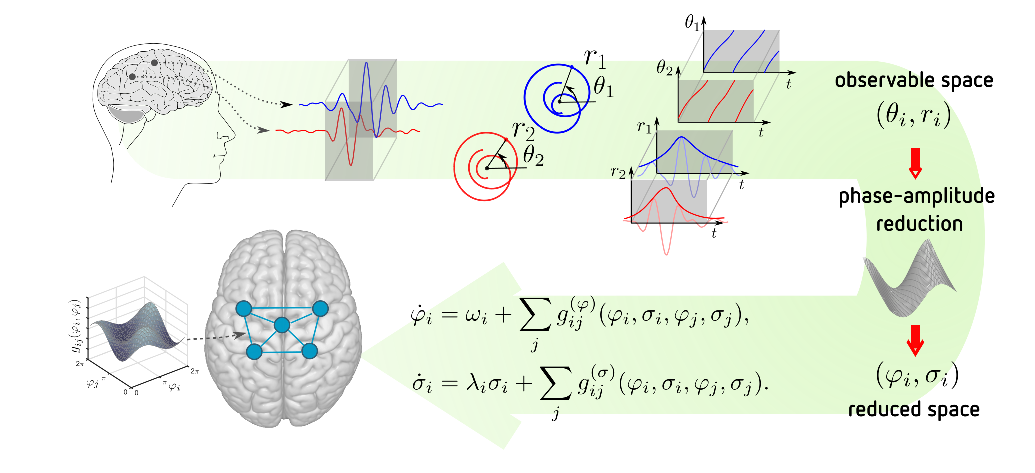}
	\caption{From electrophysiological signals to a network model. We extract the amplitude and phase of the brain signals of different sources obtained from experimental measurements during the performance of a cognitive task. The extracted (observable) amplitude and phase $(\theta_i,r_i)$ for each region $i$ are transformed into the reduced amplitude and phase $(\varphi_i,\sigma_i)$ by means of the phase-amplitude reduction. The reconstructed network is a dynamical system, where the edges are given by coupling functions. The parameters of the reconstructed system and of the transformation from observable to reduced space are obtained from the experimental data. \label{fig:introScheme}}
\end{figure*}

In this work we aim at developing a new method to build a dynamical system that will be able to describe the coupling between regions of the brain (the sources of measured signals) in the framework of transient synchronization, reflected in phase and amplitude of the measured signals (Fig.~\ref{fig:introScheme}). We achieve this by reconstructing the network of dynamical systems in a unique (invariant) form using the amplitude and phase of the measured signals by adapting the recent advances in the phase-amplitude reduction theory. In particular, the phase-amplitude reduction performs a transformation of observable variables (i.e. the phase and amplitude of the measured signals $(\theta_i,r_i)$, see Fig.~\ref{fig:introScheme}) into uniquely defined reduced variables (the reduced phase and amplitude $(\varphi_i,\sigma_i)$, see Fig.~\ref{fig:introScheme}). This reduction allows us to find a unique dynamical description of the network and the coupling terms between the regions. Thereby, we make several assumptions: (1) during the performance of a cognitive task some intrinsic processes within the cortical column of a brain region activate oscillatory activities in different frequencies, (2) in the considered interval of time this oscillating system has a limit cycle and the amplitude and phase are determinable around this limit cycle, (3) the interaction between brain areas is reflected in the transient synchronization of phase and amplitude of the measured signals.  

The outline of this article is the following. We begin by introducing the key  concepts of the phase-amplitude reduction and then divide the methods section into two parts. In the first part, we provide a brief overview of the theoretical methods of phase-amplitude reduction for uncoupled oscillators with a known vector field (VF). There, we introduce the parameterization and the Koopman operator approaches. Further, based on these approaches we present a numerical method for determining the transformation functions that perform the phase-amplitude reduction for systems with a known VF. In the second part of the methods, we introduce a method for reconstructing the phase-amplitude dynamics of an observed system using measured signals. Here, we approximate the system's VF from the measurements and find the coupling functions as a result of the phase-amplitude reduction. In the results section, we demonstrate the applicability of the phase-amplitude reduction method to several models of interconnected oscillators with different amplitude and phase properties. Finally, we conclude with a discussion of the results and outline potential future applications of the method.

\section{Methods}

\subsection{Introduction to phase-amplitude reduction}
In this section, we introduce the amplitude and phase reduction method and review several concepts that will form the basis for the methods presented further in this study. 

\subsubsection{Koopman operator and phase-amplitude reduction \label{sec:Koopman}}

Consider a dynamical system described by an autonomous system of ODEs
\begin{equation}
	\label{eq:1}
	\dot{\mathbf{x}}=\mathbf{F(\mathbf{x})}, \quad \mathbf{x} \in  \mathbb{R}^n, 
\end{equation}
which has a limit cycle $\Gamma$ of period $T$.
We assume that the limit cycle is hyperbolic stable, that is, its Floquet exponents satisfy $\Lambda_1=0$ and $\mathcal{R} \{\Lambda_j\}<0$, for $j=2,\ldots,n$. We denote by $\mathcal{M}$ the stable manifold of the limit cycle $\Gamma$, which coincides with its basin of attraction.

We introduce the phase-amplitude reduction following Ref.~\onlinecite{Mauroy2018}. Let us assume that there exists a local diffeomorphism $\mathcal{Z}: \mathcal{M} \subset \mathbb{R}^n \rightarrow \mathbb{C}^n$ that transforms the state space coordinates $\mathbf{x}$ into new coordinates $\mathbf{z}$, such that the dynamics in the new coordinates are given by
\begin{equation}
	\label{eq:zdot}
	\dot{z}_i = \mu_i z_i, \; i=1,\ldots, n,
 \end{equation}
where $\mu_1=i \omega$ and $\mu_j=\Lambda_j$, $j=2,\ldots,n$.

Then, we redefine the coordinates in the following way. If $\mu_i \in \mathbb{R}$, we define a new coordinate $\sigma_i = z_i \in \mathbb{R}$ and denote $\lambda_i \triangleq \mu_i$. If $\mu_i \notin \mathbb{R}$ we define $\sigma_i = |z_i| \in \mathbb{R}$ and $\varphi_i = \angle z_i \in [0,2\pi)$, and also denote $\lambda_i \triangleq \mathcal{R}\{\mu_i\}$ and $\omega_i \triangleq \mathcal{I}\{\mu_i\}$. By removing repeated variables due to complex conjugate pairs, we obtain an $n$-dimensional system of ODEs in the amplitude and angle variables given by
\begin{equation}
	\label{eq:phampn}
	\begin{array}{rcl}
		\dot{\varphi}_k &=& \omega_k, \\
		\dot{\sigma}_l &=& \lambda_l \sigma_l,  
	\end{array}
\end{equation}
where $k=1,\ldots,m+1$ and $l=2,\ldots,n-m-1$, and $m$ is the number of pairs of complex conjugate Floquet exponents $\Lambda_i$'s. 
This transformation of the system in Eq.~(\ref{eq:1}) into the linear system in Eq.~(\ref{eq:phampn}) is called phase-amplitude reduction \cite{Mauroy2018} and corresponds to the parameterization method presented in Refs.~\onlinecite{Cabre2005, Guillamon2009, Perez-Cervera2020}. To perform the phase-amplitude reduction, we should find a transformation from the state space coordinates  $\mathbf{x} \in \mathcal{M} \subset \mathbb{R}^n$ to a new phase and amplitude coordinates $\chi= (\varphi_1,\ldots,\varphi_{m+1},\sigma_2,\ldots,\sigma_{n-m-1}) \in  [0,2\pi)^{m+1} \times \mathbb{R}^{n-m-1}$, where the dynamics is linear as given in Eq.~(\ref{eq:phampn}). We will refer to this space as the reduced space.

This transformation is closely related to the Koopman operator as shown in Refs.~\onlinecite{Budisic2012,Mauroy2012,Mauroy2018}. 
The group of Koopman operators associated with system in Eq.~(\ref{eq:1}) is a linear one-parameter semi-group of operators $\{U^t\}_{t\geq 0}: \mathcal{F} \rightarrow \mathcal{F}$ acting on the space $\mathcal{F}$ of scalar functions (observables) $g: \mathcal{M} \subset \mathbb{R}^n \rightarrow \mathbb{R}$ as follows \cite{Mauroy2018, Mauroy2020}:
\begin{equation}\label{eq:Ut}
	U^t g  = g \circ s^t ,
\end{equation} 
where $s^t \triangleq s(t,\cdot): \mathcal{M} \rightarrow \mathcal{M}$ is the flow of system in Eq.~(\ref{eq:1}), i.e $s^t(\mathbf{x})=s(t,\mathbf{x})$ is the solution of system in Eq.~(\ref{eq:1}) with  initial condition $\mathbf{x} \in \mathcal{M}$. In other words, if the state $\mathbf{x}$ after time $t$ transfers into $s(t,\mathbf{x})$ due to the dynamics of system in Eq.~(\ref{eq:1}), then the observable $g(\mathbf{x})$ of the system transfers into $g(s(t,\mathbf{x}))$ due to the linear operator $U^t$. Thus, the Koopman operator \textit{lifts} the dynamics from the state space ($\mathbf{x}$) to the observable space ($g(\mathbf{x})$) such that the dynamics in the reduced space is linear but infinite-dimensional. 

Due to linearity we can define the eigenfunctions of the Koopman operator as the non-zero functions $\phi_{\mu} \in \mathcal{F}$ satisfying
\begin{equation}\label{eq:eigen}
	U^t \phi_{\mu}  = e^{\mu t} \phi_{\mu}, \; \forall t\geq 0,
\end{equation}
where $\mu \in \mathbb{C}$ is the corresponding eigenvalue \cite{Mauroy2018, Mauroy2020}. If we define the new variable $z(t) = \phi_{\mu}(s(t,\mathbf{x}))$, then using Eq.~\eqref{eq:Ut} and Eq.~\eqref{eq:eigen}, we have
\begin{multline}
	\dot{z}(t) = \frac{d}{dt} \phi_{\mu}(s(t,\mathbf{x})) = \frac{d}{dt} U^t \phi_{\mu}(\mathbf{x}) = \frac{d}{dt} e^{\mu t} \phi_{\mu}(\mathbf{x}) \\ = \mu  e^{\mu t} \phi_{\mu}(\mathbf{x}) = \mu z(t).
\end{multline}
Thus, the eigenfunctions of the Koopman operator transform the dynamics of the system in Eq.~\eqref{eq:1} into the linear dynamics given in Eq.~\eqref{eq:zdot}, and therefore serve as the transformation functions that perform the phase-amplitude reduction. 
Indeed, the Koopman operator associated with the system in Eq.~(\ref{eq:1}) is well characterized \cite{Mauroy2020}. The principal eigenvalues of the Koopman operator are given by $\mu_1=i \omega$ and the Floquet exponents $\mu_j=\Lambda_j$ ($j=2,\ldots,n$) with associated eigenfunctions $\phi_{\mu_j}$ that have support on the basin of attraction of the limit cycle and are continuously differentiable in the interior. Thus, the diffeomorphism $\mathcal{Z}$ is defined as $\mathcal{Z}(\mathbf{x})=(\phi_{\mu_1}(\mathbf{x}),\ldots, \phi_{\mu_n}(\mathbf{x}))$. 

The phase variable $\varphi_1$ in Eq.~(\ref{eq:phampn}) describes the periodic dynamics along the limit cycle $\Gamma$, whereas the other phase $\varphi_k$ ($k>1$) and amplitude $\sigma_l$ ($l>1$) variables describe the transient dynamics towards the limit cycle. Moreover, on the limit cycle, we have $\sigma_l=0$ for all $l$. In our study we assume that the dynamics of the considered system can be approximated as planar, meaning that all amplitude variables except one, e.g. $\sigma_2$, decay rapidly to the limit cycle (i.e. $0>\mathcal{R}\{\Lambda_2\} \gg \mathcal{R}\{\Lambda_j\}$ for $j>2$), thereby leaving only two variables - the phase $\varphi_1$ and amplitude $\sigma_2$ variable - that describe the state of the system. Further in the text we use these variables without indices. 

The dynamics in the reduced space then reads as
\begin{equation}
	\label{eq:phisig}
	\left\{
	\begin{array}{cc}
		\dot{\varphi} = \omega, \\
		\dot{\sigma} = \lambda \sigma,
	\end{array}	 
	\right.
\end{equation}
where $\omega=2 \pi/T$ is the radial frequency of the oscillation and $\lambda$ is the real part of the Floquet exponent with smallest modulus. 

It needs to be mentioned that the definition of the coordinates $(\varphi,\sigma)$ in the reduced space is not unique. Any phase variable $\hat{\varphi} = \varphi + \vartheta$, where $\vartheta \in [0,2\pi)$ and any amplitude variable $\hat{\sigma} = b\sigma$, where $b\in \mathbb{R}$ also satisfy Eq.~(\ref{eq:phisig}). This uncertainty can be solved if we fix $\varphi=0$ anywhere on the limit cycle, and set a condition to the gradient of the transformation function with respect to $\sigma$, as it will be mentioned later.

\subsubsection{parameterization method and the invariance equation} 
Here we present the parameterization method which was introduced in Ref.~\onlinecite{Cabre2005} and used in Refs.~\onlinecite{Guillamon2009, Castejon2013, Perez-Cervera2020} as a phase-amplitude reduction method. Consider a planar system (for higher dimensions see Ref.~\onlinecite{Perez-Cervera2020}) for which we know the transformation from the phase $\varphi$ and amplitude $\sigma$ variables into the state variables $\mathbf{x} \triangleq (x,y) \in \mathbb{R}^2$, i.e. $\mathbf{x} = K(\varphi,\sigma)$. Using the chain rule we have
\begin{equation}
	\label{eq:invEq}
	\dot{\mathbf{x}} = \frac{\partial K}{\partial \varphi} \dot{\varphi} + \frac{\partial K}{\partial \sigma} \dot{\sigma},
\end{equation}
and, using Eq.~(\ref{eq:1}) with $n=2$ and Eq.~(\ref{eq:phisig}), we obtain
\begin{equation}
	\label{eq:invK}
	\omega\frac{\partial K}{\partial \varphi}  + \lambda \sigma \frac{\partial K}{\partial \sigma}  = \mathbf{F}(K(\varphi,\sigma)).
\end{equation}
The equation above is known as the invariance equation \cite{Cabre2005, Guillamon2009}. One can also derive the inverse invariance equation as follows
\begin{eqnarray}
	\label{eq:invPhi}
	F^{(x)}(\mathbf{x})\frac{\partial \Phi(\mathbf{x})}{\partial x}  + F^{(y)}(\mathbf{x})\frac{\partial \Phi(\mathbf{x})}{\partial y}  = \omega, \\ \label{eq:invSig}
	F^{(x)}(\mathbf{x})\frac{\partial \Sigma(\mathbf{x})}{\partial x}  + F^{(y)}(\mathbf{x})\frac{\partial \Sigma(\mathbf{x})}{\partial y}  = \lambda \Sigma(\mathbf{x}), 
\end{eqnarray}
where $F^{(x)}$, $F^{(y)}$ are the corresponding components of $\mathbf{F}$, while $\Phi(\mathbf{x})$ and $\Sigma(\mathbf{x})$ are the phase and amplitude transformations from the state space into the reduced space. 

Note that Eqs.~(\ref{eq:invK}),(\ref{eq:invPhi}) and (\ref{eq:invSig}) can also be expressed in polar coordinates, i.e. $\mathbf{x}\triangleq (\theta,r)$. Then, Eqs.~(\ref{eq:invPhi}) and (\ref{eq:invSig}) read as 
\begin{eqnarray}
	\label{eq:invPhiPol}
    F^{(\theta)}(\mathbf{x})\frac{\partial \Phi(\mathbf{x})}{\partial \theta}  + F^{(r)}(\mathbf{x})\frac{\partial \Phi(\mathbf{x})}{\partial r}= \omega, \\ \label{eq:invSigPol}
    F^{(\theta)}(\mathbf{x})\frac{\partial \Sigma(\mathbf{x})}{\partial \theta} + F^{(r)}(\mathbf{x})\frac{\partial \Sigma(\mathbf{x})}{\partial r}  = \lambda \Sigma(\mathbf{x}), 
\end{eqnarray}
where $F^{(r)}(\mathbf{x})$ and $F^{(\theta)}(\mathbf{x})$ are the components of $\mathbf{F}$ in polar coordinates. 

\begin{figure*}
	\centering
	\includegraphics[width=0.75\textwidth]{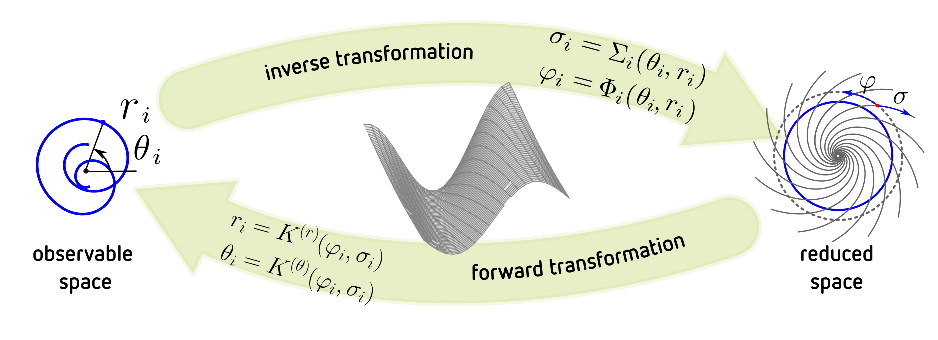}
	\caption{Phase-amplitude transformations. The transformation that performs the phase-amplitude reduction from the observable space $(\theta,r)$ to the reduced space $(\varphi,\sigma)$ is \textit{the inverse transformation}. The transformation that changes the coordinates from the reduced space to the observable space is \textit{the forward transformation}. \label{fig:expTransf}}
\end{figure*}

To further distinguish the transformations, we will use the following notation. The transformation $K(\varphi,\sigma)$, that maps the variables $(\varphi,\sigma)$ from the reduced space to the variables $\mathbf{x}=(\theta,r)$ of the state space (see Fig.~\ref{fig:expTransf}) is called \textit{forward transformation}. We use superscripts to indicate its components, for instance $K^{(x)}(\varphi,\sigma)$. The transformation $(\Phi(\mathbf{x}), \Sigma(\mathbf{x}))$ that maps the state space variables into the reduced space variables is called \emph{inverse transformation}. The latter is associated to the eigenfunctions $\phi(\mathbf{x})$ of the Koopman operator defined up to some constant multiplier or shift as discussed above.

It was shown in Refs.~\onlinecite{Guillamon2009, Huguet2013} that the forward transformation (the parameterization) can be approximated around the limit cycle by Fourier-Taylor series, that is,
\begin{equation}
	\label{eq:Kapprox}
	K(\varphi,\sigma) = \sum^{\infty}_{n=0} K_n(\varphi) \sigma^n,
\end{equation}
where $K_n(\varphi)$ are $2\pi$-periodic functions in $\varphi$, which are approximated using Fourier series as follows
\begin{equation}
	\label{eq:Knapprox}
	 K_n(\varphi)= \sum_{k=-\infty}^{\infty} q_{n,k} e^{\mathrm{i} k \varphi}. 
\end{equation}

Using the invariance equation in Eq.~(\ref{eq:invK}), it is easy to see that $K_0(\varphi)$ is a parameterization of the limit cycle $\Gamma$. 

Following Ref.~\onlinecite{Guillamon2009}, we insert Eq.~(\ref{eq:Kapprox}) into Eq.~(\ref{eq:invK}) and expand the expression using Taylor series around the limit cycle, i.e. $\sigma =0$. Then, by collecting coefficients of $\sigma$ on both sides, we obtain a homogeneous linear  differential equation for $K_1(\varphi)$:
\begin{equation}
	\label{eq:K1}
	\omega\frac{d}{d \varphi} K_1(\varphi) + \lambda  K_1(\varphi)  = D\mathbf{F}(K_0(\varphi)) \, K_1(\varphi), 
\end{equation} 
with $\omega=2 \pi/T$. Thereby, $K_1$ is an eigenfunction of the linear differential operator 
$$
\mathcal{L} \triangleq \omega \frac{d}{d \varphi} -  D\mathbf{F}(K_0(\varphi)),$$
with eigenvalue $-\lambda$. 
As it is shown in Refs.~\onlinecite{Cabre2005}~and~\onlinecite{Guillamon2009}, we have that $K_1(\varphi)$ is a periodic eigenfunction of $\mathcal{L}$ with eigenvalue $-\lambda $ if and only if $K_1(0)$ is an eigenvector of the monodromy matrix $\mathbf{M}_{\varphi=2\pi}$ with eigenvalue $e^{\lambda T}$. The matrix $\mathbf{M}_{\varphi=2 \pi}$ is a fundamental solution of the first variational equation
\begin{equation}
	\label{eq:monmatrix}
	\frac{d}{d \varphi} \mathbf{M}_{\varphi} = \frac{T}{2\pi} D\mathbf{F}(K_0(\varphi)) \mathbf{M}_{\varphi},
\end{equation}
evaluated at $\varphi=2\pi$, with initial condition $\mathbf{M}_{\varphi=0}=\mathbf{I}$. Once $K_1(0)$ is chosen as an eigenvector, $K_1(\varphi)$ is determined for all $\varphi \in [0,2\pi)$ by solving Eq.~\eqref{eq:K1}, that is, $K_1(\varphi) = e^{-\lambda T/(2 \pi) \varphi} M_{\varphi}K_1(0)$.

By repeating the procedure of inserting Eq.~(\ref{eq:Kapprox}) into Eq.~(\ref{eq:invK}) and collecting terms with $\sigma^n$, one can derive linear non-homogeneous differential equations for $K_n(\varphi)$, $n\geq2$ which depend on $K_j(\varphi), j<n$, thereby finding the higher terms of Eq.~(\ref{eq:Kapprox}) as shown in Ref.~\onlinecite{Guillamon2009}. Moreover, as it is shown in Refs.~\onlinecite{Castejon2013, Perez-Cervera2020}, using the expression for the parameterization $K(\varphi,\sigma)$ (the forward transformation function) one can find the expression for the inverse transformation functions $\Phi(\mathbf{x})$ and $\Sigma(\mathbf{x})$.

\subsection{Phase-amplitude reduction for an uncoupled oscillatory system}

In this section we describe how to find the transformation that performs the phase-amplitude reduction of an uncoupled system. For this computation, we assume that a good approximation of the 
vector field (VF)
of the uncoupled system is available in a certain region (not necessarily small) around the limit cycle in the observable space $(\theta,r)$, where $\theta$ and $r$ represent the angle and the radius of a signal measured from the uncoupled system. The system dynamics are given by: 
\begin{eqnarray*}
	\dot{\theta}&= F ^{(\theta)}(\theta,r),\\
	\dot{r}&= F ^{(r)}(\theta,r), 
\end{eqnarray*}
where $F ^{(\theta)}(\theta,r)$ and $F^{(r)}(\theta,r)$ are the two components of the vector field $\mathbf{F}(\theta,r)$ expressed in polar coordinates. The method for obtaining a good approximation of the vector field will be discussed in the next section.  It is important to note that knowing the vector field in polar coordinates is equivalent to knowing it in Cartesian coordinates $(x,y)$ via standard transformations. 

Furthermore, we assume that the approximated system has a limit cycle, which could be represented as $\{(\theta,r(\theta)), \theta \in [0,2\pi)\}$, where $r(\theta)$ is a periodic function that can be approximated using Fourier series by
\begin{equation}\label{eq:gam0rth}
 \gamma^{(r)}(\theta) = \sum_{k=-N_{g}}^{N_g} g_k e^{\mathrm{i}k\theta}.
\end{equation}

We emphasize that the observable variables are the radius $r$ and angle $\theta$ of the measured signal. However, given the vector field $\mathbf{F}$ describing the system's dynamics, the state space (in polar coordinates) is equivalent to the observable space. To distinguish the variables further in the text, we refer to $(\theta,r)$ as \textit{observable variables} and $(\varphi,\sigma)$ as \textit{reduced variables} (see Fig.~\ref{fig:expTransf}).  

We search for forward and inverse transformation functions in  Fourier-Taylor series form as given in Eq.~(\ref{eq:Kapprox}) and Eq.~(\ref{eq:Knapprox}). Namely, for the observable angle and radius we have:
\begin{equation}
	\label{eq:thappx}
	\theta =  K^{(\theta)}(\varphi,\sigma) 
     =\varphi + \sum_{n=0}^{\infty} \sum_{k=-\infty}^{\infty} q^{(\theta)}_{n,k} \sigma^n e^{\mathrm{i}k\varphi},
\end{equation}
\begin{equation}
	\label{eq:rappx}
	r  = K^{(r)}(\varphi,\sigma) 
    = \sum_{n=0}^{\infty} \sum_{k=-\infty}^{\infty} q^{(r)}_{n,k} \sigma^n e^{\mathrm{i}k\varphi}.
\end{equation}
Note that for $n=0$, i.e., on the limit cycle, we have:
\begin{eqnarray*}
 \theta &=& \varphi + \sum_{k=-\infty}^{\infty} q^{(\theta)}_{0,k} e^{\mathrm{i}k\varphi} 
\end{eqnarray*}
where the form of the phase transformation coincides with the forward transformation of the phase reduction given in Refs.~\onlinecite{Kralemann2008,Yeldesbay2019}. 

For the reduced angle and amplitude we have:
\begin{equation}
	\label{eq:phiappx}
	\varphi = \Phi(\theta,r) = \theta + \sum_{n=0}^{\infty} \sum_{k=-\infty}^{\infty} q^{(\varphi)}_{n,k} r^n e^{\mathrm{i}k\theta},
\end{equation}
and
\begin{equation}
	\label{eq:sigappx}
	\sigma = \Sigma(\theta,r) = \sum_{n=0}^{\infty}  \sum_{k=-\infty}^{\infty} q^{(\sigma)}_{n,k} r^n e^{\mathrm{i}k\theta},
\end{equation}
where $q^{(\theta)}_{n,k}$, $q^{(r)}_{n,k}$, $q^{(\varphi)}_{n,k}$, $q^{(\sigma)}_{n,k}$ are the coefficients that need to be determined. 
Here again, if we set $r=\gamma^{(r)}(\theta)$ from Eq.~\eqref{eq:gam0rth} into Eq.~\eqref{eq:phiappx} the Fourier-Taylor series of the phase transformation turns into Fourier series, which coincides with the inverse transformation of the phase reduction given in Ref.~\onlinecite{Yeldesbay2019}.

\subsubsection{Fourier-Laplace averaging integrals}
 
In section \ref{sec:Koopman} we showed that the transformation we are seeking can be obtained from the eigenfunctions of the Koopman operator.
In Ref.~\onlinecite{Mauroy2018}, it is demonstrated that an eigenfunction of the Koopman operator that corresponds to the eigenvalue $\mu$ can be found using the time averaging integral over an observable function $g(\mathbf{x})$ along the trajectory, provided that the averaging is not equal to zero everywhere and the observable function $g(\mathbf{x})$ fulfills specific properties. Let
\begin{equation}
	\label{eq:avg}
	g^*_{\mu}(\mathbf{x}_0) = \lim\limits_{\mathcal{T}\rightarrow \infty} \frac{1}{\mathcal{T}} \int_{0}^{\mathcal{T}} g(s(\tau,\mathbf{x}_0)) e^{-\mu \tau} d\tau.
\end{equation}
  If the eigenvalue $\mu$ is real, the integral \eqref{eq:avg} is called \textit{Laplace average} and the observable function should vanish on the limit cycle (i.e. $g(\mathbf{x})=0, \forall \mathbf{x}\in \Gamma$). Whereas, if the eigenvalue is imaginary (i.e. $\mu=i\omega$, $\omega \in \mathbb{R}$), then the integral is called \textit{Fourier average} and it is equal to the eigenfunction $\phi_{i \omega}$, for almost all choices of the observable function $g(\mathbf{x})$ (for details see Refs.~\onlinecite{Mauroy2012, Mauroy2018}).

The integral in Eq.~(\ref{eq:avg}) can be used to find the inverse phase-amplitude reduction transformations $\Phi(\mathbf{x})$ and $\Sigma(\mathbf{x})$. Thus, we apply the Laplace average integral to the observable function defined as $g(s(t,\mathbf{x_0})) = r(t)-\gamma^{(r)}(\theta(t))$, which vanishes on the limit cycle. 

For an uncoupled planar system, the solution of system in Eq.~(\ref{eq:phisig}) in the reduced space is 
\begin{eqnarray}
	\varphi(t) &=& \omega t + \varphi_0, \label{eq:phit}\\
	\sigma(t) &=& \sigma_0 e^{\lambda t}, \label{eq:sigt}
\end{eqnarray}
where $\varphi_0$ and $\sigma_0$ are the initial conditions. Using this solution and the approximation given in Eq.~(\ref{eq:rappx}) we can obtain the Laplace average for the initial condition $\mathbf{x}_0$ in the observable space:
\begin{multline*}
	g^*_{\lambda}(\mathbf{x}_0) = \lim\limits_{\mathcal{T}\rightarrow \infty} \frac{1}{\mathcal{T}} \int_{0}^{\mathcal{T}} \left[r(\tau)- \gamma^{(r)}(\theta(\tau))\right] e^{-\lambda \tau} d\tau \\
	= \lim\limits_{\mathcal{T}\rightarrow \infty} \frac{1}{\mathcal{T}} \int_{0}^{\mathcal{T}} \left[K^{(r)}(\varphi(\tau),\sigma(\tau))-K^ {(r)}(\varphi(\tau),0)\right] e^{-\lambda \tau} d\tau \\
	= \lim\limits_{\mathcal{T}\rightarrow \infty} \frac{1}{\mathcal{T}} \int_{0}^{\mathcal{T}} \sum_{n=1}^{\infty} \sum_{k=-\infty}^{\infty} q^{(r)}_{n,k} \sigma_0^n e^{\mathrm{i}k\varphi_0} e^{(n-1)\lambda \tau + \mathrm{i}k\omega \tau} d\tau.
\end{multline*}
 Since we have considered a hyperbolic stable limit cycle, we have $\lambda<0$. Therefore, taking the limit $\mathcal{T} \rightarrow \infty$ all terms tend to 0 except the one for $n=1$. Hence,
\begin{multline*}
	g^*_{\lambda}(\mathbf{x}_0) = \sigma_0 \sum_{k=-\infty}^{\infty} q^{(r)}_{1,k} e^{\mathrm{i}k\varphi_0} \lim\limits_{\mathcal{T}\rightarrow \infty} \frac{1}{\mathcal{T}} \int_{0}^{\mathcal{T}} e^{\mathrm{i}k\omega \tau} d\tau = \\ 
	= \sigma_0 \left(q^{(r)}_{1,0} + \lim\limits_{\mathcal{T}\rightarrow \infty} \frac{1}{\mathcal{T}}\sum_{\substack{k \in \mathbb{Z} \\ k\neq 0}}q^{(r)}_{1,k}  \frac{e^{\mathrm{i}k\varphi_0}}{\mathrm{i}k\omega} \left(e^{\mathrm{i}k\omega \mathcal{T}}-1 \right)\right) =  \sigma_0 q^{(r)}_{1,0}. 
\end{multline*}
Therefore, 
\begin{equation}
	\label{eq:avg_sigma}
	g^*_{\lambda}(\mathbf{x}_0) = q^{(r)}_{1,0}\sigma_0 = q^{(r)}_{1,0}\Sigma(\mathbf{x}_0),
\end{equation}
which indicates that the Laplace average provides the inverse transformation function $\Sigma$ scaled by a constant factor $q^{(r)}_{1,0}$.

In the same way, we can show that the angle of a Fourier average along the trajectory over an arbitrary observable function (for example, the $x$-component of the state space or the measured signal) results in the inverse transformation function for $\varphi$:
\begin{multline}
    \label{eq:avgvarphi}
	g^*_{\mathrm{i}\omega}(\mathbf{x}_0) = \lim\limits_{\mathcal{T}\rightarrow \infty} \frac{1}{\mathcal{T}} \int_{0}^{\mathcal{T}} x(\tau) e^{-\mathrm{i}\omega \tau} d\tau \\
 = \lim\limits_{\mathcal{T}\rightarrow \infty} \frac{1}{\mathcal{T}} \int_{0}^{\mathcal{T}} K^{(x)}(\varphi(\tau),\sigma(\tau)) e^{-\mathrm{i}\omega \tau} d\tau \\
	= \lim\limits_{\mathcal{T}\rightarrow \infty} \frac{1}{\mathcal{T}} \int_{0}^{\mathcal{T}} \sum_{n=0}^{\infty} \sum_{k=-\infty}^{\infty} q^{(x)}_{n,k} \sigma_0^n e^{\mathrm{i}k\varphi_0} e^{n\lambda \tau + \mathrm{i}(k-1)\omega \tau} d\tau \\
	= q^{(x)}_{0,1}e^{\mathrm{i}\varphi_0} 
    + \lim\limits_{\mathcal{T}\rightarrow \infty} \frac{1}{\mathcal{T}} \sum_{\substack{k \in \mathbb{Z} \\ k\neq 1}} q_{0,k}^{(x)} e^{\mathrm{i} k \varphi_0} \frac{e^{\mathrm{i}(k-1) \omega \mathcal{T}}-1}{\mathrm{i}(k-1) \omega}
    \\
	+ \lim\limits_{\mathcal{T}\rightarrow \infty} \frac{1}{\mathcal{T}} \sum_{n=1}^{\infty} \sum_{k \in \mathbb{Z}} q^{(x)}_{n,k}\sigma_0^n e^{\mathrm{i}k\varphi_0} \frac{e^{(n\lambda+\mathrm{i}(k-1)\omega)\mathcal{T}}-1 }{(n\lambda + \mathrm{i}(k-1)\omega)}  \\
	= q^{(x)}_{0,1} e^{\mathrm{i}\varphi_0}.
\end{multline}
Thus,
\begin{equation}
	\label{eq:avg_phi}
	\angle g^*_{\mathrm{i}\omega}(\mathbf{x}_0) = \angle q^{(x)}_{0,1} e^{\mathrm{i}\varphi_0} = \varphi_0 = \Phi(\mathbf{x}_0).
\end{equation}

\subsubsection{Numerical calculation of the inverse and forward transformations}\label{sec:transfapprox}
The numerical calculation of the Laplace average Eq.~(\ref{eq:avg}) presents certain challenges. First, we need the values of $\lambda$ and $\omega$ for the reduced system in Eq.~\eqref{eq:phisig}. Second, it is necessary to determine an optimal time range $\mathcal{T}$ for computing the time averaging integral in Eq.~(\ref{eq:avg}) (see. Ref.~\onlinecite{Mauroy2018}).  In this section, we outline the detailed steps involved in these computations.

Given the VF $\mathbf{F}(\mathbf{x})$ describing the system's dynamics, we estimate $\omega$ by determining the period of the limit cycle, and $\lambda$ by calculating the eigenvalues of the monodromy matrix obtained from Eq.~\eqref{eq:monmatrix}. Namely, we simulate the system dynamics for a sufficiently long time until the radius coordinate converges to the limit cycle $\gamma^{(r)}(\theta)$. To obtain this expression, we determine the coefficients $g_k$ of the Fourier series in Eq.~\eqref{eq:gam0rth} by fitting the time courses $\theta(t)$ and $r(t)$ on the limit cycle using the least squares error method (see Appendix~\ref{sec:rlsq}). 

Then, we compute the derivatives of $\mathbf{F}$ with respect to $r$ and $\theta$ along the limit cycle and solve the linear variational equation in Eq. (\ref{eq:monmatrix}) for one period $T$ using the initial condition $\mathbf{M}_{\varphi=0}=\mathbf{I}$. The eigenvalues of the resulting monodromy matrix $\mathbf{M}_{\varphi=2\pi}$ are equal to $\hat{\lambda}_k = e^{\lambda_k T}=e^{\lambda_k 2\pi/\omega}$. The largest negative $\lambda_k$ corresponds to the value $\lambda$ we are looking for.

\begin{figure}[ht]
	\includegraphics[width=\columnwidth]{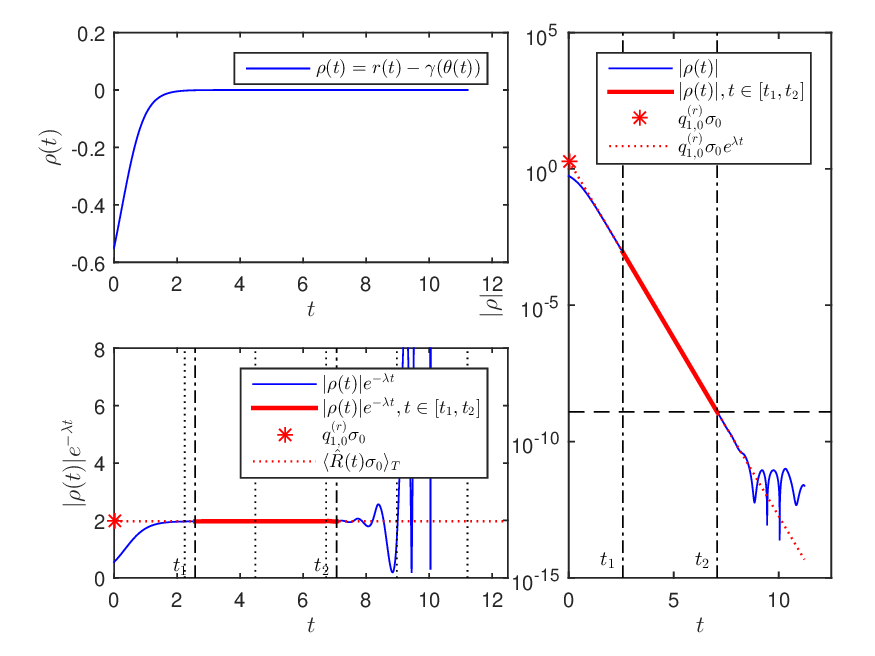}
	\caption{Illustration of how to find the amplitude variable $\sigma_0$ corresponding to an initial point $(\theta_0,r_0)$ in the observable space for the canonical model close to the numerical instabilities. \textit{Upper left panel}: Time course of the observable $\rho$ introduced in Eq.\eqref{eq:rhot} (solid blue curve) that vanishes on the limit cycle $\Gamma$. \textit{Right panel}: The absolute value of $\rho$ in semi-logarithmic scale, where the expression $\sigma_0 e^{\lambda t}$ appears as a line (dotted red line) with constant part equal to $\sigma_0$ (red star, the crossing with the vertical axes). The horizontal dashed black line is the threshold of numerical precision, and the right vertical dot-dashed black line marks the time point $t_2$ at this threshold. The left vertical dot-dashed black line at $t_1$ is in two-period-distance to the right one. The solid red line denotes the time course of $|\rho(t)|, t\in [t_1,t_2]$. \textit{Bottom left panel}: The time course of $|\rho(t)|e^{-\lambda t}$ (solid blue curve) approaches the horizontal line $\langle|\hat{R}(t)\sigma_0|\rangle_{T} \sim |q^{(r)}_{1,0} \sigma_0|$ (dotted red line, red star) before the numerical instabilities. The solid red line and the two vertical dot-dashed black lines are the same as in the right panel. The vertical dotted lines indicate the oscillation periods. \label{fig:sig0calc}}
\end{figure}

The calculation of the Laplace averaging integral is limited by numerical precision \cite{Mauroy2018}. Consider the time evolution of an observable defined as
\begin{equation} \label{eq:rhot}
 \rho(t)= r(t)- \gamma^{(r)}(\theta(t)) .
\end{equation}
 The time course of $\rho(t)$ is shown in the upper left panel of Fig.~\ref{fig:sig0calc}. Due to the limitation of the precision of the numerical integration, the product of two competing functions  - exponentially growing $e^{-\lambda t}$ and $\rho(t)$ that exponentially approaches the limit cycle - results in numerical errors in the averaging integral. This issue is illustrated by plotting $|\rho(t)|$ in logarithmic scale (blue curve in right panel of Fig.~\ref{fig:sig0calc}), where numerical instabilities appear once the numerical precision threshold is reached (indicated by the horizontal dashed black line in right panel of Fig.~\ref{fig:sig0calc}).
 
However, the initial values $\sigma_0$ and $\varphi_0$ can be determined without directly calculating the integral. Instead, they can be found by using a finite time and logarithmic scale, as it will  be discussed further. 

The observable $\rho(t)$ can be expressed as
\begin{multline*}
	\rho(t) = r(t)- \gamma^{(r)}(\theta(t))  = 
	K^{(r)}(\varphi(t),\sigma(t))-K^{(r)}_0(\varphi(t)) \\
	= \sum_{n=1}^{\infty} \sum_{k=-\infty}^{\infty} q^{(r)}_{n,k} \sigma_0^n e^{\mathrm{i}k\varphi_0} e^{n\lambda t + \mathrm{i}k\omega t}.
\end{multline*}  
If we factor out $\sigma_0 e^{\lambda t}$, we have
\begin{multline}
    \label{eq:rhot_short}
	\rho(t) = \sigma_0 e^{\lambda t}\sum_{n=1}^{\infty} \sum_{k=-\infty}^{\infty} q^{(r)}_{n,k} \sigma_0^{n-1} e^{\mathrm{i}k\varphi_0} e^{(n-1)\lambda t + \mathrm{i}k\omega t} \\
		 \triangleq \sigma_0 e^{\lambda t} R(t).
\end{multline}
Since $\lambda<0$, after a transient time $\mathcal{T}$, the terms involving $e^{\lambda t}$ in the expression for $R(t)$  decay down to zero and, therefore, only oscillatory terms with $n=1$ remain, i.e.
\begin{equation}\label{eq:limRt}
 R(t) \sim \sum_{k=-\infty}^{\infty} q^{(r)}_{1,k} e^{\mathrm{i}k\varphi_0} e^{\mathrm{i}k\omega t} \triangleq \hat{R}(t).  
\end{equation}
The average of $\hat R(t)$ over one period $T$ is 
\[
\langle \hat{R}(t) \rangle_{T}= \left\langle\sum_{k=-\infty}^{\infty} q^{(r)}_{1,k} e^{\mathrm{i}k\varphi_0} e^{\mathrm{i}k\omega t} \right\rangle_{T} = q^{(r)}_{1,0}.  
\]
Hence, from Eqs.~\eqref{eq:rhot_short} and \eqref{eq:limRt} we have
\begin{equation}
     \langle |\rho(t)| e^{- \lambda t} \rangle_{T}= \langle|\sigma_0  R(t)|\rangle_{T} \sim |q^{(r)}_{1,0} \sigma_0|.
\end{equation}
On the other hand, we have
\begin{equation}
\log |\rho(t)| \sim \lambda t + \log (\sigma_0 \hat{R}(t)), \label{eq:avgsimRt}
\end{equation}
that is, the asymptotic behavior of the function can be described by a straight line with slope $\lambda$ plus an oscillatory term around $\log |q^{(r)}_{1,0} \sigma_0|$. This is shown in the right panel of Fig.~\ref{fig:sig0calc}, where $|\rho(t)|$ is plotted in semi-logarithmic scale. As shown, the calculation of $|\rho(t)|$ becomes unstable after approximately the time point $t_2$, which is defined using the precision threshold of the numerical integration (horizontal dashed line). If we consider an interval of several periods before reaching this threshold ($t\in[t_1,t_2]$, where $t_2=t_1+2T$), we observe that the time course of $|\rho(t)|$ in logarithmic scale oscillates around the function $q^{(r)}_{1,0} \sigma_0 e^{\lambda t}$, which appears as a straight line in the logarithmic scale plot (red dotted line).

Therefore, to find $|q^{(r)}_{1,0} \sigma_0|$,  we average the last few periods of $|\rho(t)| e^{- \lambda t}$ before reaching the numerical precision threshold, as shown in the bottom panel of Fig.~\ref{fig:sig0calc}.  In this case, an interval of two periods (red line) was used. Notice that the time course of $|\rho(t)| e^{- \lambda t}$ becomes generally oscillatory after a transient time (see Eq.~(\ref{eq:limRt})), as shown in Fig.~\ref{fig:sig0calc_vdp} in the Appendix. Moreover, if the initial estimate of $\lambda$ obtained from solving Eq.~\eqref{eq:monmatrix} is not accurate, one can determine the correction to $\lambda$ by using the average slope of the curve $\log |\rho(t)|- \lambda t$, which is also shown in Fig.~\ref{fig:sig0calc_vdp} in the Appendix.

As discussed previously, the reduced amplitude $\sigma$ can be defined uniquely up to a constant multiplier. Therefore, we can redefine the amplitude $\sigma_0$ as $q^{(r)}_{1,0}\sigma_0$, thereby fixing the scaling of $\sigma$. The sign of $\sigma_0$ is determined by the sign of $\rho(0)$.  

In a similar way, we can determine the initial value $\varphi_0$ by computing the angle of the average of $x(t)e^{-\mathrm{i}\omega t}=r(t)\cos(\theta(t))e^{-\mathrm{i}\omega t}$ over several periods of oscillation, after a transient time and before reaching the numerical precision threshold. This process is equivalent to calculating  the integral in Eq.~\eqref{eq:avgvarphi}, where the terms involving $e^{n\lambda t}$ decay to 0 after some transient time and the terms involving $e^{\mathrm{i} k \omega T}$ vanish due to time averaging over several periods. We define $\varphi_0$ as $\varphi_0=0$ when $\theta_0=0$. 

Having computed several pairs of initial conditions in the observable space $\{\theta_{j,0},r_{j,0}\}_{j \in \mathbb{N}}$, and the corresponding ones in the reduced space $\{\varphi_{j,0},\sigma_{j,0}\}_{j \in \mathbb{N}}$, we can then approximate the forward transformation functions in the form given in Eqs.~(\ref{eq:thappx})-(\ref{eq:rappx}) using the regularized least squares method (see Appendix~\ref{sec:rlsq}) with the Fourier-Taylor basis functions
\begin{multline}
	\label{eq:psichi}
	\Psi^{(\chi)} = \left\{\sigma^{n} e^{\mathrm{i}k\varphi}, n = 0,\ldots,N_n, k = -N_k,\ldots,N_k\right\},
\end{multline}
where the superscript $(\chi)$ denotes the dependence on the reduced variables $\varphi$ and $\sigma$.

Similarly, we can approximate the inverse transformations in the form given in Eqs.~\eqref{eq:phiappx} and \eqref{eq:sigappx} 
using the basis functions
 \begin{equation}
 	\label{eq:psix}
 	\Psi^{(x)} = \left\{r^{n} e^{\mathrm{i}k\theta}, n = 0,\ldots,N_n, k = -N_k,\ldots,N_k\right\},
 \end{equation}
 where the superscript $(x)$ denotes the dependence on the observable variables $\theta$ and $r$. 
    
Finally, the complete procedure for finding the transformation functions is as follows: 
\begin{enumerate}
	\item Knowing the VF of the system, or a good approximation of it in a given region, numerically integrate the system to find the limit cycle $\Gamma$, the frequency $\omega$, and the approximation of the limit cycle as given by Eq.~\eqref{eq:gam0rth} using the least squares (LSQ) method (see Appendix~\ref{sec:rlsq}). 
	\item Using the limit cycle, solve the variational equations in Eq. (\ref{eq:monmatrix}) and estimate $\lambda$ by finding the largest (in modulus) eigenvalue $e^{\lambda}$ of the monodromy matrix $\mathbf{M}_{\varphi=2\pi}$ . 
	\item Numerically integrate the system with different initial conditions on a grid $\{\theta_{j,0},r_{j,0}\}_{j \in \mathbb{N}}$ over several periods, obtaining the time courses  $r(t)$ and $\theta(t)$ of the solutions. 
    \item Compute $\rho(t)$ and refine the value of $\lambda$ by determining the average slope of $\log |\rho(t)| - \lambda t$ during the last two periods of the simulation (or two periods before reaching the numerical precision threshold). See Appendix~\ref{sec:sigma0calc} for details.
	\item Determine $|\sigma_0|$ for a given $(\theta_0,r_0)$ by averaging $|\rho(t)|e^{- \lambda t}$ over the last two periods of the simulation (or two periods before reaching the numerical precision threshold) and setting $\mathrm{sign}(\sigma_0) =\textrm{sign}(\rho(0))$.
	\item Determine $\varphi_0$ for a given $(\theta_0,r_0)$  as the angle of the time average of $r(t)\cos(\theta(t))$ over the last two periods of the simulation (or two periods before reaching the numerical precision threshold).
	\item Use the initial conditions in the observable space $\{\theta_{j,0},r_{j,0}\}_{j \in \mathbb{N}}$ and the reduced space $\{\varphi_{j,0},\sigma_{j,0}\}_{j \in \mathbb{N}}$ to approximate the forward transformations given in Eqs.~(\ref{eq:rappx}) and (\ref{eq:thappx}) by means of the regularized least squares method (see Appendix~\ref{sec:rlsq}) with the basis functions in Eq.~(\ref{eq:psichi}).  
	\item Similarly, approximate the inverse transformations in Eqs.~(\ref{eq:sigappx}) and (\ref{eq:phiappx}) by means of the regularized least squares method (see Appendix~\ref{sec:rlsq}) with the basis functions in Eq.~(\ref{eq:psix}).  
\end{enumerate}

\subsection{Amplitude and phase reduction for coupled oscillatory systems}  

In the previous section, we discussed how to determine the transformation functions assuming that the VF of the system, or a good approximation of it, is known. In this section, we present a method for finding the VF and performing the phase-amplitude reduction for coupled systems using the signals measured from these systems.

\subsubsection{Approximation of the VF for coupled oscillatory systems}
Let us assume that we have a set of $M$ trials of measured signals from $N$ different sources with high time resolution, 
$\{s_i(t)\}_p, \, i=1,\ldots,N, \, p=1,\ldots M$. Then, using methods to extract the amplitude and angle of the signals, for example the Hilbert transformation \cite{Rosenblum2021, Busch2022}, we obtain a set of trials with radius and angle time courses $\{\theta_i(t),r_i(t)\}_p$, which we refer to as observable amplitude and phase. 
Using finite differences we obtain the time derivatives of the observable radius and phase $\{\dot{\theta}_i(t),\dot{r}_i(t)\}_p$. 

The dynamics of a source $i$ \footnote{A brain region can have several sources of oscillatory activity around different frequencies.} can be represented as
\begin{align}
	\dot{x}_i = \mathbf{F}^{(x)}_i(\theta_i,r_i,\ldots, v_j(\theta_j,r_j), \ldots), 
\end{align}
where $x$ is either $\theta$ or $r$, and $v_j$ is the input from other sources coupled to source $i$, e.g. $v_j=r_j\cos(\theta_j)$. As we have shown in the previous sections, in order to find the transformation functions we need an expression for the VF of the uncoupled system. The function $\mathbf{F}_i$ is, in general, non-linear with respect to its variables, thus the dynamics of the uncoupled system can not be separated from the dynamics induced by the coupling. However, if we approximate $\mathbf{F}_i$ using a Fourier-Taylor series, then all the terms in the series that depend only on $\theta_i$ and $r_i$ correspond to the uncoupled dynamics, and the rest to the coupling. Moreover, we assume that inputs from other sources do not interact with each other and we can approximate the coupling part by only considering pairwise interactions. 

Thus, the approximated VF describing the dynamics of source $i$ consists of its uncoupled dynamics and the sum of pairwise couplings with other sources, that is, 
\begin{equation}
	\label{eq:rhsappx}
	\dot{x}_i \approx F^{(x)}_i(\theta_i,r_i) + \sum_j G^{(x)}_{i,j}(\theta_i,r_i,\theta_j,r_j), 
\end{equation}
where $x$ is again either $\theta$ or $r$. In Eq.~(\ref{eq:rhsappx}) we refer to the function $F^{(x)}_i(\theta_i,r_i)$ representing the uncoupled dynamics as \textit{uncoupled VF}, which we approximate as
\begin{equation}
	\label{eq:Fi}
	F^{(x)}_i(\theta_i,r_i) = \sum_{n=0}^{N_n} \sum_{k=-N_k}^{N_k} q^{(x_i)}_{n,k} r_i^n e^{\mathrm{i}k\theta_i},
\end{equation}
where $q^{(x_i)}_{n,k}$ are the coefficients of the Fourier-Taylor series for the VF in source $i$. 

The pairwise coupling function, $G^{(x)}_{i,j}(\theta_i,r_i,\theta_j,r_j)$, referred to as \textit{coupling part of the VF}, depends also on the variables $\theta_j$ and $r_j$ of other sources. We approximate this function using a Fourier-Taylor series expansion in the variables $(\theta_i,r_i,\theta_j,r_j)$ as follows:
\begin{multline}
	\label{eq:Gij}
	G^{(x)}_{i,j}(\theta_i,r_i,\theta_j,r_j) = 
	\sum_{m=0}^{N_m} \sum_{l=0}^{m-l} \sum_{k_i=-N_{k_i}}^{N_{k_i}}  \\ \sum_{\substack{k_j=-N_{k_j} \\ k_j\neq 0}}^{N_{k_j}} c^{(x_{i,j})}_{m,l,k_i,k_j} r_i^l r_j^{m-l} e^{\mathrm{i}(k_i\theta_i+k_j\theta_j)},
\end{multline}
where $c^{(x_{i,j})}_{m,l,k_i,k_j}$ are the coefficients which need to be determined, $x$ is either $\theta$ or $r$, and the subscripts $i,j$ of the superscripts $(x_{i,j})$ indicate the relation to the coupling $G^{(x)}_{i,j}$.
Notice that the sum in Eq.~(\ref{eq:Gij}) excludes terms with $k_j=0$. Terms with $r_j^0$ and $k_j=0$ correspond to the uncoupled VF (Eq.~(\ref{eq:Fi})), since only terms involving $r_i$ and $\theta_i$ remain. Terms with $r_j^m$ (for $m\neq0$) and $k_j=0$ are excluded to ensure the dependence on $v_j$, as the variables $r_j$ and $\theta_j$ do not appear separately in the coupling part $G^{(x)}_{i,j}$, but only in combination via $v_j$.

In practice, the approximations in Eq.~(\ref{eq:rhsappx})-\eqref{eq:Fi}-\eqref{eq:Gij} are obtained by means of the least squares method (see Appendix~\ref{sec:rlsq}) using a combination of the following basis functions
\begin{equation}
	\Psi^{(x)}_{i} = \{r_i^n e^{\mathrm{i}k\theta_i}, n= 0,\ldots,N_n \text{ and } k=-N_k,\ldots,N_k \},\label{eq:bfunc_rth_i}
\end{equation}
and 
\begin{align}
	\Psi^{(x)}_{ij} = \left\{r_i^{m_{i}} r_j^{m_{j}} e^{\mathrm{i}(k_i\theta_i+k_j\theta_j)}\right., \label{eq:bfunc_rth_ij} \\
	m_i,m_j = 0,\ldots,N_m, \; m_i+m_j \leq N_m, \nonumber\\
	k_i=-N_{k_i},\ldots,N_{k_i}, \nonumber \\ 
    \left. k_j =-N_{k_j},\ldots,N_{k_j}, \; k_j\neq 0 \nonumber\right\},
\end{align}
for every input region $j$. The superscript $(x)$ in Eqs.~\eqref{eq:bfunc_rth_i} and \eqref{eq:bfunc_rth_ij} denotes the dependence on the state variables $\theta$ and $r$, and the lower indices $i$ or $j$ denote of which source. Notice that we have different ranges for the indices $n$ and $k$ in Eq.~\eqref{eq:bfunc_rth_i} and for the indices $m_i$, $m_j$, $k_i$ and $k_j$ in Eq.~\eqref{eq:bfunc_rth_ij}. 

Using the above equations, we can compute an approximation of the uncoupled VF $F^{(x)}_{i}$ and the coupling part of the VF $G^{(x)}_{i,j}$ for every source and the incoming input to this source. The obtained functions give a good approximation of the uncoupled VF and coupling in the domain that is covered by the measured data. Furthermore, using the regularized least squares method (see Appendix~\ref{sec:rlsq}) we escape over-fitting and provide a good approximation even in the surrounding of the domain covered by the measured data, as it will be shown in the results Section. 

\subsubsection{Computing the coupling in the reduced space}
   
The equations of the dynamics of a system with coupling in the reduced space read as:
\begin{equation}
	\label{eq:phisigcoupl}
	\left\{
	\begin{array}{cc}
		\dot{\varphi}_i = \omega_i + \sum_j g^{(\varphi)}_{i,j}(\varphi_i,\sigma_i,\varphi_j,\sigma_j), \\
		\dot{\sigma}_i = \lambda_i \sigma_i +\sum_j g^{(\sigma)}_{i,j}(\varphi_i,\sigma_i,\varphi_j,\sigma_j),
	\end{array}	 
	\right.	
\end{equation}
where $g^{(\varphi)}_{i,j}$ and $g^{(\sigma)}_{i,j}$ are the pairwise coupling functions between the oscillator $i$ and $j$. 

Using the invariance equations Eqs.~\eqref{eq:invPhiPol}-\eqref{eq:invSigPol}, taking into account the approximation of the VF in Eqs.~(\ref{eq:rhsappx})-\eqref{eq:Fi}-\eqref{eq:Gij}, and 
assuming independence of the coupling terms between different oscillators, we obtain the following equations for $g^{(\varphi)}_{i,j}$ and $g^{(\sigma)}_{i,j}$ in the reduced space separately for every pairwise coupling $\{i,j\}$:
\begin{eqnarray}
	\label{eq:gvarphi}
	g^{(\varphi)}_{i,j} \circ (\Phi_i,\Sigma_i,\Phi_j,\Sigma_j) = \frac{\partial \Phi_i}{\partial \theta_i} G^{(\theta)}_{i,j} + \frac{\partial \Phi_i}{\partial r_i} G^{(r)}_{i,j}, \\
	\label{eq:gsigma}
	g^{(\sigma)}_{i,j} \circ (\Phi_i,\Sigma_i,\Phi_j,\Sigma_j)= \frac{\partial \Sigma_i}{\partial \theta_i} G^{(\theta)}_{i,j} + \frac{\partial \Sigma_i}{\partial r_i}  G^{(r)}_{i,j}, 
\end{eqnarray}
where $\Phi_i = \Phi_i(\theta_i,r_i)$ and $\Sigma_i = \Sigma_i(\theta_i,r_i)$ are the inverse transformation functions approximated using the method presented in Section \ref{sec:transfapprox}. Notice that the terms on the right hand side of Eqs.~(\ref{eq:gvarphi}) and (\ref{eq:gsigma}) depend on the observable variables $(\theta,r)$. In the practical implementation, we define a  four dimensional grid $\{\theta_i,r_i,\theta_j,r_j\}$ for every $j=1,\ldots,N$ and find the values of $g^{(\varphi)}_{i,j}$ and $g^{(\sigma)}_{i,j}$ with respect to this grid in the observable space. Then, using the inverse transformations we find the corresponding four dimensional grid in the reduced space $\{\varphi_i,\sigma_i,\varphi_j,\sigma_j\}$ and approximate the coupling using Fourier-Taylor series as follows
\begin{multline}
g^{(\chi)}_{i,j}(\varphi_i,\sigma_i,\varphi_j,\sigma_j) \approx \sum_{n_i=0}^{N_{n_i}}\sum_{n_j=0}^{N_{n_j}} \sigma_i^{n_i} \sigma_j^{n_j} \\ \sum_{k_i=-N_{k_i}}^{N_{k_i}} \sum_{\substack{k_j=-N_{k_j} \\ k_j\neq 0}}^{N_{k_j}} p^{(\chi_{i,j})}_{n_i,n_j,k_i,k_j} e^{\mathrm{i}(k_i \varphi_i + k_j \varphi_j)}, \label{eq:coupl_appx_red}
\end{multline}
where $\chi$ is either $\varphi$ or $\sigma$, and again the subscripts $i,j$ of the superscripts $(\chi_{i,j})$ denote the reference to the coupling $g^{(\chi)}_{i,j}$. 

Then, the coefficients $p^{(\chi_{i,j})}_{n_i,n_j,k_i,k_j}$ together with $\lambda_i$ and $\omega_i$ provide the whole set of parameters of the amplitude and phase equations Eq.~(\ref{eq:phisigcoupl}) in the reduced space. 

The approximation in Eq.~(\ref{eq:coupl_appx_red}) can be written in a more suitable form for the regularized least squares method (see Appendix~\ref{sec:rlsq}) using the basis functions
\begin{align}
	\Psi^{(\chi)}_{{i,j}} = \left\{\sigma_i^{n_i} \sigma_j^{n_j} e^{\mathrm{i}(k_i \varphi_i + k_j \varphi_j)}\right., \label{eq:bfunc_coupl_red} \\
	n_i,n_j = 0,\ldots,N_m, k_i=-N_{k_i},\ldots,N_{k_i}, \nonumber \\
	\left. k_j=-N_{k_j},\ldots,N_{k_j}, \; k_j\neq 0 \right\}. \nonumber
\end{align}


    
\section{Results}

\subsection{Description of the models }

To test our method we used four models: (i) the radial isochron clock, (ii) the canonical model for an oscillator, also known as Stuart-Landau oscillator, (iii) the van der Pol oscillator and (iv) the Wilson-Cowan model (see Appendix~\ref{sec:models}). The first two models have circular limit cycles and known transformation functions (see Eqs.~(\ref{eq:transf_inv_radisocl}) and (\ref{eq:transf_inv_canmod}), and Supplementary Materials for their derivations), whereas the last two models have non-circular limit cycles and there are no known analytical expressions for their transformation functions. Moreover, the phase variable of the uncoupled radial isochron clock does not depend on the amplitude (see Eq.~(\ref{eq:radisocl})) and the input in the Wilson-Cowan model is non-linear (see Eq.~(\ref{eq:wcmod})). The properties of the models are summarized in Table 1.
\begin{table*}[ht!]
	\label{tab:models}
	\centering
    \begin{tabular}{|l|c|c|c|c|}
		\hline
		Model name & analytical transformation & phase depends on amplitude  & circular limit cycle & linear input \\
		\hline 
		Radial isochron clock & known & no & yes & yes \\
		Canonical model & known & yes & yes & yes \\
		van der Pol oscillator & unknown & yes & no & yes \\
		Wilson-Cowan & unknown & yes & no & no \\
		\hline
	\end{tabular}
	\caption{Table of the properties of the models studied.}
\end{table*} 
 
 We simulated these models as pairs of two coupled oscillators of the same kind but with different parameter values (see Appendix~\ref{sec:models} for the parameter values). We generated 100 trials for every pair of oscillators with different initial conditions. The initial conditions for the simulations were chosen randomly to cover the region around the limit cycles. The trials lasted several periods of the oscillations. In all simulations the first oscillator was driving the second one.
    
\subsection{Reconstruction of the VF}    
In this Section, we present results of the reconstruction of the VF using synthetically simulated signals (trials) from the different models presented above.

From the trials, we first determine the time courses of the radius $r(t)$ and the angle $\theta(t)$. Then, we compute a numerical approximation of their time derivatives $\dot{r}(t)$ and $\dot{\theta}(t)$. Notice that, in practice, experimental signals are typically filtered around a specific frequency before analyisis. As a result, the time courses of both the radius and angle, as well as their derivatives, will exhibit smooth behavior.

Using the approximated derivatives and the time courses of the variables, we approximate the VF with the regularized least squares method (see Appendix~\ref{sec:rlsq}) using the basis functions in Eq.~(\ref{eq:bfunc_rth_i}) and Eq.~(\ref{eq:bfunc_rth_ij}). This allows us to represent the VF as a sum of the uncoupled VF and the pairwise coupling from other sources.

\begin{figure}[h!]
	\centering
	\includegraphics[width=0.48\textwidth]{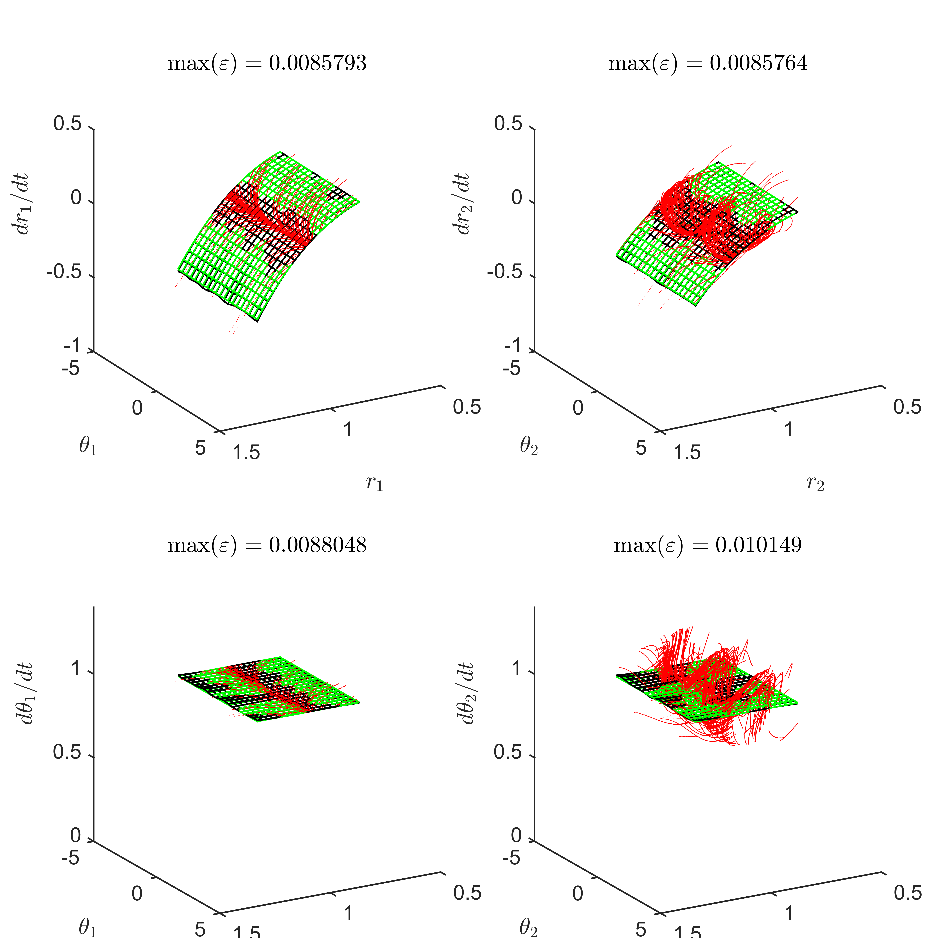}
	\caption{An example of the uncoupled VF reconstruction for two non-identically coupled radial isochron clocks. The left panels correspond to the first, and the right ones to the second oscillator. The first oscillator drives the second one. Red dots are the data points (trials). Approximated uncoupled VFs are depicted as black mesh surfaces. The corresponding theoretical VFs are plotted as green mesh surfaces. The numbers above each panel correspond to the maximal deviation of the green and black surfaces at the data points. \label{fig:rhs_radisocl}}
\end{figure}

\begin{figure}[h!]
	\centering
	\includegraphics[width=0.48\textwidth]{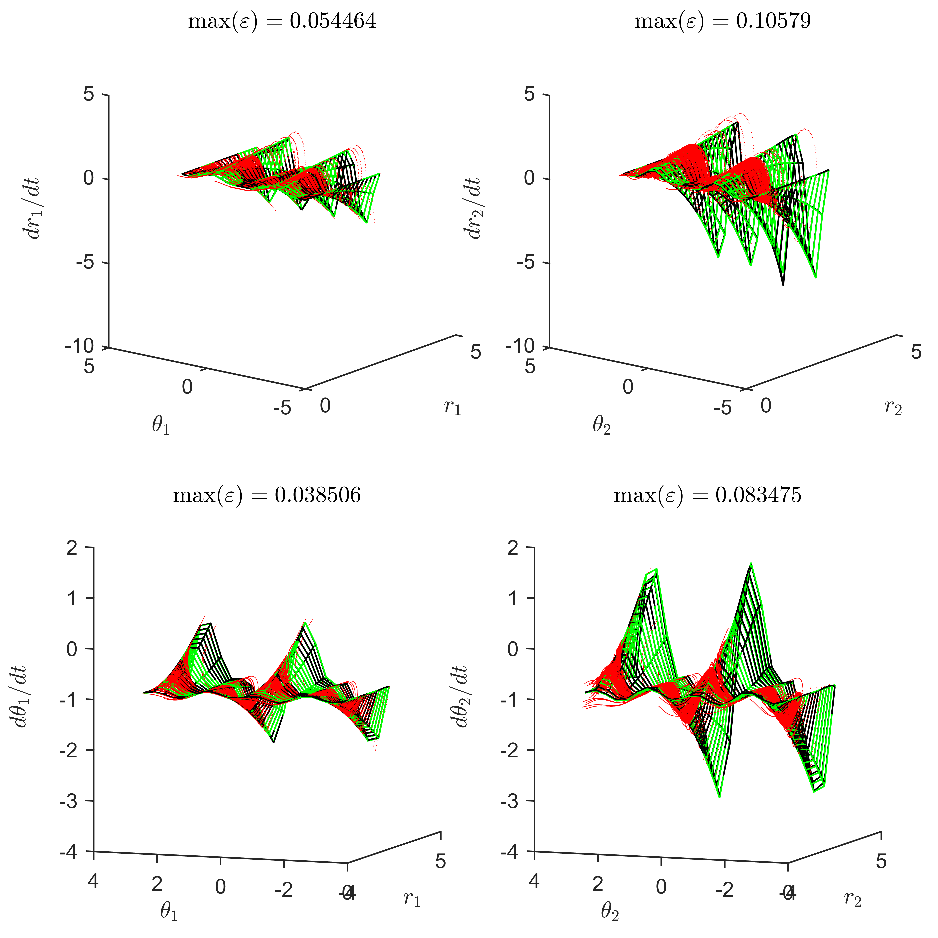}
	\caption{An example of the uncoupled VF reconstruction for two non-identically coupled van der Pol oscillators. The notations are the same as in Fig.~\ref{fig:rhs_radisocl}. \label{fig:rhs_vdp}}
\end{figure}

The results of the reconstruction of the uncoupled VF for all four different models are given in the Supplementary Materials. Examples of the reconstruction of the VF for two coupled radial isochron clocks and two coupled van der Pol oscillators are given in Figs.~\ref{fig:rhs_radisocl} and \ref{fig:rhs_vdp}, respectively. In both figures, as goodness of approximation, the maximal deviation of the reconstructed uncoupled VF from the theoretical one among all coordinates of the data points is given on the top of every panel. Overall, one can see that the approximations are close to the theoretical VFs. 

Observe in Fig.~\ref{fig:rhs_radisocl} that the trials (data points in red) for the uncoupled oscillator (left panels) gather around the surface of the uncoupled VF, whereas for the coupled oscillator (right panels), the trials strongly deviate due to the input received from the first oscillator. Nevertheless, the VF approximation successfully distinguishes the uncoupled VF from the coupling with high accuracy, as indicated by the maximal deviation between the theoretical (green mesh) and approximated (black mesh) surfaces at the coordinates of data points, which is quantified above each panel. 

Since the trials do not uniformly cover the space $(\theta,r)$, some regions on the plane $(\theta,r)$ have a higher density of data points than others. This uneven distribution causes slight deviations between the values of the approximated VF and the theoretical one in areas with lower density or absence of data points. Notably, thanks to the regularization (ridge regression) used in the least squares method, the approximation does not deviate much from the theoretical one, even outside the region covered by the data points, with only minor  deviations in these uncovered areas, which still maintain a close fit to the theoretical model.

Even for a highly nonlinear VF function, such as for the van der Pol oscillator (see Fig.~\ref{fig:rhs_vdp}), the VF can be well approximated outside of the regions covered by the data points. 
    
\subsection{Reconstruction of the transformations}    
After determining the uncoupled VF of the oscillatory systems, we apply the approximation method presented in Section \ref{sec:transfapprox} to find the individual inverse and forward transformations of the oscillators. Specifically, we construct a grid of initial conditions in the observable space. For each initial condition $\{\theta_{0,k},r_{0,k}\}_{k\in \mathbb{N}}$ in the observable space, we simulate the system over several oscillation periods and determine the corresponding initial conditions in the reduced space $\{\varphi_{0,k},\sigma_{0,k}\}_{k \in \mathbb{N}}$ using the linear approximation method discussed in the Methods Section. With these pairs of initial conditions in both the observable and reduced spaces, the transformation functions are then approximated using the regularized least squares method (see Appendix~\ref{sec:rlsq}).
 
In Fig.~\ref{fig:theor_vs_approx_transf} we show the comparison between the reconstructed inverse transformations and the theoretical ones for the radial isochron clock and the canonical model (see Appendix~\ref{sec:radisocl} and \ref{sec:canmod}). Since for both models the theoretical expression of the inverse transformation $\sigma=\Sigma(\theta,r)$ does not depend on the phase $\theta$ (see Eqs.~(\ref{eq:transf_inv_radisocl}) and (\ref{eq:transf_inv_canmod})), the comparison is performed only for the observable radius $r$ and the reduced amplitude $\sigma$. The blue dots correspond to the approximated initial values of the amplitude in the reduced space $\sigma_0$. 

\begin{figure}
    \includegraphics[width=0.49\columnwidth]{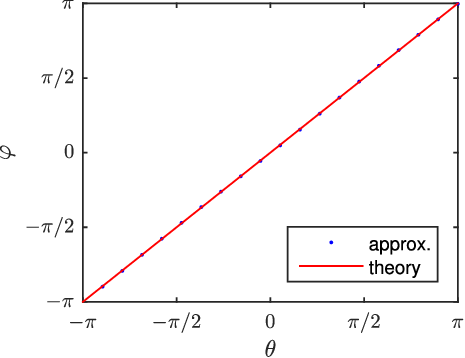}
    \includegraphics[width=0.49\columnwidth]{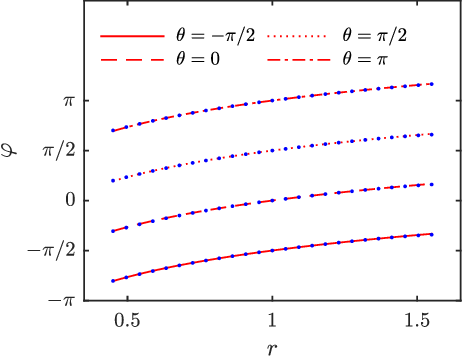}
    \caption{Comparison of the reconstructed inverse phase transformations and the theoretical ones for the radial isochron clock (left panel) and the canonical model (right panel). \textit{Left panel}: Blue dots are the approximated values of the inverse phase transformation (the initial values of the phase in the reduced space $\varphi$) for the radial isochron clock with respect to different initial phases $\theta$ and with a fixed initial value of the radius $r=1.2$ in the observable space. Red curve represents the corresponding theoretical values of the inverse phase transformation for the radial isochron clock. \textit{Right panel}: Blue dots are the approximated values of the inverse phase transformations (the initial values of the phase in the reduced space $\varphi$) for the canonical model with respect to the initial values of the radius $r$ for different fixed values of the initial phase $\theta$. Red lines represent the corresponding theoretical values of the inverse phase transformation for different fixed values of $\theta$. \label{fig:theor_vs_approx_phase_transf}}
    \end{figure}

In Fig.~\ref{fig:theor_vs_approx_phase_transf} we show the comparison between the approximated inverse phase transformation and the theoretical one for the radial isochron clock (left panel) and the canonical model (right panel). For the radial isochron clock, the phase transformation is $\varphi=\theta$ (Eq.~\eqref{eq:transf_inv_radisocl}) and does not depend on the radius (see the left panel in Fig.~\ref{fig:theor_vs_approx_phase_transf}). For the canonical model, the inverse transformation depends logarithmically on the radius (Eq.~\eqref{eq:transf_inv_canmod}). In the right panel of Fig.~\ref{fig:theor_vs_approx_phase_transf} the approximated initial values of $\varphi$ are compared with the theoretical ones for different initial values of radius $r$ and different fixed values of $\theta$. As one can see, the approximations of the initial amplitude and phase, and therefore the approximated amplitude and phase inverse transformation functions, coincide with the theoretical values for the radial isochron clock and the canonical model.

\begin{figure*}[ht!]
	\centering
	\includegraphics[width=0.45\textwidth]{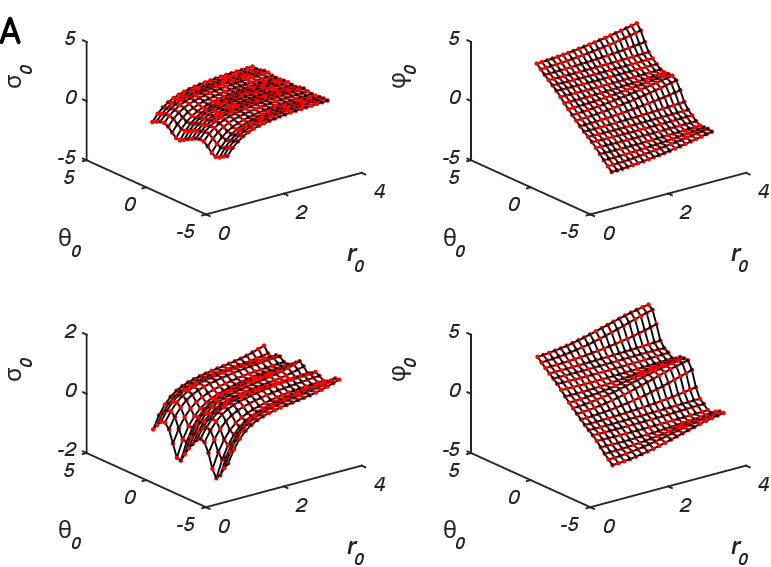} 
	\includegraphics[width=0.45\textwidth]{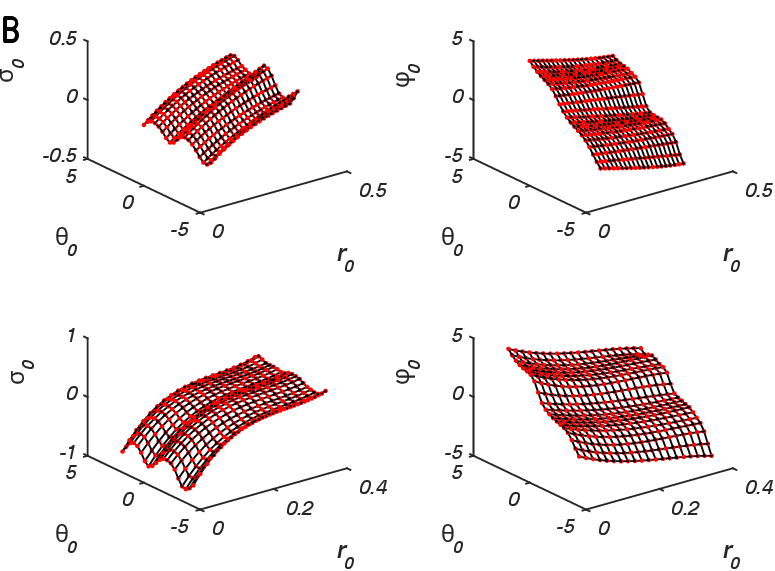}
	\caption{An example of the reconstructed inverse transformation functions $\Sigma(\theta,r)$ and $\Phi(\theta,r)$ for (A) the van der Pol oscillators and (B) the Wilson-Cowan models. The left column of each panel corresponds to the amplitude transformation $\sigma_0=\Sigma(\theta_0,r_0)$ and the right column corresponds to the phase transformation $\varphi_0=\Phi(\theta_0,r_0)$ from the observable to the reduced space. The upper panels correspond to the first oscillator and the lower panels to the second oscillator. The red dots are the initial values of the amplitude and phase in the observable $(\theta_0,r_0)$ and reduced ($\varphi_0$, $\sigma_0$) spaces. The black surfaces are the approximated inverse transformations using the regularized least squares method. \label{fig:transf_wcmod}}
	\end{figure*}

In Fig.~\ref{fig:transf_wcmod}, we show an example of the reconstructed inverse transformations for the van der Pol oscillator (see Appendix~\ref{sec:vdp}) and the Wilson-Cowan model (see Appendix~\ref{sec:wcmod}). In both models, the inverse transformations depend non-linearly on both amplitude and phase. The results of the reconstruction of the transformation functions for all oscillator models, together with the derivation of the theoretical expressions for the transformation functions are given in the Supplementary Materials.
	
	\begin{figure}[h]
		\centering
		\includegraphics[width=0.49\columnwidth]{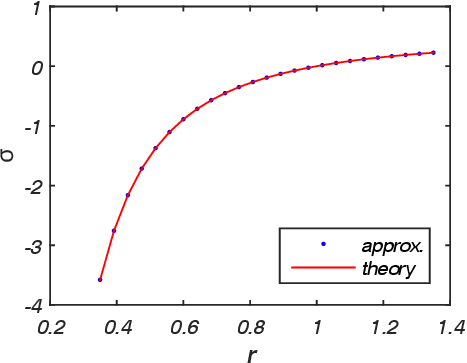}
		\includegraphics[width=0.49\columnwidth]{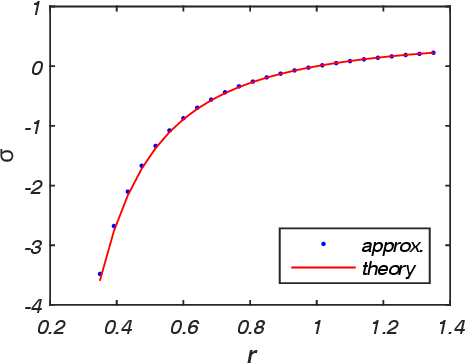}
		\caption{Comparison of the reconstructed inverse amplitude transformations and theoretical ones for the radial isochron clock (left panel) and the canonical model (right panel). Blue dots are the approximated initial values of the amplitude in the reduced space $\sigma$ for the simulated trials with the initial values of the radius $r$ in the observable space. Red curves are the corresponding theoretical values of the inverse amplitude transformations.\label{fig:theor_vs_approx_transf}}
	\end{figure}

\subsection{Reconstruction of the coupling part}    
In this Section, we present the results of reconstructing the coupling part of the VF for the four models under consideration. As indicated in the Methods Section, we first determine the coupling in the observable space using the approximated coupling component of the VF and the gradients of the inverse transformations. To do this, we generate random points in the four dimensional space $(\theta_i,r_i,\theta_j,r_j)$ for every pair of coupling $ij$. The range of random points is selected between 0.8 times the maximum and 1/0.8 times the minimum radius of the data points. This ensures a more accurate approximation, since both the VF and the transformation functions are better approximated near the limit cycle. 
    
After estimating the approximation of the coupling component $G^{(x)}_{i,j}$ in the observable space according to Eq.~(\ref{eq:Gij}) with the basis functions in  Eq.~(\ref{eq:bfunc_rth_ij}), we determine the corresponding coordinates of the random points in the reduced space $(\varphi_i,\sigma_i,\varphi_j,\sigma_j)$ using the inverse transformation. Then, we approximate the coupling component $g^{(\chi)}_{i,j}$ in the reduced space according to Eq.~(\ref{eq:coupl_appx_red}) with the basis functions in Eq.~(\ref{eq:bfunc_coupl_red}). The results of the approximation are the real-valued Fourier-Taylor coefficients $p^{(\chi_{i,j})}_{n_i,n_j,k_i,k_j}$ that define the coupling in the amplitude and phase variables (see Eqs.~(\ref{eq:phisigcoupl})) for every pair of oscillators $ij$. 
    
\begin{figure}
    	\centering
    	\includegraphics[width=0.45\textwidth]{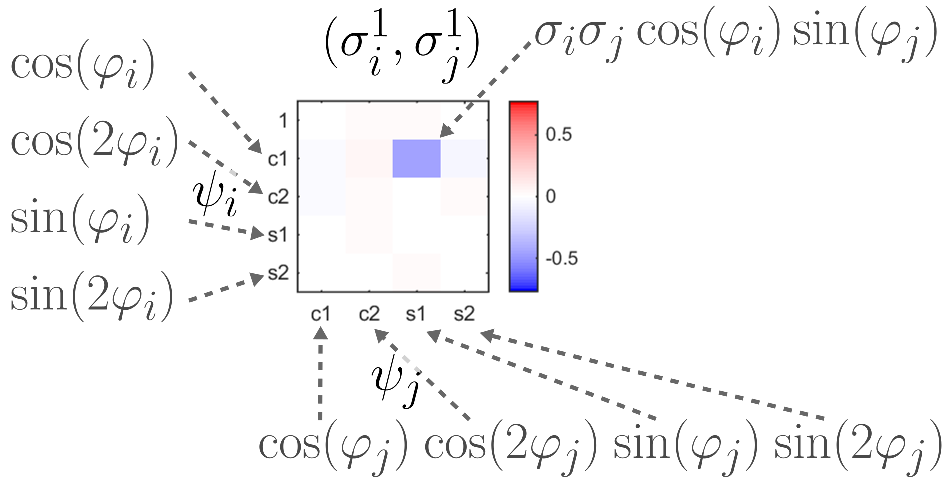}
    	\caption{Explanation of the meaning of the Fourier-Taylor coefficients panels for a particular example.  See text for more details. \label{fig:expcoefs}}
    \end{figure}

We present the approximated Fourier-Taylor coefficients $p^{(\chi_{i,j})}_{n_i,n_j,k_i,k_j}$ as a heat map (see Fig.~\ref{fig:expcoefs}). The coefficients are sorted in panels based on the amplitude terms of the corresponding basis functions, that is, the combinations of powers of $\sigma_i$ and $\sigma_j$ are indicated by labels at the top of each panel. In all panels, the vertical axis corresponds to the components of the basis functions that depend on the oscillator's own amplitude and phase (i.e. $\Psi_i$), while the horizontal axis corresponds to those that depend on the amplitude and phase of an input oscillator (i.e. $\Psi_j$). Each intersection on the heat map corresponds to the value of a Fourier-Taylor coefficient associated to the product of the vertical and horizontal basis functions. This visualization effectively arranges the approximated Fourier-Taylor coefficients into a matrix, constructed as the product of basis functions as described in Eq.~\eqref{eq:bfunc_coupl_red}.
    
\begin{figure*}
    	\centering
    	\includegraphics[width=0.65\textwidth]{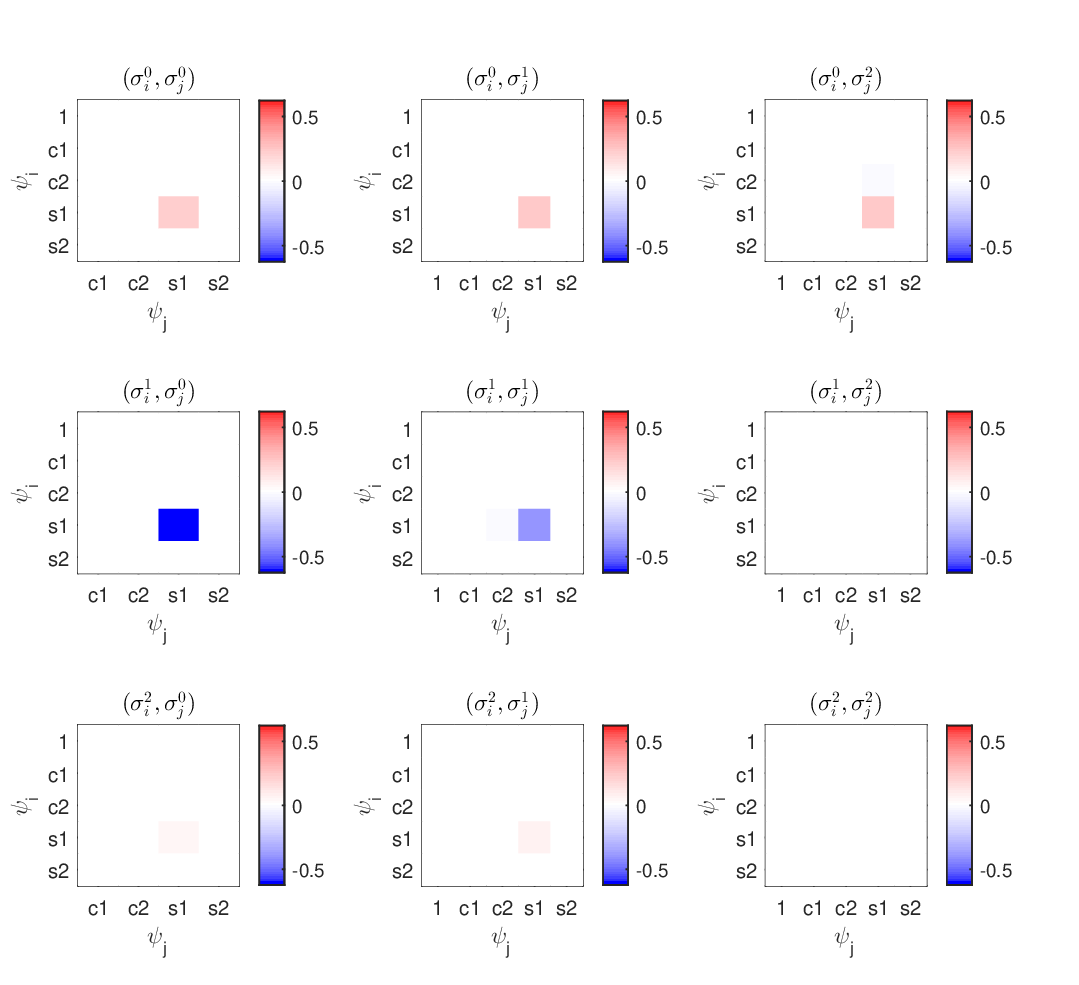}
    	\centering
    	\includegraphics[width=0.65\textwidth]{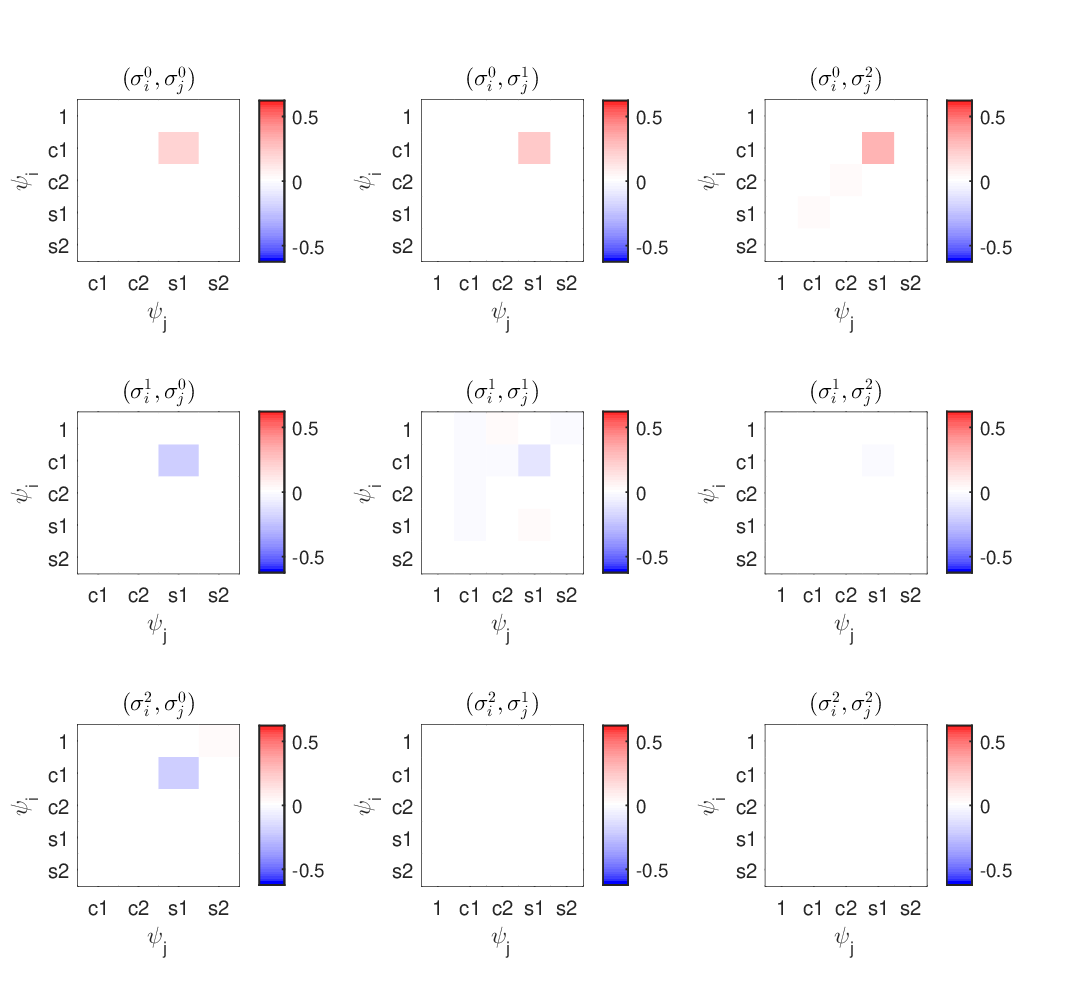}
    	\caption{An example of the computed Fourier-Taylor coefficients of the reconstructed coupling in the amplitude equation (top panels) and phase equation (bottom panels) for two uni-directionally coupled radial isochron clocks (see Appendix~\ref{sec:radisocl}). \label{fig:coupl_RadIsoCl_ampphase}}
    \end{figure*}

\begin{figure*}
	\centering
	\includegraphics[width=0.65\textwidth]{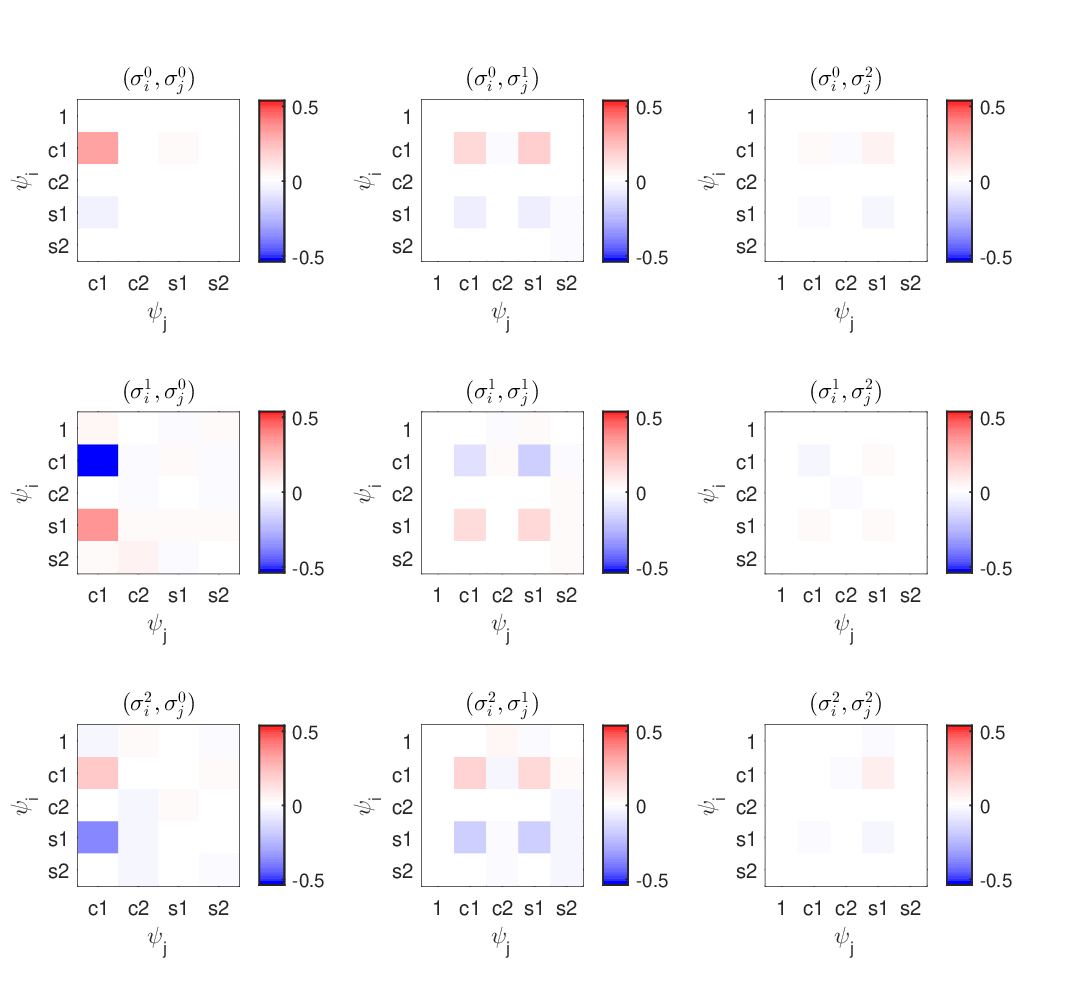}
	\includegraphics[width=0.65\textwidth]{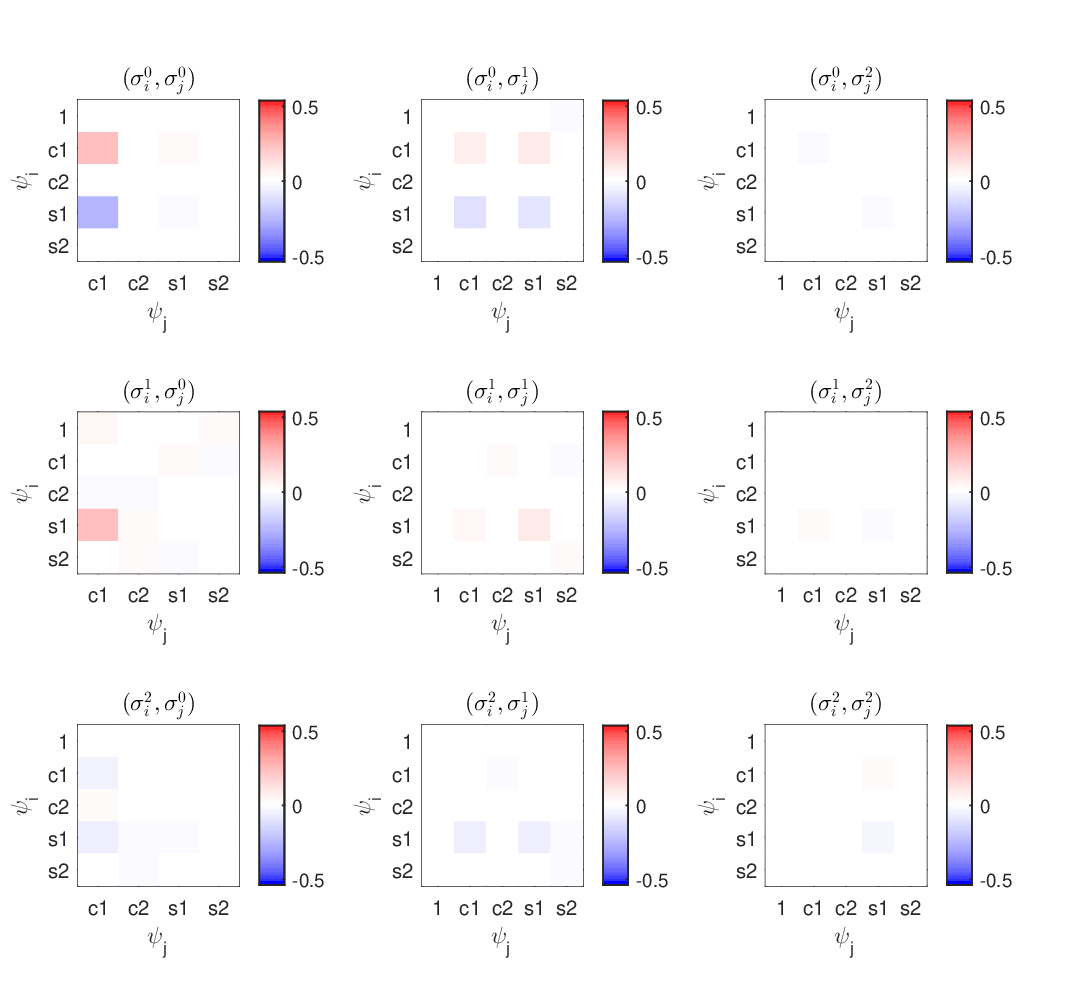}
	\caption{An example of the computed Fourier-Taylor coefficients of the reconstructed coupling in the amplitude equation (top panels) and phase equation (bottom panels) for two uni-directionally coupled canonical models (see Appendix~\ref{sec:canmod}). \label{fig:coupl_ampphase_canmod}}
\end{figure*}

Examples of the reconstructed coupling coefficients for the radial isochron clock and for the canonical model are presented in Figs.~\ref{fig:coupl_RadIsoCl_ampphase} and \ref{fig:coupl_ampphase_canmod}. For these models, the analytic expression of the transformation functions, as well as the coupling functions can be derived (see Appendix~\ref{sec:models} and  Supplementary Materials).  The reconstructed coupling  coefficients for these models  are non-zero only for the terms that have a specific combination of the cosine and sine functions (Figs.~\ref{fig:coupl_RadIsoCl_ampphase} and ~\ref{fig:coupl_ampphase_canmod}). For example, for the radial isochron clock, the coupling coefficients for the phase equation are non-zero for the terms involving $\cos(\varphi_i)\sin(\varphi_j)$ (see Fig~\ref{fig:coupl_RadIsoCl_ampphase} bottom), which agrees with the theoretical expression of the coupling  functions in Eq.~\eqref{eq:radisocl_coupl}. Moreover, the signs of the coefficients are also in accordance with the theoretical ones, e.g. for the amplitude coupling we have a positive sign for the terms with $\sigma^0_i\sigma_j^0=1$, $\sigma_j$, $\sigma_j^2$ and $\sigma_i^2$, and a negative sign for the terms with $\sigma_i$, $\sigma_i\sigma_j$. We observe a similar agreement between the theoretical expression of the coupling coefficients (Eq.~\eqref{eq:canmod_coupl_phase}) and the reconstructed coupling coefficients for the canonical model (Fig.~\ref{fig:coupl_ampphase_canmod}). The coupling coefficients for the other models are given in the Supplementary Materials as well. 
        
    \begin{figure*}
    	\centering
    	\includegraphics[width=0.65\textwidth]{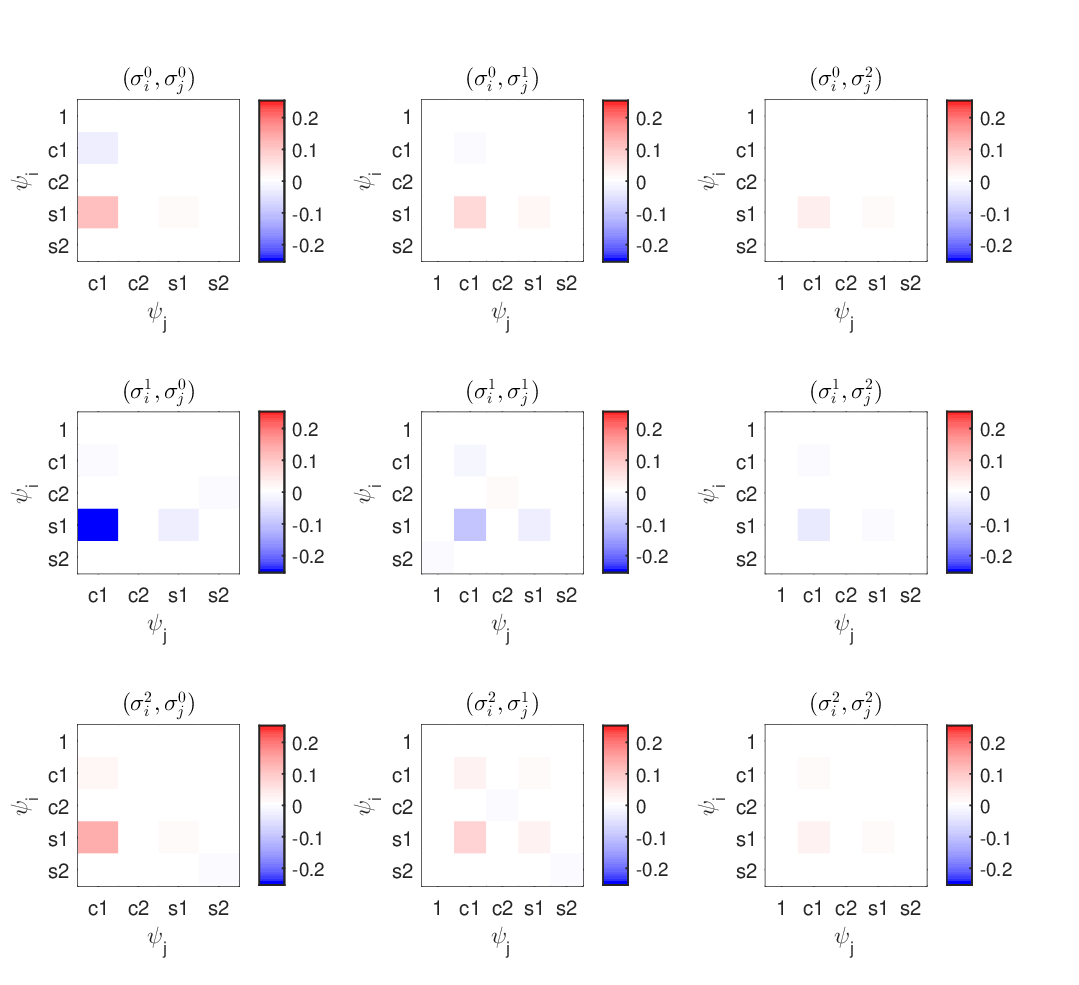}
    	\includegraphics[width=0.65\textwidth]{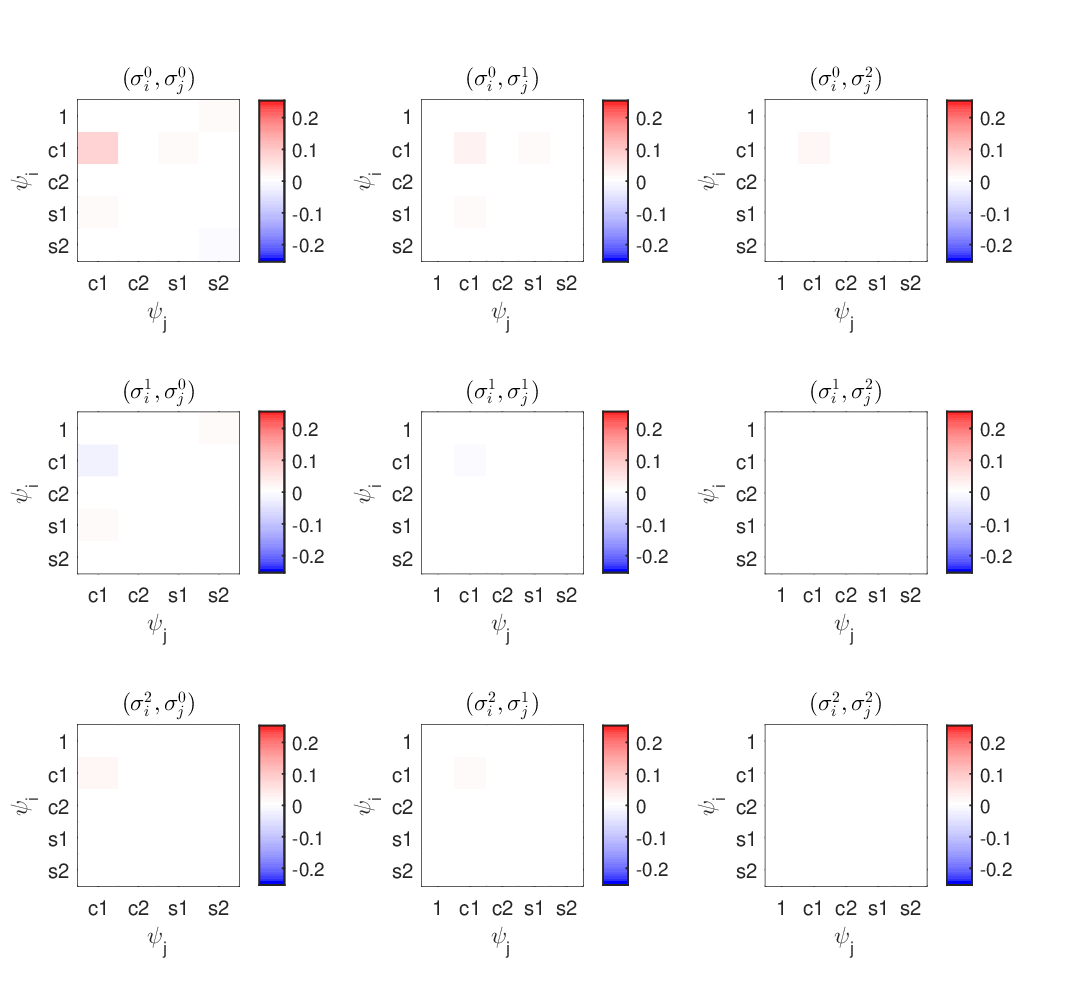}
    	\caption{An example of the computed Fourier-Taylor coefficients of the reconstructed coupling in the amplitude equation (top panels) and phase equation (bottom panels) for two uni-directionally coupled van der Pol oscillators (see Appendix~\ref{sec:vdp}). \label{fig:coupl_vdp_ampphase}}
    \end{figure*}

For the van der Pol oscillator, the reconstructed coupling coefficients are shown in Fig.~\ref{fig:coupl_vdp_ampphase}. In this case, we do not have an analytic expression neither for the transformation functions nor for the coupling part of the VF in the reduced space. However, it is noteworthy that the non-zero coefficients of the coupling part  form a grid-like pattern, suggesting that the basis functions can be sparsified. Moreover, the grid form suggests that the coupling part can be approximated (factorized) as the product of two functions: one  consisting of the basis functions that depend only on the oscillator's own variables (i.e. $\sigma_i$ and $\varphi_i$, represented on the vertical axes), and the other consisting of basis functions that depend on the input variables (i.e. $\sigma_j$ and $\varphi_j$, represented on the horizontal axes). Thus, the former function corresponds to either the amplitude or the phase response functions of the oscillator, while the latter corresponds to the input. This representation of the coupling part as a product of response functions and input is possible due to the linear coupling in the original van der Pol oscillator system (see Eq.~\eqref{eq:vdp}). In contrast, the coefficients of the coupling part for the Wilson-Cowan models do not exhibit a clear grid form (see 
corresponding figures in the Supplementary Materials), since the system non-linearly depends on the input (see Eq.~\eqref{eq:wcmod}). 
    
For all models, the calculated coefficients of the coupling part of the first oscillator (without input) were negligibly small, which demonstrates the ability of the presented method to detect the directionality of the coupling. The results of the estimated coupling coefficients for all models and oscillators are given in the Supplementary Materials. 
    
\section{Discussion}
 
We have presented a new method for the reconstruction of the phase-amplitude dynamics from high time resolution signals measured   from oscillatory systems. This method is intended to be applied to   electrophysiological signals from the brain, such as Magneto- or Electroencephalography (M/EEG), Electrocorticography (ECoG), and intracranial EEG (iEEG), to build a dynamical system that reproduces the transient synchronization properties between brain regions within and between different frequency bands observed in the experiments. 

The presented method combines recent advances in the field of phase-amplitude reduction. Namely, the parameterization method \cite{Cabre2003a, Cabre2003b, Cabre2005, Guillamon2009, Huguet2013, Perez-Cervera2019, Perez-Cervera2020} and the time averaging approach for determining the transformation functions \cite{Mauroy2012, Mauroy2013, Mauroy2018} based on the Koopman operator framework \cite{Mezic2004, Mezic2005, Budisic2012}. The approximation of the parameterization functions from Refs.~\onlinecite{Guillamon2009, Perez-Cervera2020} allowed us to clarify the relationship between the eigenfunctions of the Koopman operator and the transformation functions. Moreover, using this approximation of the parameterization, we derived a method for evaluating the time averaging integral that is well-suited for numerical estimation (Eq.~(\ref{eq:avgsimRt}) and Fig.~\ref{fig:sig0calc}).

The results of the phase-amplitude reconstruction across different models show an accurate VF approximation for both uncoupled and coupled oscillatory systems. Furthermore, the approximated transformation functions show a robust agreement with the theoretical ones for the coupled radial isochron clocks and the coupled canonical models, which have an anlytical expression to compare with. Indeed, the approximated coupling coefficients show good agreement with the theoretical values. When our method was applied to the van der Pol oscillators and the Wilson-Cowan models, which exhibit  non-linear interdependence of the phase and amplitude in their transformation functions (not known analytically), compelling results were also obtained. 

The goal of reconstructing the phase-amplitude dynamics in the reduced space is to develop a representation of the system's dynamics that is independent of the specific observables. 
As shown in previous studies on phase reduction \cite{Kralemann2008,Yeldesbay2019}, when reconstructing the reduced phase dynamics from measured signals, it is necessary to take into account the non-circularity of the limit cycle, or equivalently, the asymmetry of the phase distribution. Neglecting the transformation between observable and reduced phases \footnote{Observable and reduced phases were termed proto- and true phases in Kralemann~et~al.,~2008, and observable and theoretical phases in Yeldesbay~et~al.,~2019.} could cause spurious coupling terms. These principles also apply when reconstructing phase-amplitude dynamics from the amplitude and phase of measured signals. Indeed, the system in Eqs.~(\ref{eq:phisigcoupl}) implies that, in the absence of coupling ($g^{(\varphi)}_{i,j}=0$ and $g^{(\sigma)}_{i,j}=0$), the phase grows linearly and the amplitude decays exponentially. If we dismiss the non-circularity of the limit cycle, e.g. by taking the experimentally measured phase and defining the amplitude variable as the radial distance to the limit cycle, the deviation of the observed amplitude decay from the exponential one and the variation of the rotation speed of the observed phase along the limit cycle and in the surrounding region, will appear in the coupling terms. Therefore, the reconstructed transformation functions serve the purpose of describing  these dependencies of the amplitude and the phase on the domain surrounding the limit cycle. 

As previously mentioned, the transformation functions are unique up to a constant factor for the amplitude and a constant shift for the phase. In the method presented, these constants are fixed during the numerical approximation of the transformation functions, ensuring the reconstruction of a unique and observable independent dynamical system in the reduced space.

In comparison to several data-driven methods presented in recent studies (e.g. Refs.~\onlinecite{Wilson2021b,Kaiser2021,Wilson2023}), this study utilizes only polar coordinates and Fourier-Taylor approximations. The crucial part of the presented method is determining the VF with the coupling using basis functions from Fourier-Taylor series. The VF approximation method resembles the ones presented in Ref.~\onlinecite{Kralemann2008} for phase coupling and in the sparse indentification of nonlinear dynamics (SINDY) method \cite{Brunton2016,brunton2019data}, but without sparsification of the VF coefficients matrix. The representation of the VF as the sum of the basis functions has another advantage. Albeit we do not make any assumptions about the VF of the uncoupled system and assume pairwise coupling, the approximation with Fourier-Taylor basis functions allows the coupling component of the VF to be treated as an additive (linear) input. This approach simplifies the calculation of the coupling functions in the reduced space.

Another distinction of our method is the approximation of the coupling part as a sum of pairwise interactions, each  approximated using two dimensional Fourier-Taylor series that depend on the  amplitude and phase of both the driving and driven systems in a general form. This approach is requested by the specificity of the high time resolution of brain signals. Unlike data-driven phase-amplitude reduction methods based on dynamic mode decomposition (DMD) with input and control \cite{Tu2014, Williams2015,Proctor2016, Proctor2018}, SINDY \cite{Brunton2016,brunton2019data,Kaiser2021}, first and higher order approximation of the isochron equations \cite{Wilson2020b, Wilson2021b, Ahmed2022, Wilson2023} and the methods for estimating the response functions \cite{Cestnik2022}, our approach does not have access to the input or driving signals coming into the region of interest. Instead, we can only measure the resulting electrical/electromagnetic signals emitted by different brain regions. Therefore, we assume that input to a brain region, along with the interaction of the neuronal populations within that brain region, generates the oscillations at different frequencies observed in the measurements. Thus, the reconstruction method presented here is intended to be used as a tool to describe the transient synchronization between these oscillations, both within and across distant brain regions, in the form of a dynamical system. Our method thereby extends conventional methods -- such as phase synchronization measures in the phase domain, correlation analyses and event-related (de-)synchronization (ERD/ERS) in the amplitude domain, or Phase-Amplitude Coupling (PAC) and similar methods in the phase and amplitude domain -- by not only providing information about the directionality of coupling, but also revealing the intrinsic mechanisms of the dynamical interactions.

The presented reconstruction method has several requirements for successful application that we next discuss in detail. (1) A good approximation of the VF in the domain outside of the coverage of the data points is essential. This was achieved by applying the regularized least squares method with generalized cross-validation (GCV) (see Appendix~\ref{sec:rlsq}). The regularization constrains the coefficient values by means of the $L^2$ norm, helping to prevent overfitting of the data points. However, effective reconstruction of the VF with coupling requires a good coverage of the observable amplitude and phase domains. From a practical point of view, this requires a data set (trials) that starts from a range of different amplitude and phase values. (2) Our method can be applied to signals originated from any network of oscillatory systems, provided that a limit cycle exists for the sytem under consideration. Thus, clear oscillatory patterns must be present in the measurements for meaningful interpretation of the results. (3) Our method is intended to be applied to oscillatory signals filtered from the actual measurements around a given frequency. Thus, the noise in the signals should be filtered out as much as possible in order to obtain a good approximation. (4) For the extraction of the phase and amplitude information from the signal, one should use adapted methods that reduce the dependence of the extracted information from the future and the past of the signal, as discussed in Refs.~\onlinecite{Rosenblum2021, Busch2022, Cestnik2022}. 

In conclusion, we present a method for deriving a dynamical system from high time resolution signals generated by a network of interacting oscillators,  particularly, from interacting distant brain regions. The analysis of the resulting dynamical system provides a valuable tool for constructing  mathematical models and investigating neural pathologies and brain disorders. This includes conducting \textit{in silico} experiments with targeted manipulations of activity within specific brain regions or their couplings. Such research can contribute to a deeper understanding of the neurobiological underpinnings of cognitive processes and how network synchronization is impaired by (non-)pathological aging.

\begin{acknowledgments}
This work was funded by the Deutsche Forschungsgemeinschaft (DFG, German Research Foundation) — Project-ID 431549029, SFB 1451 and 491111487. GH was supported by PID-2021-122954NB-I00 funded by MCIN/AEI/ 10.13039/501100011033 and “ERDF: A way of making Europe”, the Maria de Maeztu Award for Centers and Units of Excellence in R\&D (CEX2020-001084-M) and the AGAUR project 2021SGR1039. 
\end{acknowledgments}

\section{Supplementary Material}

In the Supplementary Material section, we provide detailed information about the models used to test our method, together with the derivation of the theoretical expressions for the transformation and the coupling functions for the radial isochron clock and the canonical model. We also present the results of the VF and the reconstruction of the transformation functions, and the estimated coefficients of the coupling functions for all models.

\appendix

\section{Regularized least squares method \label{sec:rlsq}}

In this section, we present the regularized least squares method in the ridge regression formulation \cite{Golub1979,Bates1983,Fenu2017}. 

Let us assume that we have the column $N$-vectors $\mathbf{x}=\{x_i\}^{T}$, $\mathbf{y}=\{y_i\}^{T}$, $i =1,\ldots,N$ and basis functions $\psi_j : \mathbb{R} \rightarrow \mathbb{R}$, $j = 1,\ldots,M$. Let $\boldmath{\Psi}$ be the $N\times M$ matrix of basis functions with elements $\psi_j(x_i)$. We want to approximate $\mathbf{y}$ using the basis functions as
\begin{equation}
\mathbf{y} \approx \boldmath{\Psi}\mathbf{q},
\end{equation}
where $\mathbf{q} = \{q_j\}^{T}$ is the column $M$-vector of coefficients that can be obtained by solving the minimization problem
\begin{equation}
	\label{eq:regminprob}
	\mathbf{q} = \arg \min_{\mathbf{q} \in \mathbb{R}^M} \left\{\|\boldmath{\Psi}\mathbf{q}-\mathbf{y}\|^2 + \kappa\|\mathbf{q}\|^2\right\},
\end{equation}
where $\kappa$ is the regularization parameter. Then, the solution of Eq.~(\ref{eq:regminprob}) is
\begin{equation}
	\mathbf{q} = \left(\kappa \mathbf{I} + \boldmath{\Psi}^{T}\boldmath{\Psi} \right)^{-1} \boldmath{\Psi}^{T} \mathbf{y}.
\end{equation}
We find the optimal regularization parameter using the generalized cross-va\-li\-da\-tion method (GCV) \cite{Golub1979,Bates1983,Fenu2017}. The optimal $\kappa$ is found by minimizing
\begin{equation}
	GCV(\kappa) = \frac{\|\boldmath{\Psi}\mathbf{q}-\mathbf{y}\|^2}{\tau^2}.
\end{equation}
Here, $\tau$ is the effective number of the degree of freedom:
\begin{equation}
	\tau = N - \sum_{l=1}^{r}\frac{\varrho_l^2}{\varrho_l^2+\kappa},
\end{equation}
where $\varrho_l$ are the singular values of $\boldmath{\Psi}$ (or $\varrho_l^2$ are the eigenvalues of $\boldmath{\Psi}^{T}\boldmath{\Psi}$), and $r$ is the rank of the matrices. 

\section{Models}\label{sec:models}
In all simulations, models of two uni-directionally coupled oscillators were used, that is, the second oscillator receives an input from the first oscillator. 

\subsection{Radial isochron clock}\label{sec:radisocl}
The system of equations for the coupled radial isochron clock in polar coordinates reads as \cite{winfree1980, strogatz2000nonlinear}
\begin{equation}\label{eq:radisocl}
\begin{aligned}
	\frac{d\theta_{i}}{dt} & =1+\epsilon_{ij}\frac{r_{j}}{r_{i}}\cos(\theta_{i})\sin(\theta_{j}),\\
	\frac{dr_{i}}{dt} & =a_{i} r_{i}(1-r_{i}^{2})+\epsilon_{ij}\sin(\theta_{i})r_{j}\sin(\theta_{j}), 
\end{aligned}
\end{equation}
where $a_i>0$ is the parameter controlling the decay rate to the limit cycle and $\epsilon_{ij}$ is the parameter controlling the coupling strength.

The theoretical expressions for the forward transformation functions are
\begin{equation}
\begin{aligned}
	K^{(\theta)}_i(\varphi_i,\sigma_i) = & \varphi_i, \\
    K^{(r)}_i(\varphi_i,\sigma_i) = &\sqrt{\frac{1}{1-2a_i\sigma_i}},
\end{aligned}
\end{equation}
and for the inverse transformation functions are
\begin{equation}
    \begin{aligned}
      \Phi_i(\theta_i,r_i) = & \theta_i, \\
	   \Sigma_i(\theta_i,r_i) = & \frac{1}{2a_i}\frac{r_i^2-1}{r_i^2}. \label{eq:transf_inv_radisocl}
    \end{aligned}
\end{equation}

The system of equations for the coupled radial isochron clock in the reduced space approximated using Fourier-Taylor series on both amplitude variables $\sigma_i$ and $\sigma_j$ is 
\begin{multline}
    \frac{d\varphi_{i}}{dt}\approx 1+\epsilon_{ij}\cos(\varphi_{i})\sin(\varphi_{j})\Bigl(1-a_i\sigma_i+a_j\sigma_j-\frac{1}{2}a_i^2\sigma_i^2 \\
    - a_i a_j \sigma_i\sigma_j + \frac{3}{2} a_j^2\sigma_j^2 + \ldots\Bigr),\\
	\frac{d\sigma_{i}}{dt} \approx -2a_{i}\sigma_{i}+\epsilon_{ij}\frac{\sin(\varphi_{i})\sin(\varphi_{j})}{a_{i}}\Bigl(1-3a_{i}\sigma_{i}+a_{j}\sigma_{j} \\
    +\frac{3}{2}a_i^2\sigma_i^2-3a_i a_j \sigma_i \sigma_j+\frac{3}{2} a_j^2 \sigma_j^2 + \ldots\Bigr),\label{eq:radisocl_coupl}
\end{multline}
where we have provided the terms up to order 2 in the amplitude variables $\sigma_i$ and $\sigma_j$. For the derivation of these equations we refer the reader to the Supplementary Materials. 

\subsection{Canonical model}\label{sec:canmod}
The system of equations for the canonical model for an oscillator, also known as Stuart-Landau oscillator \cite{Landau:1944ibh,Stuart_1958,Pikovsky2001}, in polar coordinates reads as follows
\begin{equation}\label{eq:canmod}
\begin{aligned}
    \dot{\theta_{i}} = & 1+\alpha_{i} a_{i}r_{i}^{2}-\epsilon_{ij}\frac{r_{j}}{r_{i}}\sin(\theta_{i})\cos(\theta_{j}), \\
	\dot{r_{i}}	= & \alpha_{i} r_{i}(1-r_{i}^{2})+\epsilon_{ij} r_{j}\cos(\theta_{i})\cos(\theta_{j}),
\end{aligned}
\end{equation}
where oscillator $i$ receives input from  oscillator $j$. The parameter values used in the simulations are: $a_1=1.2$, $a_2=1.0$, $\alpha_1=1.5$, $\alpha_2=2.0$, $\epsilon_{12} = 0.0$, $\epsilon_{21} = 0.3$. 

The theoretical expressions for the forward transformation functions are
\begin{equation}
    \begin{aligned}
	   K^{(\theta)}_i(\varphi_i,\sigma_i)=& \varphi_i + \frac{a_i}{2}\ln(1-2\alpha_i\sigma_i), \\
	   K^{(r)}_i(\varphi_i,\sigma_i)=&\sqrt{\frac{1}{1-2\alpha_i\sigma_i}},
    \end{aligned}
\end{equation}
and for the inverse transformation functions are
\begin{equation}
    \begin{aligned}
	    \Phi_i(\theta_i,r_i)	= & \theta_i+a_i\ln(r_i), \\
        \Sigma_i(\theta_i,r_i)= &\frac{1}{2\alpha_i}\left(1-\frac{1}{r_i^{2}}\right). \label{eq:transf_inv_canmod}
    \end{aligned}
\end{equation}

The system of equations for the canonical model in the reduced space using Fourier-Taylor series on both amplitude variables $\sigma_i$ and $\sigma_j$ is
\begin{widetext}
\begin{align}
	\frac{d\sigma_{i}}{dt} = & -2\alpha_{i} \sigma_{i}+\epsilon_{ij}\Bigl((\cos(\varphi_{i})\cos(\varphi_{j}))/\alpha_{i} + (\sin(\varphi_{i})\cos(\varphi_{j})a_{i}-3\cos(\varphi_{i})\cos(\varphi_{j}))\sigma_{i}
    \nonumber \\
	&+ (\cos(\varphi_{i})\sin(\varphi_{j})a_{j}+\cos(\varphi_{i})\cos(\varphi_{j}))\alpha_{j}\sigma_{j}/\alpha_{i} 
    -((a_{i}^{2}-3)\cos(\varphi_{i})\cos(\varphi_{j})+4a_{i}\sin(\varphi_{i})\cos(\varphi_{j}))\alpha_{i}\sigma_{i}^{2}/2 
    \nonumber\\
	& +(\sin(\varphi_{i})\sin(\varphi_{j})a_{i}a_{j}+\sin(\varphi_{i})\cos(\varphi_{j})a_{i}
    - 3\cos(\varphi_{i})\sin(\varphi_{j})a_{j}-3\cos(\varphi_{i})\cos(\varphi_{j}))\alpha_{j}\sigma_{j}\sigma_{i} 
    \nonumber \\
	&  -((a_{j}^{2}-3)\cos(\varphi_{i})\cos(\varphi_{j})-4a_{j}\cos(\varphi_{i})\sin(\varphi_{j}))\alpha_{j}^{2}\sigma_{j}^{2}/(2\alpha_{i})
    + \ldots \Bigr). 
    \nonumber \\
	\frac{d\varphi_{i}}{dt}  = & 1+\epsilon_{ij}\Bigl(- \sin(\varphi_{i})\cos(\varphi_{j})+\cos(\varphi_{i})\cos(\varphi_{j})a_{i}
    +  
    (a_i^{2}+1)\sin(\varphi_{i})\cos(\varphi_{j})\alpha_i \sigma_{i} 
    \nonumber \\
	& +\bigl[a_{i}a_{j}\cos(\varphi_{i})\sin(\varphi_{j})+a_{i}\cos(\varphi_{i})\cos(\varphi_{j})  
    -a_{j}\sin(\varphi_{i})\sin(\varphi_{j})-\sin(\varphi_{i})\cos(\varphi_{j})\bigr]\alpha_{j}\sigma_{j}
    \nonumber \\
	& 
    +\bigl[- \cos(\varphi_{i})\cos(\varphi_{j})a_{i}^{3}+\sin(\varphi_{i})\cos(\varphi_{j})a_{i}^{2}
    -\cos(\varphi_{i})\cos(\varphi_{j})a_{i}+\sin(\varphi_{i})\cos(\varphi_{j})\bigr]\alpha_{i}^{2}\sigma_{i}^{2}/2 
    \nonumber \\
	& +\bigl[a_{i}^{2}a_j\sin(\varphi_{i})\sin(\varphi_{j})+a_{i}^{2}\sin(\varphi_{i})\cos(\varphi_{j})
    +a_{j}\sin(\varphi_{i})\sin(\varphi_{j})+\sin(\varphi_{i})\cos(\varphi_{j})\bigr]\alpha_{i}\alpha_{j}\sigma_{j}\sigma_{i}    
\nonumber \\
	&  -\bigl[(a_{j}^{2}-3)a_{i}\cos(\varphi_{i})\cos(\varphi_{j})+(3-a_{j}^{2})\sin(\varphi_{i})\cos(\varphi_{j}) -4a_{i} a_{j}\cos(\varphi_{i})\sin(\varphi_{j})+4a_{j}\sin(\varphi_{i})\sin(\varphi_{j})\bigr]\alpha_{j}^{2}\sigma_{j}^{2}/2+\ldots \Bigr), \label{eq:canmod_coupl_phase}
\end{align}
\end{widetext}
where we have provided the terms up to order 2 in the amplitude variables $\sigma_i$ and $\sigma_j$. For the derivation of the equations we refer the reader to the Supplementary Materials.

\subsection{Van der Pol oscillator}\label{sec:vdp}
The system of equations for the van der Pol oscillator in Cartesian coordinates \cite{vanderPol1920,Bogolubov1961,strogatz2000nonlinear} reads as follows
\begin{equation}\label{eq:vdp}
    \begin{aligned}
        	\frac{dx_{i}}{dt}=& y_{i}, \\
	        \frac{dy_{i}}{dt}=& \mu_{i} (1-x_{i}^{2})y_{i} - x_{i} + \epsilon_{ij} x_{j},
    \end{aligned}
\end{equation}
where an oscillator $i$ receives an input from an oscillator $j$. The parameter values used in the simulations are: $\mu_1=0.3$, $\mu_2=0.5$, $\epsilon_{12} = 0.0$, $\epsilon_{21} = 0.1$. 

\subsection{Wilson-Cowan model}\label{sec:wcmod}
The system of equations for the Wilson-Cowan model \cite{Wilson1972, Hoppensteadt1997} in Cartesian coordinates reads
\begin{equation}\label{eq:wcmod}
\begin{aligned}
	\frac{dx_{i}}{dt}=& -x_{i}+S(\rho_{x,i}+a_{i}x_{i}-b_{i}y_{i}), \\
	\frac{dy_{i}}{dt}=& -y_{i}+S(\rho_{y,i}+c_{i}x_{i}-d_{i}y_{i}+\epsilon_{ij}(x_j-x_{0,j})), 
\end{aligned}
\end{equation}
where oscillator $i$ receives input from oscillator $j$, $x_{0,j}$ is the $x$-component of the unstable equilibrium point for the WC model corresponding to uncoupled oscillator $j$, that can be found by solving the equation:
\begin{align}
 	0 = &-x_{0,j}+S(\rho_{x,j}+a_j x_{0,j}-b_j y_{0,j}), \\
 	0 = & -y_{0,j}+S(\rho_{y,j}+c_j x_{0,j}-d_j y_{0,j}).
\end{align}
The parameter values used in the simulations are: $a_1=a_2=b_1=b_2=c_1=c_2=10$, $d_1=d_2=-2$, $\rho_{x,1}=\rho_{x,2}=0.0$, $\rho_{y,1}=-6.75$,
$\rho_{y,2}=-7.0$, $\epsilon_{12} = 0.0$, $\epsilon_{21} = 1.0$. 

\section{Finding initial conditions in the reduced space}\label{sec:sigma0calc}

\begin{figure}[ht!]
	\includegraphics[width=\columnwidth]{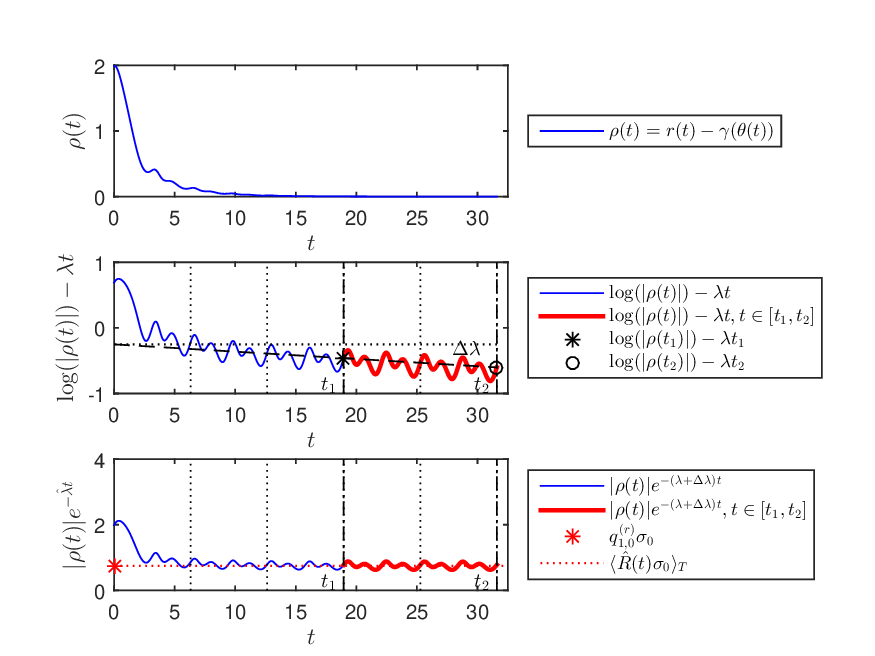}
	\caption{An example of finding the correction to the value of $\lambda$ and  the initial condition $\sigma_0$ for the van der Pol oscillator. \textit{Upper panel}: the time course of an observable $\rho(t)=r(t)-\gamma^{(r)}_0(\theta(t))$ (solid blue line) that vanishes on the limit cycle $\Gamma$. \textit{Middle panel}: the expression $\log(|\rho(t)|)-\lambda t$ (solid blue line) relaxes to a periodic oscillation. The slope of the oscillation $\Delta\lambda$ between the black dashed line and the horizontal black dotted line gives us the correction to $\lambda$. The two vertical dot-dashed lines are in two-period-distance to each other ($t_2=t_1+2T$). The red line marks the time course within $t\in[t_1,t_2]$. The black star and black circle are the first and the last point of this interval, respectively. Vertical dotted lines denote the periods of the oscillation. \textit{Bottom panel}: The blue line is the time course of $|\rho(t)|e^{-(\lambda+\Delta\lambda)t}$. The vertical dotted lines and the two vertical dot-dashed lines ($t_1$ and $t_2$) are the same as in the middle panel. The time course within the interval $t\in[t_1,t_2]$ is marked by the red line, the mean of which is equal to $|q^{(r)}_{1,0} \sigma_0|$ (red star). The horizontal red dotted line is $\langle|\sigma_0  \hat{R}(t)|\rangle_{T}$.\label{fig:sig0calc_vdp}}
\end{figure}

In general, in the limit when $\mathcal{T}\rightarrow \infty$ the expression for $R(t)$ in Eq.~(\ref{eq:limRt}) is oscillatory, as shown in Fig.~\ref{fig:sig0calc_vdp} for the van der Pol oscillator. The system is integrated for five periods, enough time for the trajectory to relax to the limit cycle (upper panel). The red line in the middle panel of Fig.~\ref{fig:sig0calc_vdp} marks the time course of $\log(|\rho(t)|)-\lambda t$, $t\in[t_1, t_2]$, where $t_2$ is the time point at the threshold of numerical precision and $t_1=t_2-2T$. If the initial approximation of $\lambda$, obtained using the monodromy matrix from the solution of Eq.~\eqref{eq:monmatrix}, is not enough accurate the function $\log(|\rho(t)|)-\lambda t$ at the limit oscillates around a non-zero slope line, as shown in the middle panel in Fig.~\ref{fig:sig0calc_vdp} (here for demonstration purposes the initial approximation of $\lambda$ was changed). The correction $\Delta\lambda$ to the value of $\lambda$ can be estimated as the slope of the line that passes through the points $\log(|\rho(t_1)|)-\lambda t_1$ and $\log(|\rho(t_2)|)-\lambda t_2$. After correction of the value of $\lambda$, the time course of $|\rho(t)|e^{-(\lambda+\Delta\lambda)t}$ tends to an oscillation around a zero slope line and mean equal to $|q^{(r)}_{1,0} \sigma_0|$. 

\section*{References}

\end{document}


\newcommand{\spref}{S}
\renewcommand{\thepage}{\spref\arabic{page}} 
\renewcommand{\thesection}{\spref\arabic{section}}  
\renewcommand{\thetable}{\spref\arabic{table}}  
\renewcommand{\thefigure}{\spref\arabic{figure}}
\renewcommand{\theequation}{\spref\arabic{equation}}


\title{Supplementary Material for: \\ Reconstruction of phase-amplitude dynamics from electrophysiological signals}
	
	\author{Azamat Yeldesbay}
    \affiliation{University of Cologne, Institute of Zoology, Cologne, Germany}
    \affiliation{Cognitive Neuroscience, Institute of Neuroscience and Medicine (INM-3),  Research Centre J\"{u}lich, J\"{u}lich, Germany}
    \email{a.yeldesbay@fz-juelich.de.de}
	\author{Gemma Huguet}%
	\affiliation{Universitat Polit\`{e}cnica de Catalunya, Departament de Matem\`{a}tiques, Barcelona, Spain}%
    \affiliation{Centre de Recerca Matem\`{a}tica, Barcelona, Spain}
    \email{gemma.huguet@upc.edu}
 	\author{Silvia Daun}
	\affiliation{Cognitive Neuroscience, Institute of Neuroscience and Medicine (INM-3), Research Centre J\"{u}lich, J\"{u}lich, Germany}%
    \affiliation{University of Cologne, Institute of Zoology, Cologne, Germany}
    \email{s.daun@fz-juelich.de}
    \date{\today}


\maketitle

\section{Models}
\subsection{Radial isochron clock}
\subsubsection{Isolated oscillator}
The system of equations for an isolated radial isochron clock is
\begin{align}
	\frac{dx}{dt} =& a x - y - ax(x^2+y^2), \\
	\frac{dy}{dt} =& x + a y - ay(x^2+y^2).
\end{align}
In polar coordinates $(r,\theta)$ the system reads as
\begin{equation}\label{eqsm:pol_radisocl}
\begin{array}{rl}
	\dfrac{dr}{dt} = &ar(1-r^2), \\[10pt]
	\dfrac{d\theta}{dt} = &1,
\end{array}
\end{equation}
where $a$ is a real parameter. When $a>0$ the system has a hyperbolic attracting limit cycle with $r=1$. 
\subsubsection{Transformation functions}
\label{secrm:transf_radisocl}
We find the transformations (parameterization) $K^{(r)}(\varphi,\sigma), K^{(\theta)}(\varphi,\sigma)$ using Lie brackets as presented in \cite{Castejon2013}. Denoting the vector field of Eq.~(\ref{eqsm:pol_radisocl}) as $X$ and considering the vector field
\begin{equation}
	Y = \left[
		\begin{array}{c}
			a r^3 \\
			0
		\end{array}
	\right],
\end{equation}
one can see that these two vector fields satisfy $[X,Y]=-2 a Y$, where square brackets denote Lie brackets. Thus, $\lambda = - 2 a$ and we find the transformation functions by solving the equations $\partial_{\sigma}K^{(\theta)}=Y_{2}(K^{(\theta)})$
and $\partial_{\sigma}K^{(r)}=Y_{1}(K^{(r)})$ using the stable limit cycle
solution as the boundary condition. More precisely,
\begin{align}
	\dfrac{\partial K^{(\theta)}}{\partial\sigma} & =0,\\
	\dfrac{\partial K^{(r)}}{\partial\sigma} & = a {K^{(r)}}^{3}. \label{eqrm:int_forwtransf_radisocl}
\end{align} 
From the first equation and using $K^{(\theta)}(\varphi,0)=\varphi$ we have
\[
\int_{\varphi}^{K^{(\theta)}}dK^{(\theta)}=0,
\]
\[
K^{(\theta)} - \varphi=0,
\]
\[
K^{(\theta)}(\varphi,\sigma)=\varphi.
\]
From the second equation and using $K^{(r)}(\varphi,0)=1$ we have 
\[
\int_{1}^{K^{(r)}}\frac{dK^{(r)}}{{K^{(r)}}^{3}}=a\int_{0}^{\sigma}d\sigma,
\]
\[
\left . -\frac{1}{2}\frac{1}{{K^{(r)}}^{2}} \right |_{1}^{K^{(r)}}=-\frac{1}{2} \left (\frac{1}{{K^{(r)}}^{2}}-1 \right )=a\sigma,
\]
\[
\frac{1}{{K^{(r)}}^{2}}=1-2a\sigma,
\]
\[
K^{(r)}(\varphi,\sigma)=\sqrt{\frac{1}{1-2a\sigma}}.
\]

Then, the forward transformation functions in polar coordinates are
\begin{align}
	K^{(r)}(\varphi,\sigma) = &\sqrt{\frac{1}{1-2a\sigma}}, \\
	K^{(\theta)}(\varphi,\sigma) = & \varphi,
\end{align}
and in Cartesian coordinates are
\begin{align*}
	K^{(x)}(\varphi,\sigma) = K^{(r)} (\varphi,\sigma) \cos (K^{(\theta)}(\varphi,\sigma))& =\sqrt{\frac{1}{1-2a\sigma}}\cos(\varphi),\\
	K^{(y)}(\varphi,\sigma) = K^{(r)} (\varphi,\sigma) \sin (K^{(\theta)}(\varphi,\sigma))& =\sqrt{\frac{1}{1-2a\sigma}}\sin(\varphi).
\end{align*}
The corresponding inverse transformation functions are
\begin{align}
	\Sigma(\theta,r) = & \frac{1}{2a}\frac{r^2-1}{r^2},\\
	\Phi(\theta,r) = & \theta. 
\end{align}

\subsubsection{The response functions}
The response functions in polar coordinates are
\begin{align}
\partial_{\sigma}K^{(r)}  & = a\left(\frac{1}{1-2a\sigma}\right)^{\frac{3}{2}}=a r^{3},\\
\partial_{\varphi}K^{(r)}  & = 0,\\
\partial_{\sigma}K^{(\theta)}  & = 0,\\
\partial_{\varphi}K^{(\theta)}  & = 1.
\end{align}
The response functions in Cartesian coordinates are
\begin{align}
	\partial_{\sigma}K^{(x)}=& r^{3}a\cos(\theta), \\
	\partial_{\varphi}K^{(x)}=& -r\sin(\theta), \\
    \partial_{\sigma}K^{(y)}=& r^{3}a\sin(\theta), \\
	\partial_{\varphi}K^{(y)}=& r\cos(\theta).
\end{align}

Derivating the expression $(K^{(x)}(\Phi(x,y),\Sigma(x,y)),K^{(y)}(\Phi(x,y),\Sigma(x,y)))=(x,y)$, we obtain the following relation between partial derivatives of the forward and inverse transformations (parameterizations), see Section 3.2 in \cite{Castejon2013}: 
\begin{equation}
	\label{eqrm:rel_forw_inv_grad}
	\left(\begin{array}{cc}
		\partial_{x}\Phi & \partial_{y}\Phi\\
		\partial_{x}\Sigma & \partial_{y}\Sigma
	\end{array}\right)=\left(\begin{array}{cc}
		\partial_{\varphi}K^{(x)} & \partial_{\sigma}K^{(x)}\\
		\partial_{\varphi}K^{(y)} & \partial_{\sigma}K^{(y)}
	\end{array}\right)^{-1}.
\end{equation}
Thus, we have
\begin{equation}
   \label{eq:Grad}
	\left(\begin{array}{cc}
		\partial_{x}\Phi & \partial_{y}\Phi\\
		\partial_{x}\Sigma & \partial_{y}\Sigma
	\end{array}\right)=
	\left(\begin{array}{cc}
		-r\sin \theta & r^{3}a \cos \theta\\
		r\cos \theta & r^{3}a\sin \theta
	\end{array}\right)^{-1} = 
	\left(\begin{array}{cc}
	-\dfrac{\sin \theta}{r} & \dfrac{\cos \theta }{r}\\[10pt]
	\dfrac{\cos \theta}{a r^{3}} & \dfrac{\sin \theta}{a r^{3}}
	\end{array}\right)
\end{equation}

\subsubsection{Coupled radial isochron clocks}

The system of equations for a radial isochron clock (with index $i$) coupled with another radial isochron clock (with index $j$) via the second variable $y$ is given by
\begin{eqnarray}
	\frac{dx_{i}}{dt} & = & a_{i} x_{i}-y_{i}-a_{i} x_{i}(x_{i}^{2}+y_{i}^{2}),\label{smeq:coupl_radisocl_ampl}\\
	\frac{dy_{i}}{dt} & = & x_{i}+a_{i} y_{i}-a_{i} y_{i}(x_{i}^{2}+y_{i}^{2})+\epsilon_{ij} y_{j}, \label{smeq:coupl_radisocl_phase}
\end{eqnarray}
where $\epsilon_{ij}$ is the coupling strength between the clocks with indices $i$ and $j$. 
The system in polar coordinates reads as
\begin{align*}
	\frac{dr_{i}}{dt} & =a_{i} r_{i}(1-r_{i}^{2})+\epsilon_{ij}\sin(\theta_{i})r_{j}\sin(\theta_{j}),\\
	\frac{d\theta_{i}}{dt} & =1+\epsilon_{ij}\frac{r_{j}}{r_{i}}\cos(\theta_{i})\sin(\theta_{j}).
\end{align*}

Taking into account that $\lambda = -2a$, $\omega=1$ and that the input is $(0,y_j)^T = (0,r_j \sin(\theta_j))^T$, we write the system of reduced amplitude and phase equations as
\begin{eqnarray*}
	\frac{d\sigma_{i}}{dt} & = & -2a_{i} \sigma_{i}+\epsilon_{ij}\partial_{y}\Sigma_{i}(\sigma_{i},\varphi_{i})\cdot K^{(r_{j})}(\sigma_{j},\varphi_{j})\sin(K^{(\theta_{j})}(\sigma_{j},\varphi_{j})), \\
	\frac{d\varphi_{i}}{dt} & = & 1+\epsilon_{ij}\partial_{y}\Phi_{i}(\sigma_{i},\varphi_{i})\cdot K^{(r_{j})}(\sigma_{j},\varphi_{j})\sin(K^{(\theta_{j})}(\sigma_{j},\varphi_{j})),
\end{eqnarray*}
where $\partial_{y}\Sigma_{i}(\sigma_{i},\varphi_{i})$ and $\partial_{y}\Phi_{i}(\sigma_{i},\varphi_{i})$ are the amplitude and phase response functions of Eq.~\eqref{eq:Grad}, presented as functions of $\sigma_i$ and
$\varphi_i$:
\begin{align}
	\partial_{y}\Sigma_{i}=\frac{\sin(\theta_{i})}{a_i r_{i}^{3}}=\frac{1}{a_{i}}(1-2a_{i}\sigma_{i})^{3/2}\sin(\varphi_{i}), \label{eqrm:respfunc_sigma_radisocl}\\
	\partial_{y}\Phi_{i}=\frac{\cos(\theta_{i})}{r_{i}}=(1-2a_{i}\sigma_{i})^{1/2}\cos(\varphi_{i}). \label{eqrm:respfunc_phi_radisocl}
\end{align}
Finally, we obtain the system of equations in the reduced space
\begin{align*}
	\frac{d\sigma_{i}}{dt} & =-2a_{i}\sigma_{i}+\epsilon_{ij}\frac{1}{a_{i}}\frac{(1-2a_{i}\sigma_{i})^{3/2}}{(1-2a_{j}\sigma_{j})^{1/2}}\sin(\varphi_{i})\sin(\varphi_{j}), \\
	\frac{d\varphi_{i}}{dt} & =1+\epsilon_{ij}\frac{(1-2a_{i}\sigma_{i})^{1/2}}{(1-2a_{j}\sigma_{j})^{1/2}}\cos(\varphi_{i})\sin(\varphi_{j}).
\end{align*}
We can approximate the response functions Eqs.~(\ref{eqrm:respfunc_sigma_radisocl})-(\ref{eqrm:respfunc_phi_radisocl}) and the input using Fourier-Taylor series on $\sigma_i$ and $\sigma_j$ around zero:
\begin{align}
	\partial_{y}\Sigma_{i}=&\frac{1}{a_{i}}(1-2a_{i}\sigma_{i})^{3/2}\sin(\varphi_{i})\approx \frac{1}{a_i}\sin(\varphi_{i})\left(1-3a_i\sigma_i+\frac{3}{2}a_i^2 \sigma_i^2 + \frac{1}{2} a_i^3\sigma_i^3 + \ldots\right), \\
	\partial_{y}\Phi_{i}=&(1-2a_{i}\sigma_{i})^{1/2}\cos(\varphi_{i}) \approx \cos(\varphi_{i})\left(1-a_i\sigma_i-\frac{1}{2}a_i^2 \sigma_i^2 - \frac{1}{2} a_i^3\sigma_i^3 + \ldots\right),
\end{align}
and
\begin{align}
	y_j &= K^{(r_{j})}(\sigma_{j},\varphi_{j})\sin(K^{(\theta_{j})}(\sigma_{j},\varphi_{j})) = \frac{\sin(\varphi_j)}{(1-2a_{i}\sigma_{i})^{1/2}} \\
   &\approx \sin(\varphi_j) \left(1+a_j\sigma_j+\frac{3}{2}a_j^2 \sigma_j^2 + \frac{5}{2} a_j^3\sigma_j^3 + \ldots\right).
\end{align}
Thus,
\begin{align}
	\frac{d\sigma_{i}}{dt} & \approx -2a_{i}\sigma_{i}+\epsilon_{ij}\frac{\sin(\varphi_{i})\sin(\varphi_{j})}{a_{i}}\left(1-3a_i\sigma_i+\frac{3}{2}a_i^2 \sigma_i^2 + \frac{1}{2} a_i^3\sigma_i^3 + \ldots\right)\left(1+a_j\sigma_j+\frac{3}{2}a_j^2 \sigma_j^2 + \frac{5}{2} a_j^3\sigma_j^3 + \ldots\right),  \\
	\frac{d\varphi_{i}}{dt} & \approx 1+\epsilon_{ij}\cos(\varphi_{i})\sin(\varphi_{j})\left(1-a_i\sigma_i-\frac{1}{2}a_i^2 \sigma_i^2 - \frac{1}{2} a_i^3\sigma_i^3 + \ldots\right)\left(1+a_j\sigma_j+\frac{3}{2}a_j^2 \sigma_j^2 + \frac{5}{2} a_j^3\sigma_j^3 + \ldots\right). 
\end{align}

Also we can approximate using Fourier-Taylor series of the coupling part on both variables $\sigma_i$ and $\sigma_j$ around zero as 
\begin{align}
	\frac{d\sigma_{i}}{dt} & \approx -2a_{i}\sigma_{i}+\epsilon_{ij}\frac{\sin(\varphi_{i})\sin(\varphi_{j})}{a_{i}}\left(1-3a_{i}\sigma_{i}+a_{j}\sigma_{j}+\frac{3}{2}a_i^2\sigma_i^2-3a_i a_j \sigma_i \sigma_j+\frac{3}{2} a_j^2 \sigma_j^2 + \ldots\right),  \label{smeq:radisocl_coupl_ampl}\\
	\frac{d\varphi_{i}}{dt} & \approx 1+\epsilon_{ij}\cos(\varphi_{i})\sin(\varphi_{j})\left(1-a_i\sigma_i+a_j\sigma_j-\frac{1}{2}a_i^2\sigma_i^2 - a_i a_j \sigma_i\sigma_j + \frac{3}{2} a_j^2\sigma_j^2 + \ldots\right). \label{smeq:radisocl_coupl_phase}
\end{align}

\subsection{Canonical model}
\subsubsection{Isolated oscillator}
The system of equations for an isolated canonical model is
\begin{align}
	\frac{dx}{dt} =& \alpha x(1-(x^2+y^2))-y(1+\alpha a (x^2+y^2)), \\
	\frac{dy}{dt} =& \alpha y(1-(x^2+y^2))+x(1+\alpha a (x^2+y^2)).
\end{align}
In polar coordinates the system reads as
\begin{equation}\label{eqsm:pol_canmod}
	\begin{array}{rl}
		\dfrac{dr}{dt} = &\alpha r(1-r^{2}), \\[10pt]
		\dfrac{d\theta}{dt} = &1+\alpha a r^{2}.
	\end{array}
\end{equation}

\subsubsection{Transformation functions}
Following \cite{Castejon2013} we use Lie brackets to find the transformations (parameterization) as presented in section \ref{secrm:transf_radisocl} for the radial isochron clock. Denoting the vector field of equation Eq.~(\ref{eqsm:pol_canmod}) as $X$ and taking the following vector field
\begin{equation}
	Y = \left[
	\begin{array}{c}
		\alpha r^3 \\
		-\alpha a r^2
	\end{array}
	\right],
\end{equation}
we can show that they satisfy $[X,Y]=-2 \alpha Y$, where square brackets denote Lie brackets. Thus, for the canonical model $\lambda = - 2 \alpha$. The transformation functions are found by integrating the equations $\partial_{\sigma}K^{(\theta)}=Y_{2}(K^{(\theta)})$
and $\partial_{\sigma}K^{(r)}=Y_{1}(K^{(r)})$, namely
\begin{align}
	\frac{\partial K^{(\theta)}}{\partial\sigma} & = -\alpha a {K^{(r)}}^2, \label{eqrm:int_forwtransf_canmod_th}\\
	\frac{\partial K^{(r)}}{\partial\sigma} & =\alpha {K^{(r)}}^{3},
	\label{eqrm:int_forwtransf_canmod_r}
\end{align} 
and using the boundary condition on the stable limit cycle as $K^{(r)}(\varphi,0)=1$. The integration of Eq.~(\ref{eqrm:int_forwtransf_canmod_r}) is similar to the integration of Eq.~(\ref{eqrm:int_forwtransf_radisocl}) and the obtained transformation function is 
\begin{equation}
	K^{(r)}(\varphi,\sigma)=\sqrt{\frac{1}{1-2\alpha\sigma}}. \label{eqrm:transf_Kr_canmod}
\end{equation}
From Eq.~(\ref{eqrm:int_forwtransf_canmod_th}) and Eq.~\eqref{eqrm:transf_Kr_canmod} and using that $K^{(\theta)}(\varphi,0)=\varphi$ we have
\[
\dfrac{\partial K^{(\theta)}}{\partial \sigma}=-\frac{\alpha a }{1-2\alpha\sigma},
\]
\[
\int_{\varphi}^{K^{(\theta)}}dK^{(\theta)}=-\int_0^{\sigma}\frac{\alpha a }{1-2\alpha\sigma}d\sigma,
\]
\[
K^{(\theta)}-\varphi = \frac{a}{2}\ln(1-2\alpha\sigma). 
\]
Thus,
\begin{equation}
	K^{(\theta)}(\varphi,\sigma) = \varphi + \frac{a}{2}\ln(1-2\alpha\sigma), \quad
	K^{(r)}(\varphi,\sigma)=\sqrt{\frac{1}{1-2\alpha\sigma}}. \label{eqrm:transf_Kth_canmod}
\end{equation}
The inverse transformation functions are
\begin{align}
	\Sigma(r,\theta)= &\frac{1}{2\alpha}\left(1-\frac{1}{r^{2}}\right), \\
	\Phi(r,\theta)	= & \theta+a\ln(r).
\end{align}

\subsubsection{The response functions}
The response functions in polar coordinates read
\begin{align}
	\partial_{\sigma}K^{(r)}=& \alpha \left(\frac{1}{1-2 \alpha \sigma}\right)^{\frac{3}{2}}=\alpha r^{3},\\
	\partial_{\varphi}K^{(r)}=& 0,\\
	\partial_{\sigma}K^{(\theta)}=& -\frac{a\alpha}{1-2\alpha\sigma}=-a\alpha r^{2},\\
	\partial_{\varphi}K^{(\theta)}=& 1.
\end{align}
Using that $K^{(x)}=K^{(r)} \cos K^{(\theta)}$ and  $K^{(y)}=K^{(r)} \sin K^{(\theta)}$, the response functions in cartesian coordinates are
\begin{align}
	\partial_{\sigma}K^{(x)}=& r^{3}\alpha\left(\cos \theta +a \sin \theta \right), \\
	\partial_{\sigma}K^{(y)}=& r^{3}\alpha\left(-a\cos \theta + \sin \theta \right), \\
	\partial_{\varphi}K^{(x)}=& -r\sin \theta, \\
	\partial_{\varphi}K^{(y)}=& r\cos \theta.
\end{align}

Using Eq.~(\ref{eqrm:rel_forw_inv_grad}) we find
\begin{equation}
	\left(\begin{array}{cc}
		\partial_{x}\Phi & \partial_{y}\Phi\\
		\partial_{x}\Sigma & \partial_{y}\Sigma
	\end{array}\right) = 
	\frac{1}{\alpha r^{3}}
	\left(\begin{array}{cc}
		-r^{2}\alpha\left(\sin(\theta)-a\cos(\theta)\right) & r^{2}\alpha\left(a\sin(\theta)+\cos(\theta)\right)\\
		\cos(\theta) & \sin(\theta)
	\end{array}\right),
\end{equation}
or as a functions of $\sigma$ and $\varphi$:
\begin{align}
	\partial_{x}\Phi = & -\left(1-2\alpha\sigma\right)^{\frac{1}{2}}\left(\sin(\varphi+\frac{a}{2}\ln(1-2\alpha\sigma))-a\cos(\varphi+\frac{a}{2}\ln(1-2\alpha\sigma))\right), \label{eqrm:resp_Phix_canmod}\\
	\partial_{y}\Phi = & \left(1-2\alpha\sigma\right)^{\frac{1}{2}}\left(a\sin(\varphi+\frac{a}{2}\ln(1-2\alpha\sigma))+\cos(\varphi+\frac{a}{2}\ln(1-2\alpha\sigma))\right), \\
	\partial_{x}\Sigma = & \frac{\left(1-2\alpha\sigma\right)^{\frac{3}{2}}}{\alpha}\cos(\varphi+\frac{a}{2}\ln(1-2\alpha\sigma)), \label{eqrm:resp_Sigx_canmod}\\
		\partial_{y}\Sigma = & \frac{\left(1-2\alpha\sigma\right)^{\frac{3}{2}}}{\alpha}\sin(\varphi+\frac{a}{2}\ln(1-2\alpha\sigma)). 
\end{align}
\subsubsection{Coupled canonical models}
The system of equations that describes the coupled canonical models is expressed in Cartesian coordinates as 
\begin{align}
	\dot{x_{i}} = & \alpha_{i} x_{i}(1-(x_{i}^2+y_{i}^2))-y_{i}(1+\alpha_{i} a_{i}(x_{i}^2+y_{i}^2))+\epsilon_{ij} x_{j},\\
	\dot{y_{i}} = & \alpha_{i} y_{i}(1-(x_{i}^2+y_{i}^2))+x_{i}(1+\alpha_{i} a_{i}(x_{i}^2+y_{i}^2)), 
\end{align}
where $j$ is the index of another canonical model coupled to the one with index $i$, $\epsilon_{ij}$ is the coupling strength between oscillators $i$ and $j$. The same system in polar coordinates reads
\begin{align}
	\dot{r_{i}}	= & \alpha_{i} r_{i}(1-r_{i}^{2})+\epsilon_{ij} r_{j}\cos(\theta_{i})\cos(\theta_{j}),\\
	\dot{\theta_{i}} = & 1+\alpha_{i} a_{i}r_{i}^{2}-\epsilon_{ij}\frac{r_{j}}{r_{i}}\sin(\theta_{i})\cos(\theta_{j}).
\end{align}

For the canonical model $\lambda=-2\alpha$ and $\omega=1$ and the input considered is  $(x_j,0)^T = (r_j\cos \theta_j, 0)^T $. Then, the system of equations for the amplitude and phase in the reduced space is
\begin{align}
	\frac{d\sigma_{i}}{dt} = -2\alpha_{i} \sigma_{i}+\epsilon_{ij}\partial_{x}\Sigma_{i}(\sigma_{i},\varphi_{i})\cdot K^{(r_{j})}(\sigma_{j},\varphi_{j})\cos(K^{(\theta_{j})}(\sigma_{j},\varphi_{j})), \\
	\frac{d\varphi_{i}}{dt} = 1+\epsilon_{ij}\partial_{x}\Phi_{i}(\sigma_{i},\varphi_{i})\cdot K^{(r_{j})}(\sigma_{j},\varphi_{j})\cos(K^{(\theta_{j})}(\sigma_{j},\varphi_{j})).
\end{align}
Using Eq.~(\ref{eqrm:transf_Kth_canmod}), (\ref{eqrm:resp_Phix_canmod}), and (\ref{eqrm:resp_Sigx_canmod})  we obtain
\begin{align}
	\frac{d\sigma_{i}}{dt} = & -2\alpha\sigma_{i}+ \frac{\epsilon_{ij}}{\alpha_{i}}\frac{(1-2\alpha_{i}\sigma_{i})^{3/2}}{(1-2\alpha_{j}\sigma_{j})^{1/2}}\cos(\varphi_{i}+\frac{a_{i}}{2}\ln(1-2\alpha_{i}\sigma_{i}))\cos(\varphi_{j}+\frac{a_{j}}{2}\ln(1-2\alpha_{j}\sigma_{j})), \\
	\frac{d\varphi_{i}}{dt} = & 1 - \epsilon_{ij} \frac{(1-2\alpha_{i}\sigma_{i})^{1/2}}{(1-2\alpha_{j}\sigma_{j})^{1/2}}\left(\sin(\varphi_{i}+\frac{a_{i}}{2}\ln(1-2\alpha_{i}\sigma_{i})) \right. \nonumber\\
	& \left.-a_{i}\cos(\varphi_{i}+\frac{a_{i}}{2}\ln(1-2\alpha_{i}\sigma_{i}))\right)\cos(\varphi_{j}+\frac{a_{j}}{2}\ln(1-2\alpha_{j}\sigma_{j})).
\end{align}

The Fourier-Taylor approximation of the response functions on $\sigma_i$ around zero is
\begin{align}
	\partial_{x}\Sigma_{i}(\sigma_{i},\varphi_{i}) \approx & \frac{\cos \varphi_{i} }{\alpha_{i}}-(3\cos \varphi_{i}-a\sin \varphi_{i})\sigma_{i} -\frac{\alpha_{i}}{2}(4a_{i}\sin \varphi_{i}
  +(a_{i}^{2}-3)\cos \varphi_{i})\sigma_{i}^{2} \nonumber\\
 & - \frac{\alpha_{i}^{2}}{6} (a_i^2+1) (a_i \sin \varphi_{i} -3 \cos \varphi_{i})\sigma_{i}^{3}+\ldots,  \\
	\partial_{x}\Phi_{i}(\sigma_{i},\varphi_{i}) \approx & -\sin \varphi_{i} +a_{i}\cos \varphi_{i} + \alpha_{i}(a_{i}^{2}+1)\sin \varphi_{i} \sigma_{i}  
   -\frac{\alpha_{i}^{2}}{2}(a_{i}^2+1) (a_i \cos \varphi_{i} - \sin \varphi_{i})\sigma_{i}^{2} \nonumber \\ 
   &- \frac{\alpha_{i}^{3}}{6}(a_{i}^{2}+1) (4 a_i \cos \varphi_{i} +(a_{i}^{2}-3)\sin \varphi_{i})\sigma_{i}^{3}+\ldots .
\end{align}
The Taylor expansion of the input on $\sigma_j$ around zero is
\begin{multline*}
	K^{(r_{j})}(\sigma_{j},\varphi_{j})\cos(K^{(\theta_{j})}(\sigma_{j},\varphi_{j})) \approx 
	\cos \varphi_j + \alpha_j (a_j \sin \varphi_j +\cos \varphi_j ) \sigma_j 
	-\frac{\alpha_j^2}{2} \left( (a_j^2-3) \cos \varphi_j-4 a_j \sin \varphi_j \right) \sigma_j^{2}+ \ldots . 
\end{multline*}

Thus, the Fourier-Taylor expansion of the system of equations for the canonical system in the reduced space reads as
\begin{align}
	\frac{d\sigma_{i}}{dt} = & -2\alpha_{i} \sigma_{i}+\epsilon_{ij}\Bigl((\cos(\varphi_{i})\cos(\varphi_{j}))/\alpha_{i}+ \nonumber \\
	& (\sin(\varphi_{i})\cos(\varphi_{j})a_{i}-3\cos(\varphi_{i})\cos(\varphi_{j}))\sigma_{i}+(\cos(\varphi_{i})\sin(\varphi_{j})a_{j}+\cos(\varphi_{i})\cos(\varphi_{j}))\alpha_{j}\sigma_{j}/\alpha_{i}- \nonumber\\
	&-((a_{i}^{2}-3)\cos(\varphi_{i})\cos(\varphi_{j})+4a_{i}\sin(\varphi_{i})\cos(\varphi_{j}))\alpha_{i}\sigma_{i}^{2}/2+ \nonumber\\
	& +(\sin(\varphi_{i})\sin(\varphi_{j})a_{i}a_{j}+\sin(\varphi_{i})\cos(\varphi_{j})a_{i}-3\cos(\varphi_{i})\sin(\varphi_{j})a_{j}-3\cos(\varphi_{i})\cos(\varphi_{j}))\alpha_{j}\sigma_{j}\sigma_{i}- \nonumber \\
	&  -((a_{j}^{2}-3)\cos(\varphi_{i})\cos(\varphi_{j})-4a_{j}\cos(\varphi_{i})\sin(\varphi_{j}))\alpha_{j}^{2}\sigma_{j}^{2}/(2 \alpha_{i})+\ldots \Bigr). \label{smeq:canmod_coupl_ampl} \\
	\frac{d\varphi_{i}}{dt}  = & 1+\epsilon_{ij}\Bigl(- \sin(\varphi_{i})\cos(\varphi_{j})+\cos(\varphi_{i})\cos(\varphi_{j})a_{i}+(a_i^{2}+1)\sin(\varphi_{i})\cos(\varphi_{j})\alpha_i\sigma_{i} \nonumber \\
	& +\left[a_{i}a_{j}\cos(\varphi_{i})\sin(\varphi_{j})+a_{i}\cos(\varphi_{i})\cos(\varphi_{j})-a_{j}\sin(\varphi_{i})\sin(\varphi_{j})-\sin(\varphi_{i})\cos(\varphi_{j})\right]\alpha_{j}\sigma_{j}  \nonumber \\
	& +\left[-\cos(\varphi_{i})\cos(\varphi_{j})a_{i}^{3}+\sin(\varphi_{i})\cos(\varphi_{j})a_{i}^{2}-\cos(\varphi_{i})\cos(\varphi_{j})a_{i}+\sin(\varphi_{i})\cos(\varphi_{j})\right]\alpha_{i}^{2}\sigma_{i}^{2}/2 \nonumber \\
	& +\left[a_{i}^{2}a_j\sin(\varphi_{i})\sin(\varphi_{j})+a_{i}^{2}\sin(\varphi_{i})\cos(\varphi_{j})+a_{j}\sin(\varphi_{i})\sin(\varphi_{j})+\sin(\varphi_{i})\cos(\varphi_{j})\right]\alpha_{i}\alpha_{j}\sigma_{j}\sigma_{i} \nonumber \\
	&  -\bigl[(a_{j}^{2}-3)a_{i}\cos(\varphi_{i})\cos(\varphi_{j})+(3-a_{j}^{2})\sin(\varphi_{i})\cos(\varphi_{j}) \nonumber\\
	&
	-4a_{i}a_{j}\cos(\varphi_{i})\sin(\varphi_{j})+4a_{j}\sin(\varphi_{i})\sin(\varphi_{j})\bigr]\alpha_{j}^{2}\sigma_{j}^{2}/2+\ldots \Bigr). \label{smeq:canmod_coupl_phase}
\end{align}

\subsection{Van der Pol oscillator}

The system of equations for a van der Pol (vdP) oscillator (with index $i$) which receives an input from another vdP oscillator (with index $j$) reads
\begin{align}
	\frac{dx_{i}}{dt}	= &y_{i}, \\
	\frac{dy_{i}}{dt}	= & \mu_{i} (1-x_{i}^{2})y_{i} - x_{i} + \epsilon_{ij} x_{j}.
\end{align}

\subsection{Wilson-Cowan model}
The system of equations for an isolated Wilson-Cowan model is
\begin{align}
	\frac{dx_{i}}{dt}	=-x_{i}+S(\rho_{x,i}+a_{i}x_{i}-b_{i}y_{i}), \\
	\frac{dy_{i}}{dt}	=-y_{i}+S(\rho_{y,i}+c_{i}x_{i}-d_{i}y_{i}),
\end{align}
where the function $S(\cdot)$ is defined as 
\begin{equation}
	S(z)=\frac{1}{1+e^{-z}}.
\end{equation}

The system of equations for two coupled Wilson-Cowan models is
\begin{align}
	\frac{dx_{i}}{dt}	=& -x_{i}+S(\rho_{x,i}+a_{i}x_{i}-b_{i}y_{i}), \\
	\frac{dy_{i}}{dt}	=& -y_{i}+S(\rho_{y,i}+c_{i}x_{i}-d_{i}y_{i}+\epsilon_{ij}(x_j-x_{0,j})), \\
	\frac{dx_{j}}{dt}	=& -x_{j}+S(\rho_{x,j}+a_{j}x_{j}-b_{j}y_{j}), \\
	\frac{dy_{j}}{dt}	=& -y_{j}+S(\rho_{y,j}+c_{j}x_{j}-d_{j}y_{j}+\epsilon_{ji}(x_i-x_{0,i})), 
\end{align}
where $x_{0,i}$ and $x_{0,j}$ are the $x$-components of the unstable equilibrium points for the corresponding isolated WC models, that can be found by solving the equation:
\begin{align}
 	0 = &-x_{0}+S(\rho_{x}+a x_{0}-by_{0}), \\
 	0 = & -y_{0}+S(\rho_{y}+c x_{0}-dy_{0}).
\end{align}

Taking into account that
\begin{equation}
	\frac{S(z)}{dz} = S(z)(1-S(z)),
\end{equation}
we can linearize the Wilson-Cowan model for weak coupling ($\epsilon\ll 1$) as follows
\begin{align}
	\frac{dx_{i}}{dt}	\approx & -x_{i}+S(\rho_{x,i}+a_{i}x_{i}-b_{i}y_{i}), \\
	\frac{dy_{i}}{dt}	\approx & -y_{i}+S(\rho_{y,i}+c_{i}x_{i}-d_{i}y_{i}) + S(\rho_{y,i}+c_{i}x_{i}-d_{i}y_{i})\left(1-S(\rho_{y,i}+c_{i}x_{i}-d_{i}y_{i})\right)\epsilon_{ij}(x_j-x_{0,j}).
\end{align}
\clearpage
\newpage
\section{Right hand side approximation}
\begin{figure}[h]
	\centering
	\includegraphics[width=0.5\textwidth]{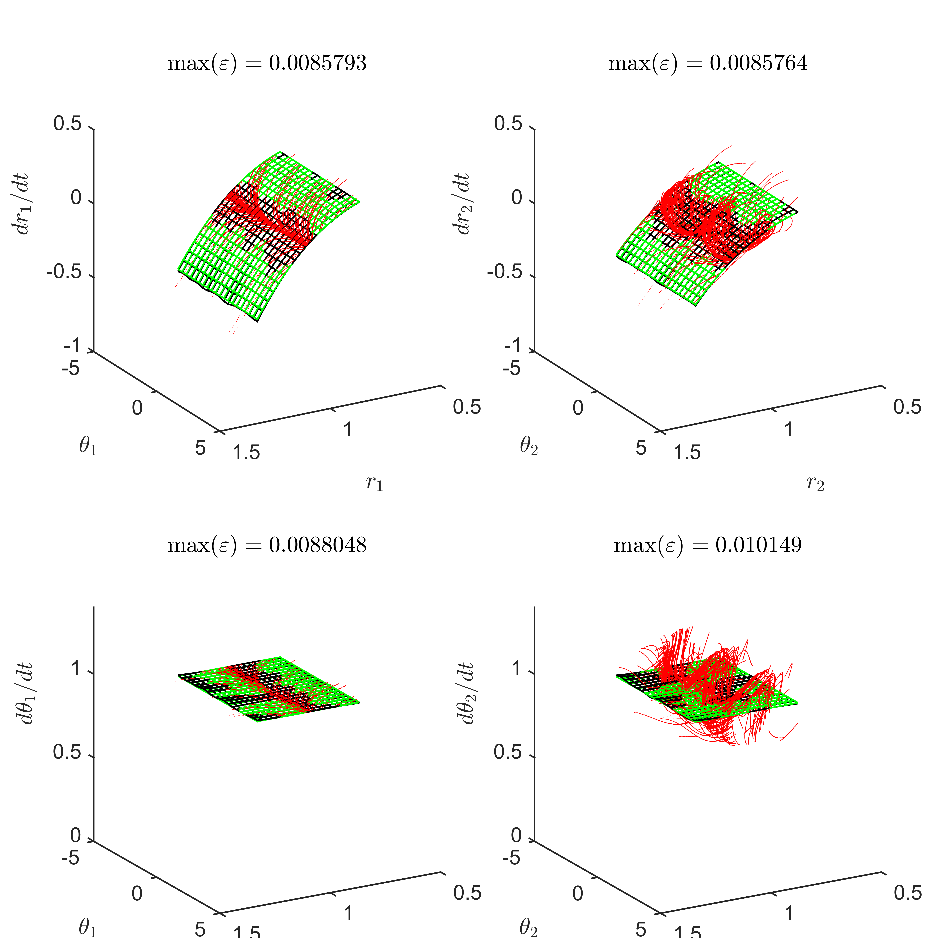}
	\caption{An example of the reconstructed uncoupled VF from simulation of two coupled radial isochron clocks. The red dots are the data points, the green mesh is the theoretical and the black mesh is the reconstructed uncoupled VF. The second oscillator (right column) receives input from the first oscillator (first column). \label{smfig:rhs_RadIsoCl}}
\end{figure}

\begin{figure}
	\centering
	\includegraphics[width=0.5\textwidth]{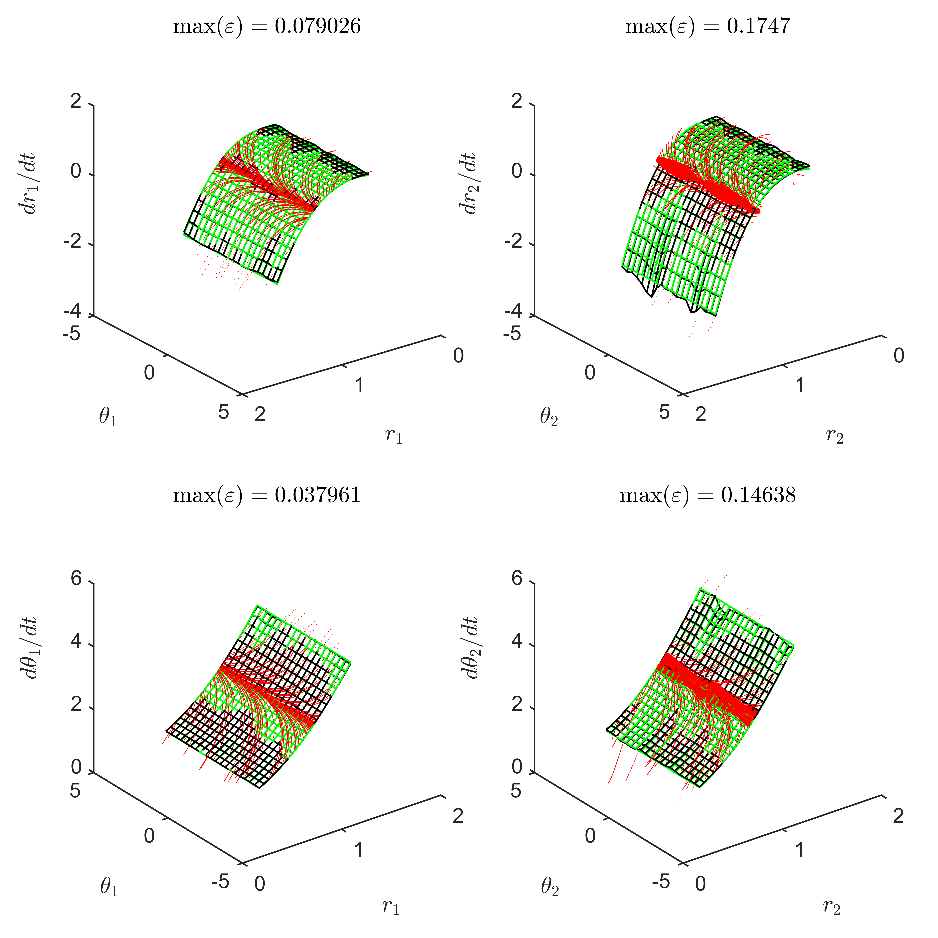}
	\caption{An example of the reconstructed uncoupled VF from simulation of two coupled canonical models. The red dots are the data points, the green mesh is the theoretical and the black mesh is the reconstructed uncoupled VF. The second oscillator (right column) receives input from the first oscillator (first column). \label{smfig:rhs_canmod}}
\end{figure}

\begin{figure}
	\centering
	\includegraphics[width=0.5\textwidth]{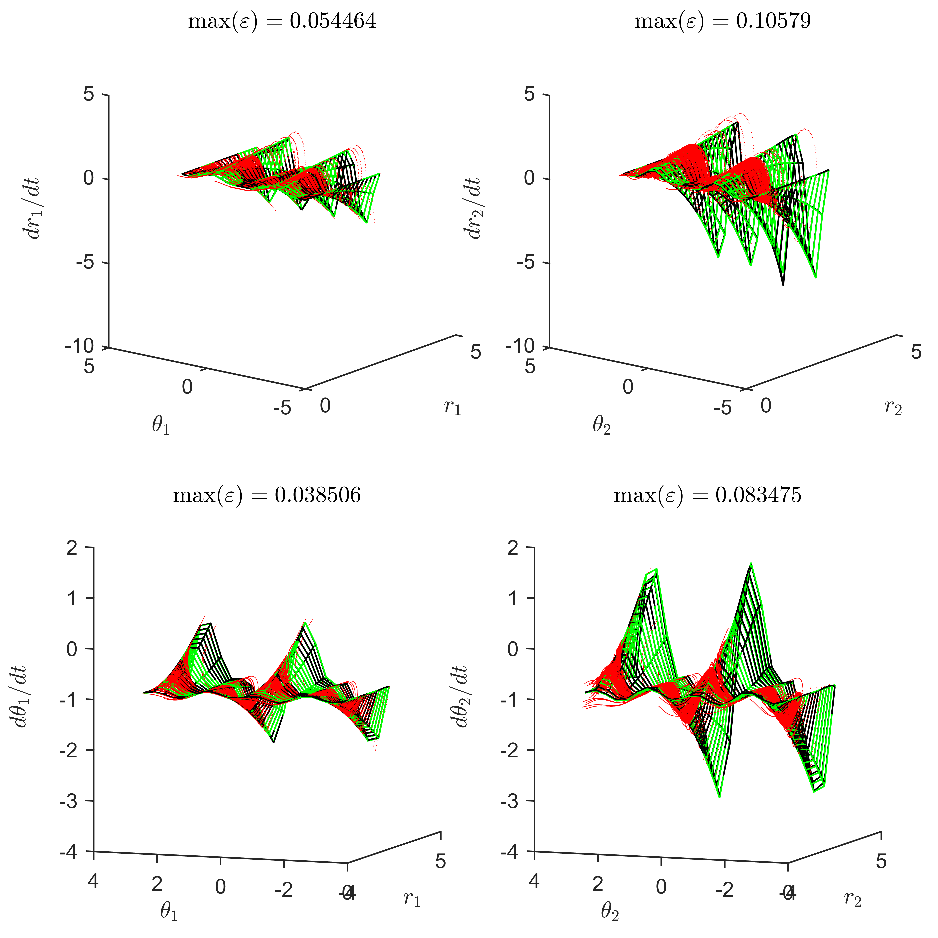}
	\caption{An example of the reconstructed uncoupled VF from simulation of two coupled van der Pol oscillators. The red dots are the data points, the green mesh is the theoretical and the black mesh is the reconstructed uncoupled VF. The second oscillator (right column) receives input from the first oscillator (first column). \label{smfig:rhs_vdp}}
\end{figure}

\begin{figure}
	\centering
	\includegraphics[width=0.5\textwidth]{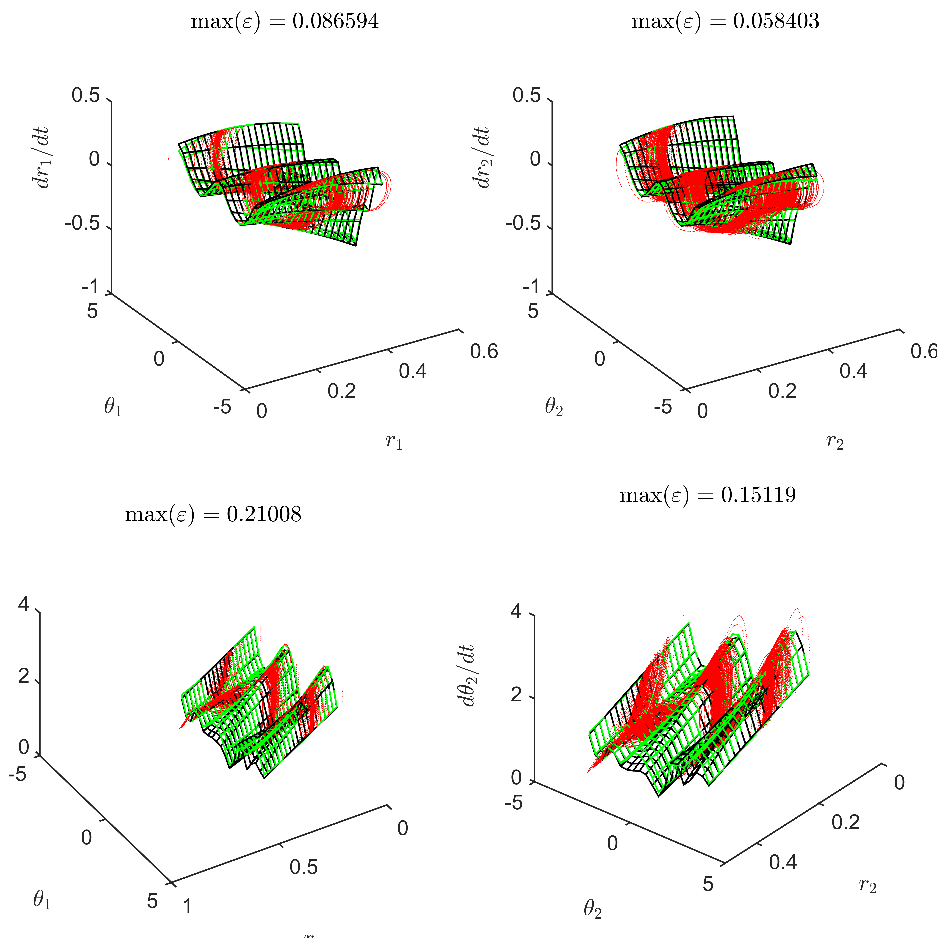}
	\caption{An example of the reconstructed uncoupled VF from simulation of two coupled Wilson-Cowan models. The red dots are the data points, the green mesh is the theoretical and the black mesh is the reconstructed uncoupled VF. The second oscillator (right column) receives input from the first oscillator (first column). \label{smfig:rhs_wcmod}}
\end{figure}

Here, we present examples of the reconstruction of the uncoupled vector field (VF) using synthetically simulated signals (i.e. observables) from different models. All simulated models consist of two uni-directionally coupled oscillators, where the first oscillator is connected to the second one. To assess the goodness of the approximation, the maximal deviation of the reconstructed uncoupled VF from the theoretical one for all coordinates evaluated on the data points is given on the top of every panel of Figs.~\ref{smfig:rhs_RadIsoCl}-\ref{smfig:rhs_wcmod}.

\newpage
\section{Reconstruction of the transformation}
Here, we present the results of the reconstruction of the transformation functions for the different models using the reconstructed uncoupled VF functions. 

\begin{figure}[hb!]
	\centering
	\includegraphics[width=0.65\textwidth]{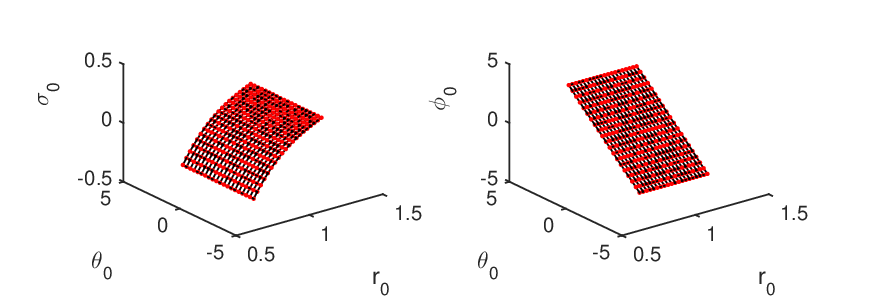}
	\includegraphics[width=0.65\textwidth]{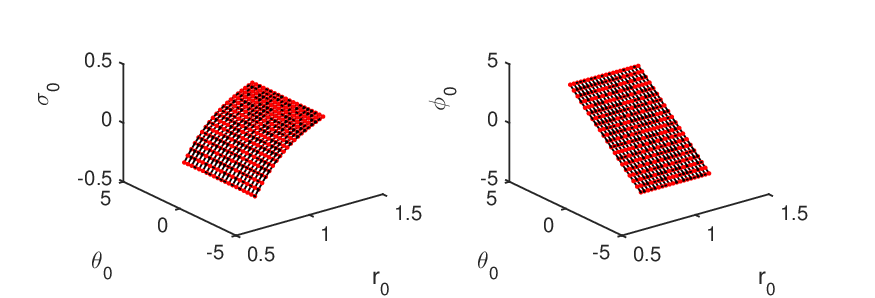}
	\caption{An example of the reconstruction of the inverse transformation functions $\Sigma(\theta,r)$ and $\Phi(\theta,r)$ from the uncoupled VF of two coupled radial isochron clocks. The red dots are the  values of the initial conditions in the reduced space calculated using the averaging method and the black mesh is the transformation function approximated using the initial values (the red dots). Upper panels - the first oscillator, lower - the second. \label{smfig:transf_radisocl}}
\end{figure}

\begin{figure}
	\centering
	\includegraphics[width=0.65\textwidth]{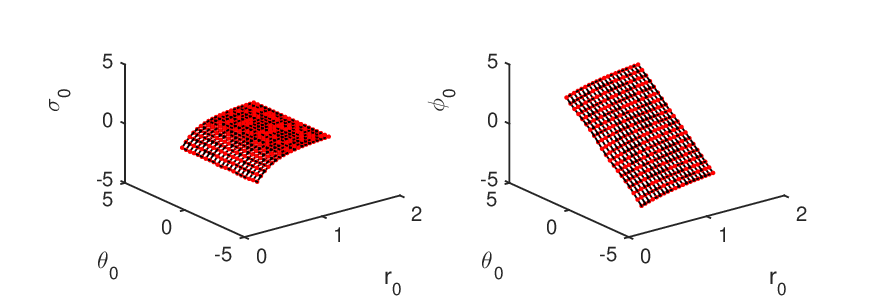}
	\includegraphics[width=0.65\textwidth]{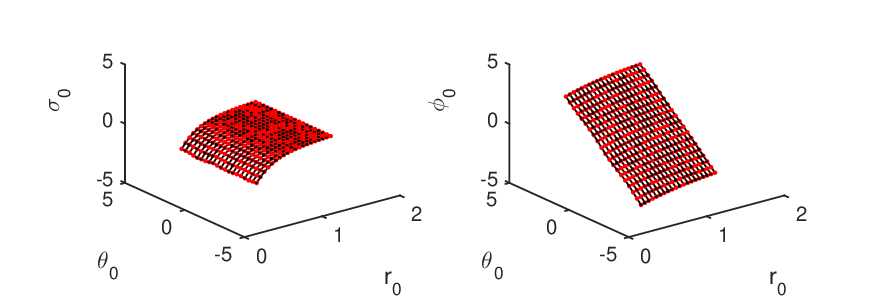}
	\caption{An example of the reconstruction of the inverse transformation functions $\Sigma(\theta,r)$ and $\Phi(\theta,r)$ from the uncoupled VF of two coupled canonical models. The red dots are the values of the initial conditions in the reduced space calculated using the averaging method and the black mesh is the transformation function approximated using the initial values (the red dots). Upper panels - the first oscillator, lower - the second. \label{smfig:transf_canmod}}
\end{figure}

\begin{figure}
	\centering
	\includegraphics[width=0.65\textwidth]{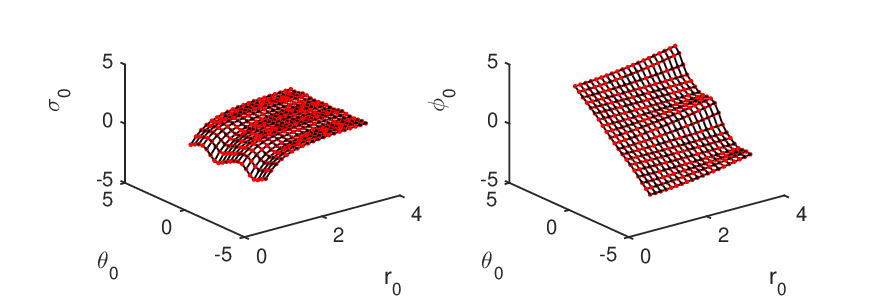}
	\includegraphics[width=0.65\textwidth]{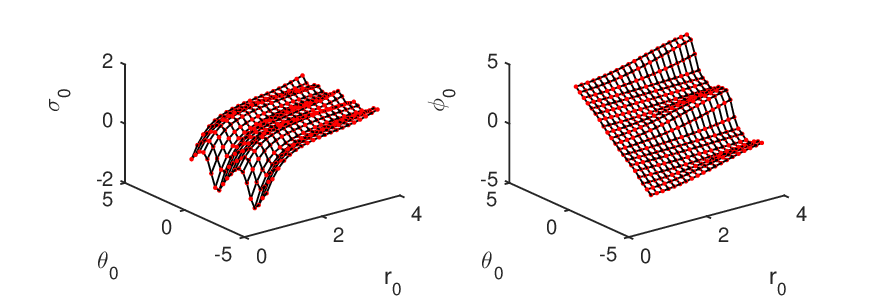}
	\caption{An example of the reconstruction of the inverse transformation functions $\Sigma(\theta,r)$ and $\Phi(\theta,r)$  from the uncoupled VF of two coupled van der Pol oscillators. The red dots are the values of the initial conditions calculated using the averaging method in the reduced space and the black mesh is the transformation function approximated using the initial values (the red dots). Upper panels - the first oscillator, lower - the second. \label{smfig:transf_vdp}}
\end{figure}

\begin{figure}
	\centering
	\includegraphics[width=0.65\textwidth]{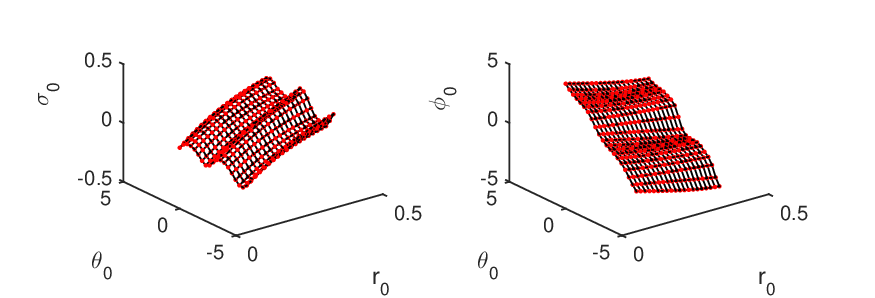}
	\includegraphics[width=0.65\textwidth]{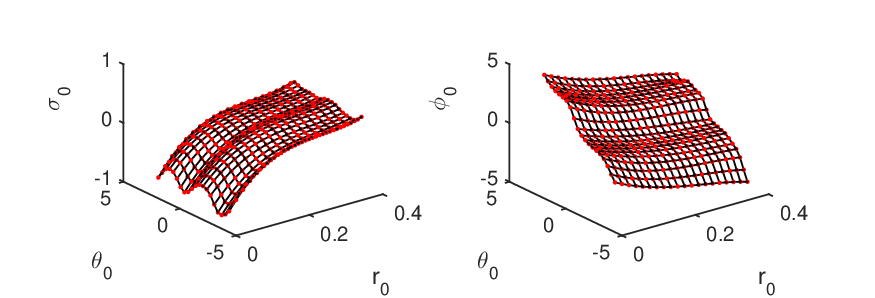}
	\caption{An example of the reconstruction of the inverse transformation functions $\Sigma(\theta,r)$ and $\Phi(\theta,r)$ from the uncoupled VF of two coupled Wilson-Cowan models. The red dots are the values of the initial conditions in the reduced space calculated using the averaging method and the black mesh is the transformation function approximated using the initial values (the red dots). Upper panels - the first oscillator, lower - the second. \label{smfig:transf_wcmod}}
\end{figure}
\clearpage
\section{Reconstruction of the coupling}
Here, we present the results of the reconstruction of the coupling functions for the different models using the reconstructed coupling part of VF and the inverse transformations.

\begin{figure}[hb!]
	\centering
	\includegraphics[width=0.6\textwidth]{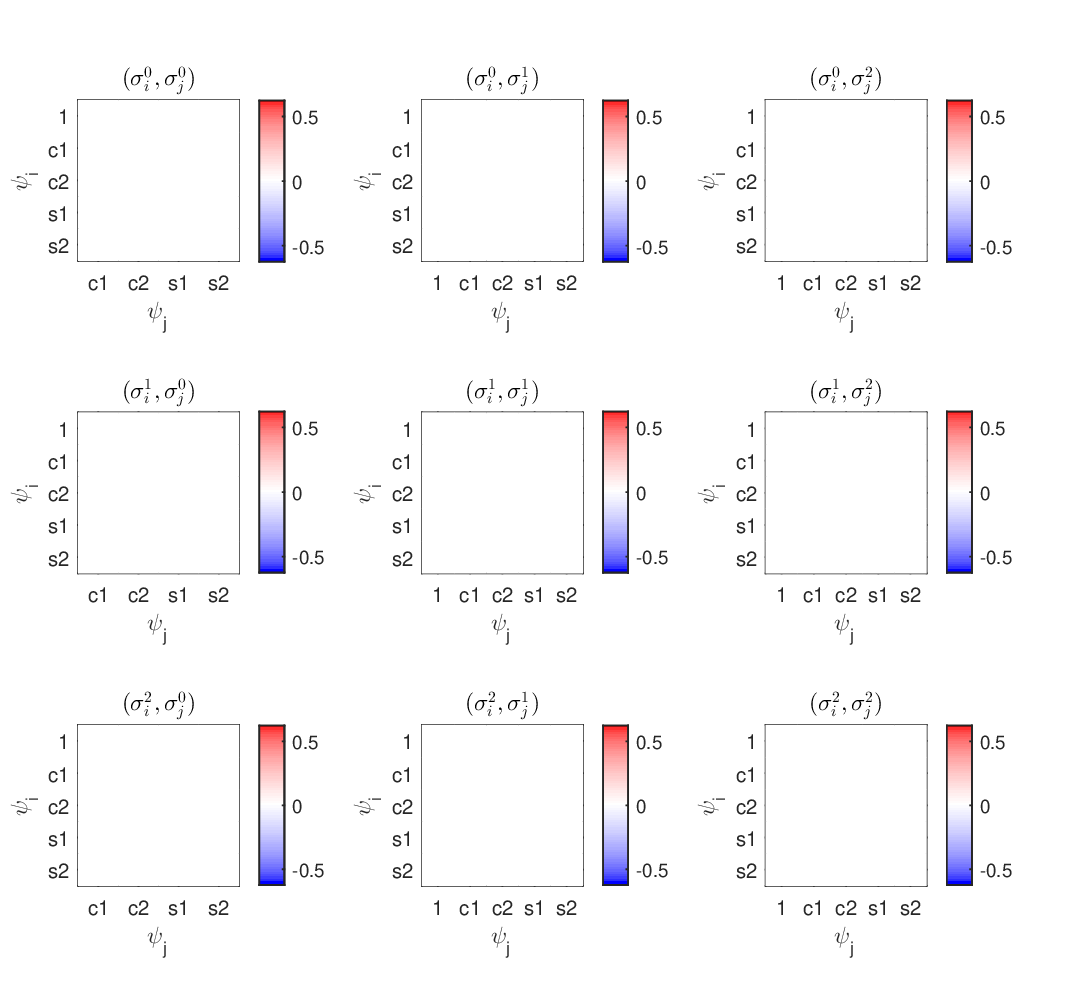}
	\includegraphics[width=0.6\textwidth]{Coupl_radisocl_61_Amp2w1.eps}
	\caption{An example of the reconstructed coupling coefficients in the amplitude equation for two uni-directionally coupled radial isochron clocks. The upper panels correspond to the first oscillator without input, and the bottom panels correspond to the second oscillator, that receives input from the first one. The coefficients are presented as a heat-map (see description in Fig.~\ref{fig:expcoefs} in the main text). \label{smfig:coupl_amp_radisocl}}
\end{figure}

\begin{figure}[h!]
	\centering
	\includegraphics[width=0.6\textwidth]{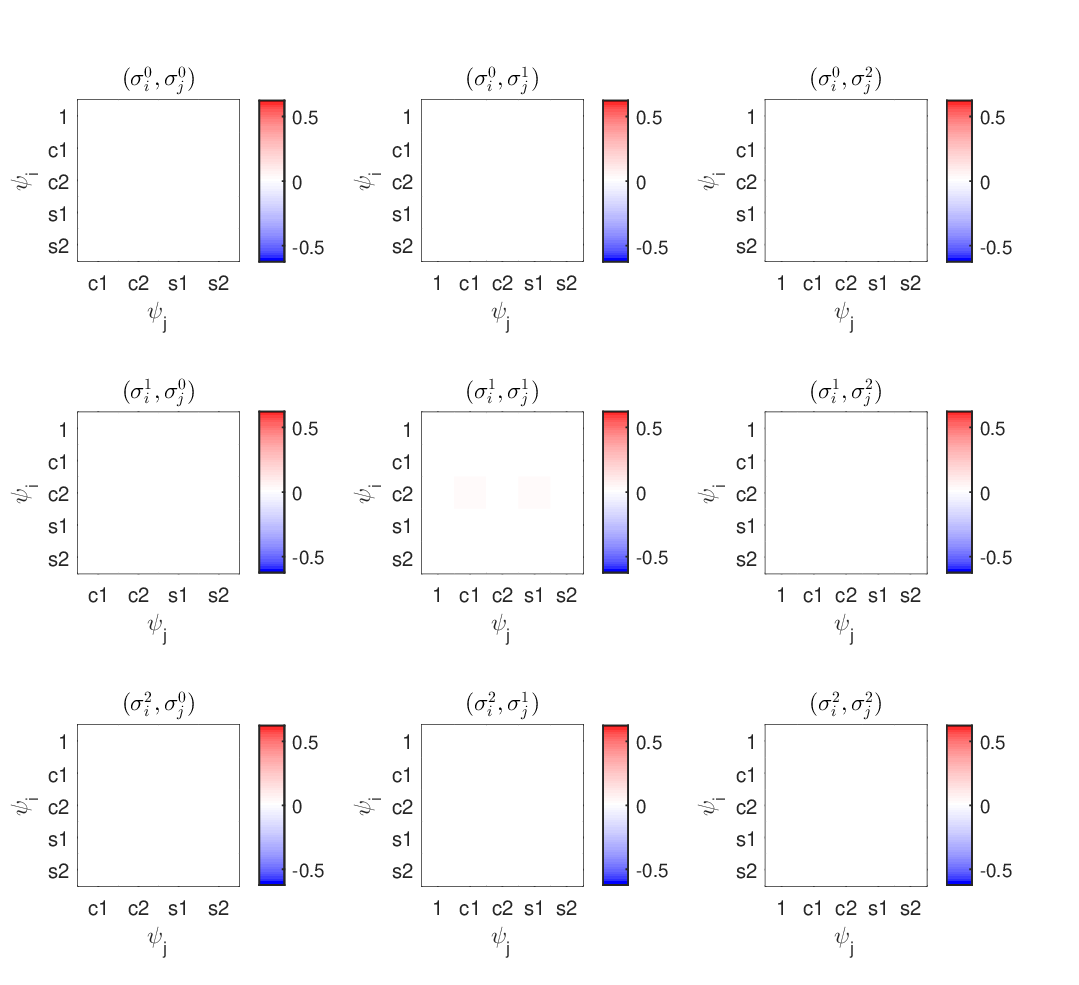}
	\includegraphics[width=0.6\textwidth]{Coupl_radisocl_61_Phase2w1.eps}
	\caption{An example of the reconstructed coupling coefficients in the phase equation for two uni-directionally coupled radial isochron clocks. The upper panels correspond to the first oscillator without input, and the lower panels correspond to the second oscillator, that receives input from the first oscillator. The coefficients are presented as a heat-map (see description in Fig.~\ref{fig:expcoefs} in the main text). \label{smfig:coupl_phase_radisocl}}
\end{figure}

\begin{figure}[h!]
	\centering
	\includegraphics[width=0.6\textwidth]{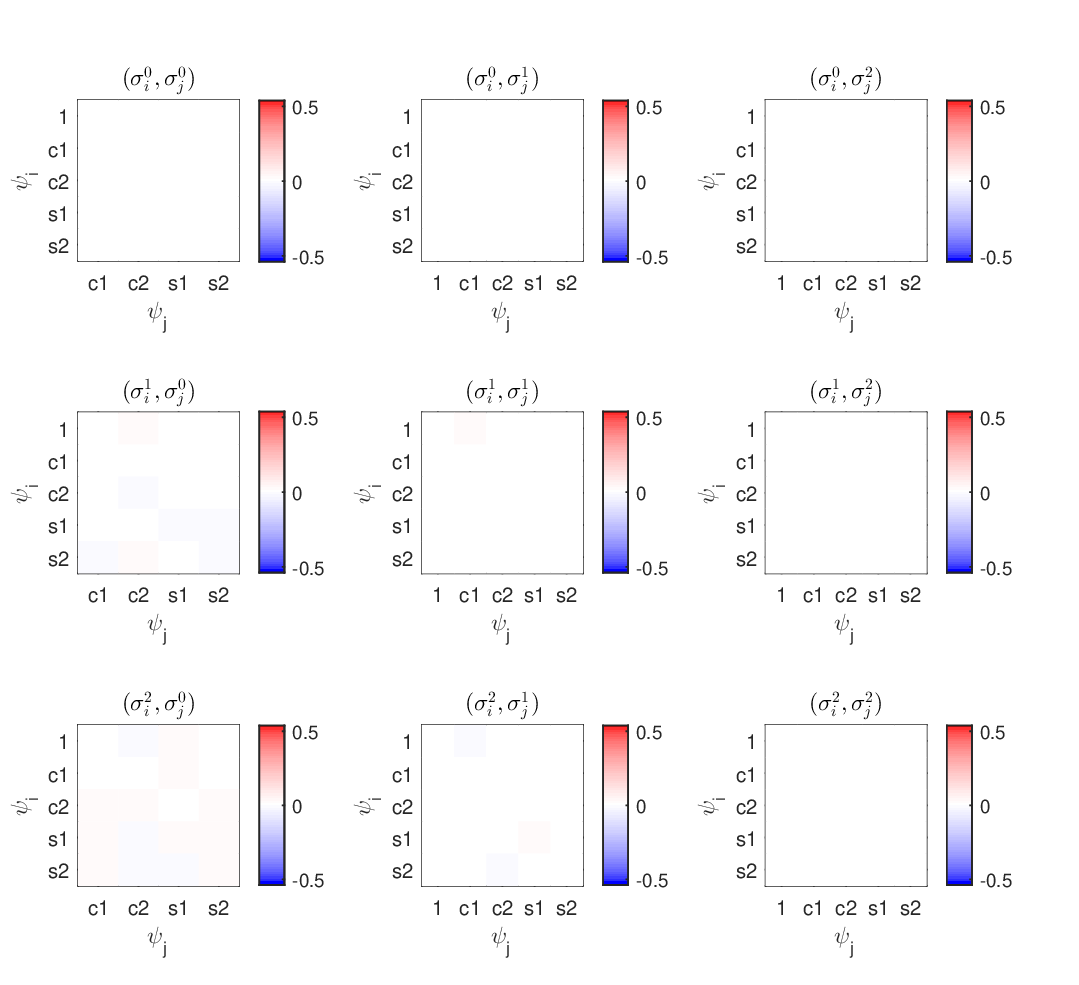}
	\includegraphics[width=0.6\textwidth]{Coupl_canmod_61_Amp2w1.eps}
	\caption{An example of the reconstructed coupling coefficients in the amplitude equation for two uni-directionally coupled canonical models. The upper panels correspond to the first oscillator without input, and the lower panels correspond to the second oscillator, that receives input from the first one. The coefficients are presented as a heat-map (see description in Fig.~\ref{fig:expcoefs} in the main text). \label{smfig:coupl_ampl_canmod}}
\end{figure}

\begin{figure}[h!]
	\centering
	\includegraphics[width=0.6\textwidth]{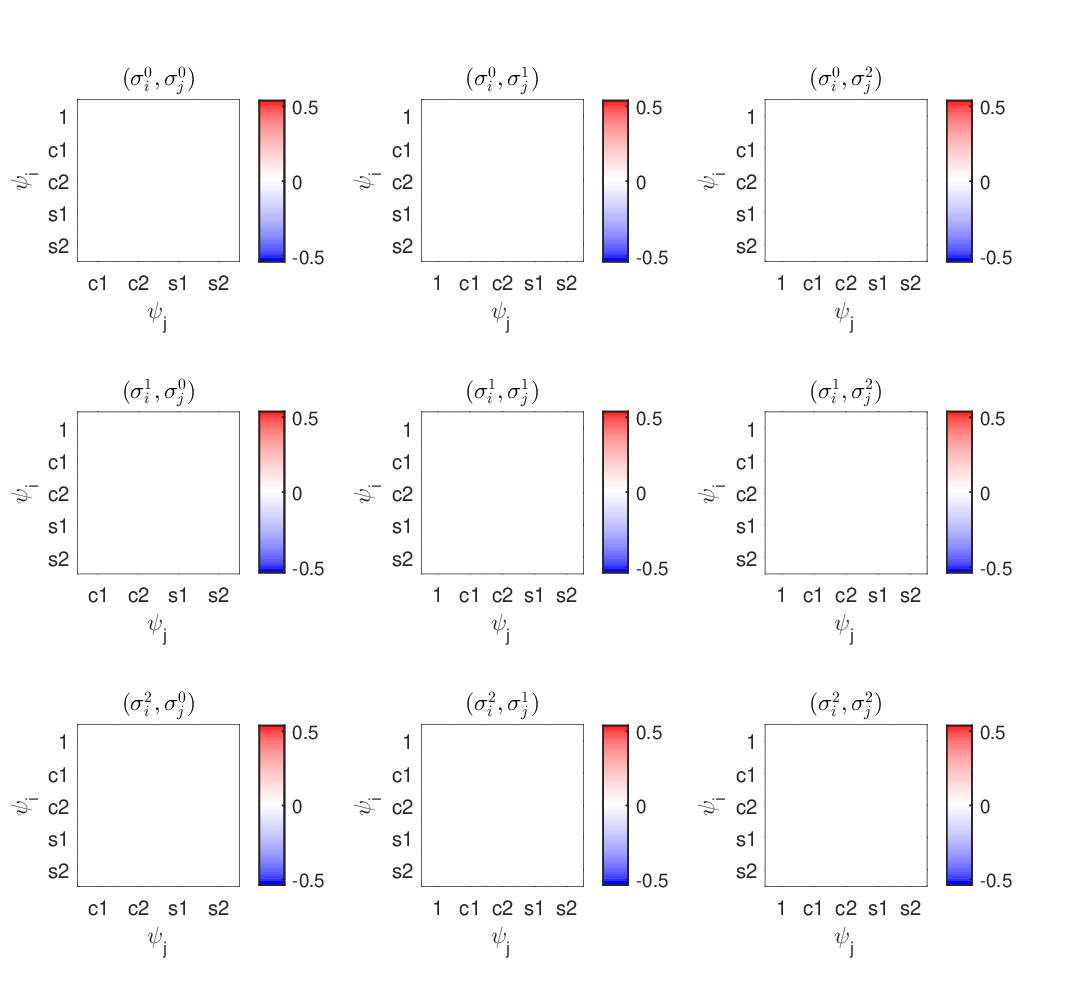}
	\includegraphics[width=0.6\textwidth]{Coupl_canmod_61_Phase2w1.eps}
	\caption{An example of the reconstructed coupling coefficients in the phase equation for two uni-directionally coupled canonical models. The upper panels correspond to the first oscillator without input, and the lower panels correspond to the second oscillator, that receives input from the first one. The coefficients are presented as a heat-map (see description in Fig.~\ref{fig:expcoefs} in the main text). \label{smfig:coupl_phase_canmod}}
\end{figure}

\begin{figure}[h!]
	\centering
	\includegraphics[width=0.6\textwidth]{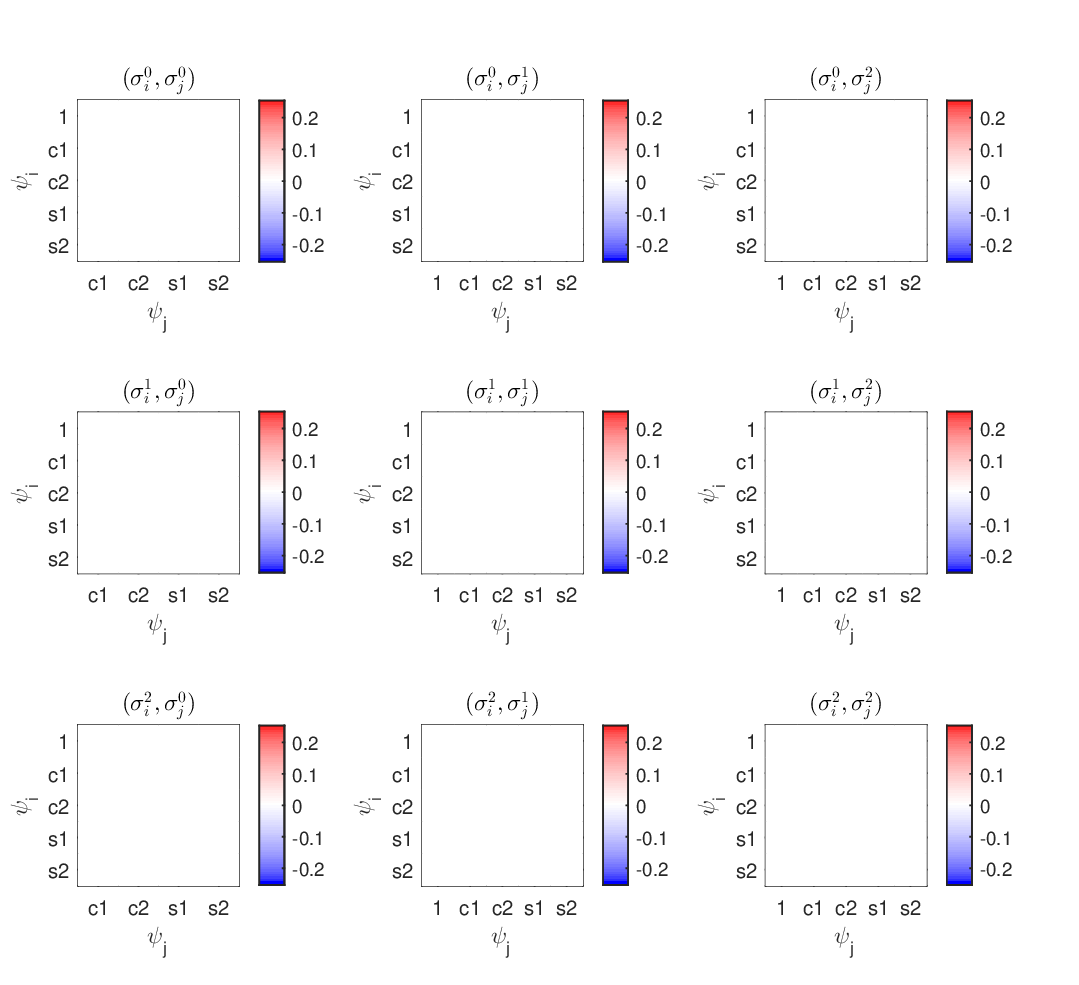}
	\includegraphics[width=0.6\textwidth]{Coupl_vdp_59_Amp2w1.eps}
	\caption{An example of the reconstructed coupling coefficients in the amplitude equation for two uni-directionally coupled van der Pol oscillators. The upper panels correspond to the first oscillator without input, and the lower panels correspond to the second oscillator, that receives input from the first one. The coefficients are presented as a heat-map (see description in Fig.~\ref{fig:expcoefs} in the main text). \label{smfig:coupl_ampl_vdp}}
\end{figure}

\begin{figure}[h!]
	\centering
	\includegraphics[width=0.6\textwidth]{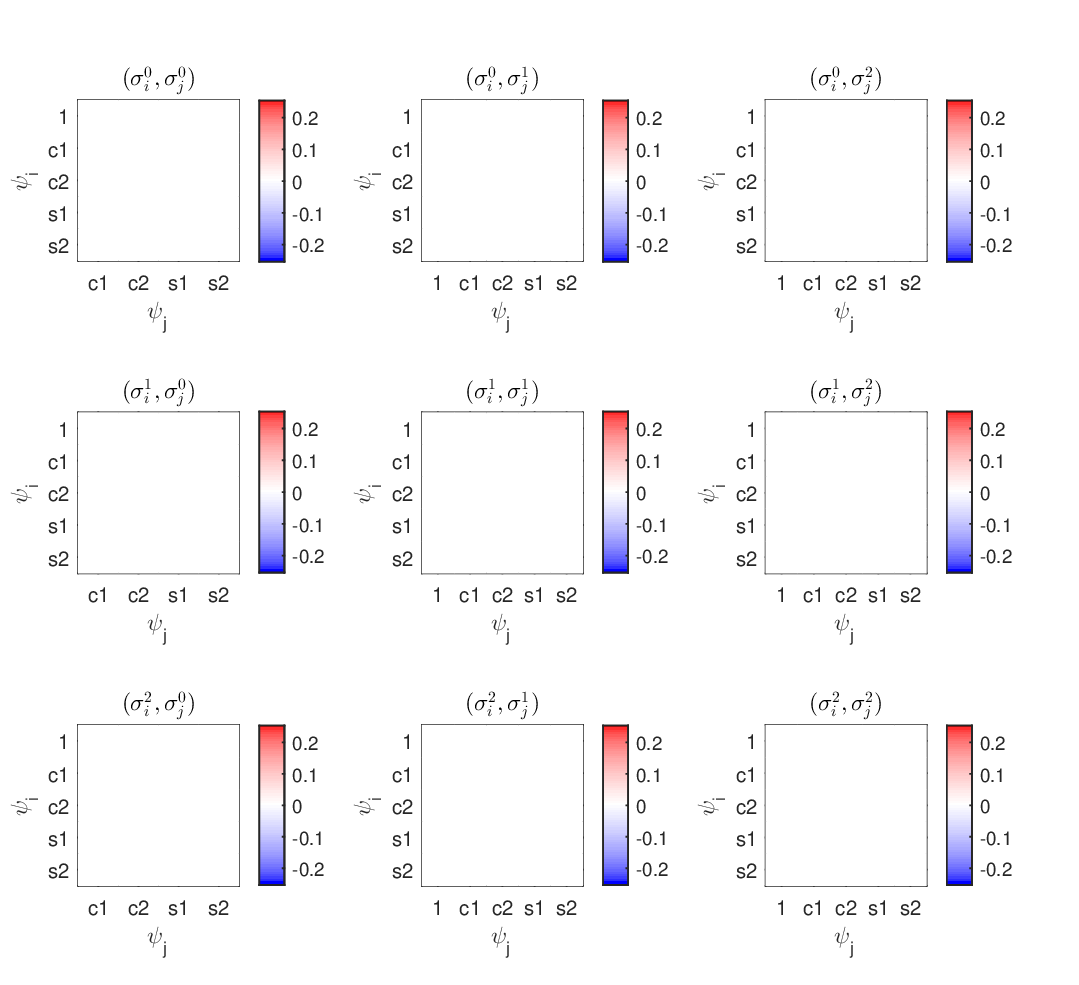}
	\includegraphics[width=0.6\textwidth]{Coupl_vdp_59_Phase2w1.eps}
	\caption{An example of the reconstructed coupling coefficients in the phase equation for two uni-directionally coupled van der Pol oscillators. The upper panels correspond to the first oscillator without input, and the lower panels correspond to the second oscillator, that receives input from the first one. The coefficients are presented as a heat-map (see description in Fig.~\ref{fig:expcoefs} in the main text). \label{smfig:coupl_phase_vdp}}
\end{figure}

\begin{figure}[h!]
	\centering
	\includegraphics[width=0.6\textwidth]{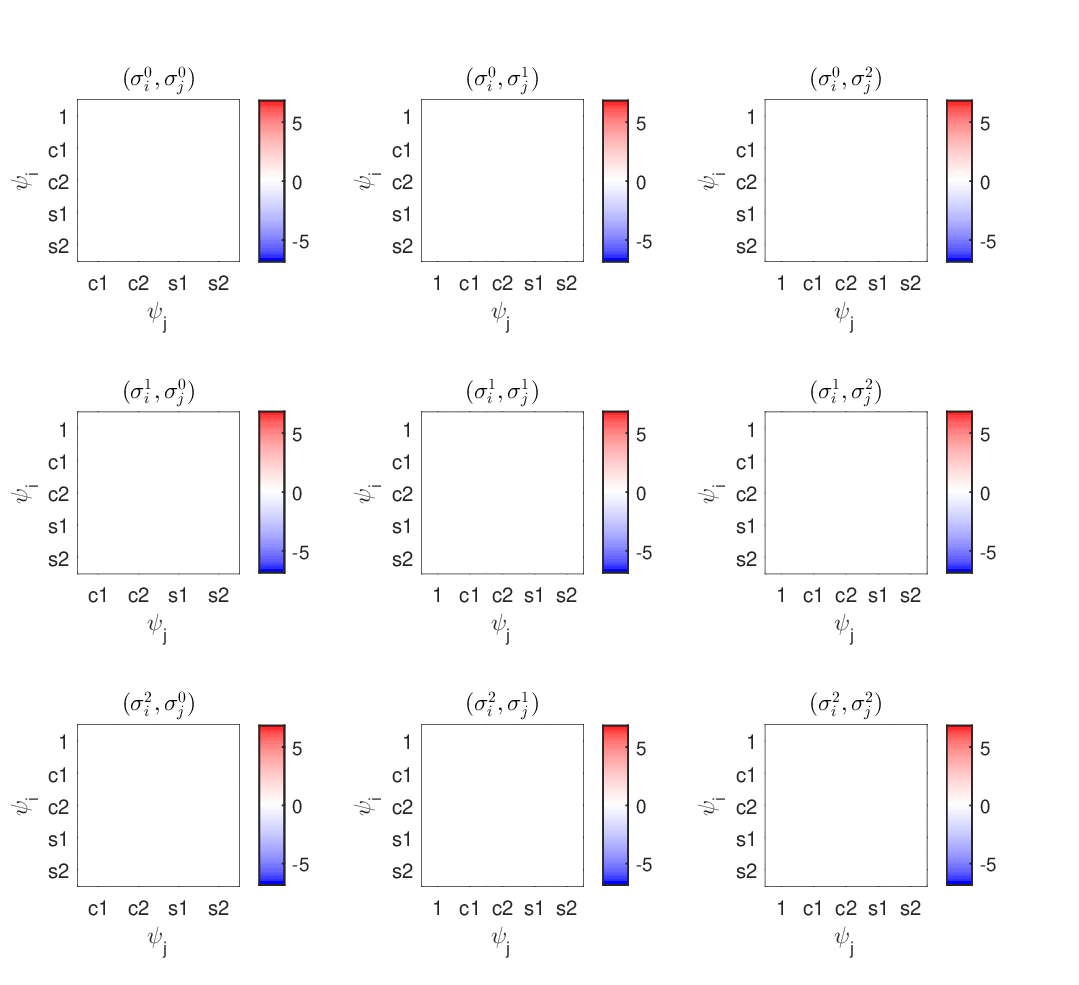}
	\includegraphics[width=0.6\textwidth]{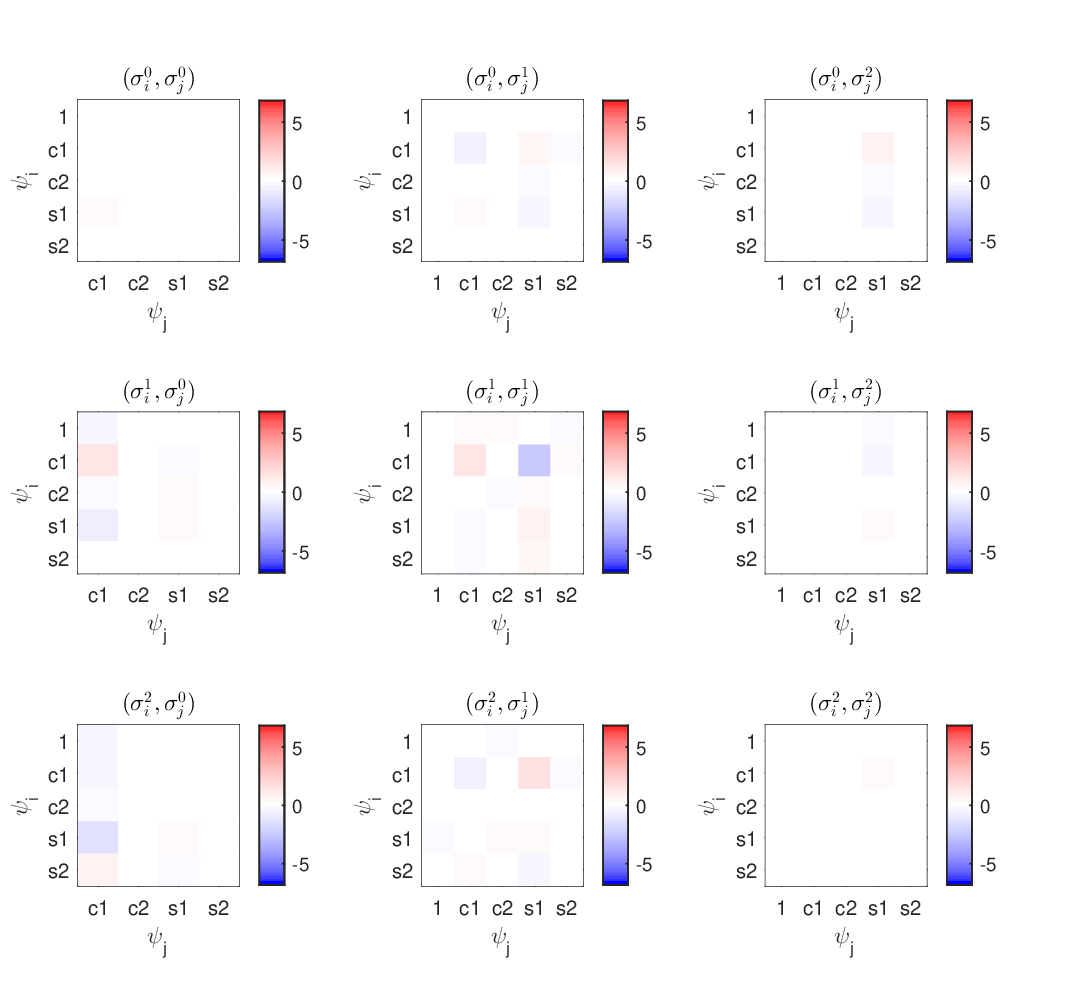}
	\caption{An example of the reconstructed coupling coefficients in the amplitude equation for two uni-directionally coupled Wilson-Cowan models. The upper panels correspond to the first oscillator without input, and the lower panels correspond to the second oscillator, that receives input from the first one. The coefficients are presented as a heat-map (see description in Fig.~\ref{fig:expcoefs} in the main text). \label{smfig:coupl_ampl_wcmod}}
\end{figure}

\begin{figure}[h!]
	\centering
	\includegraphics[width=0.6\textwidth]{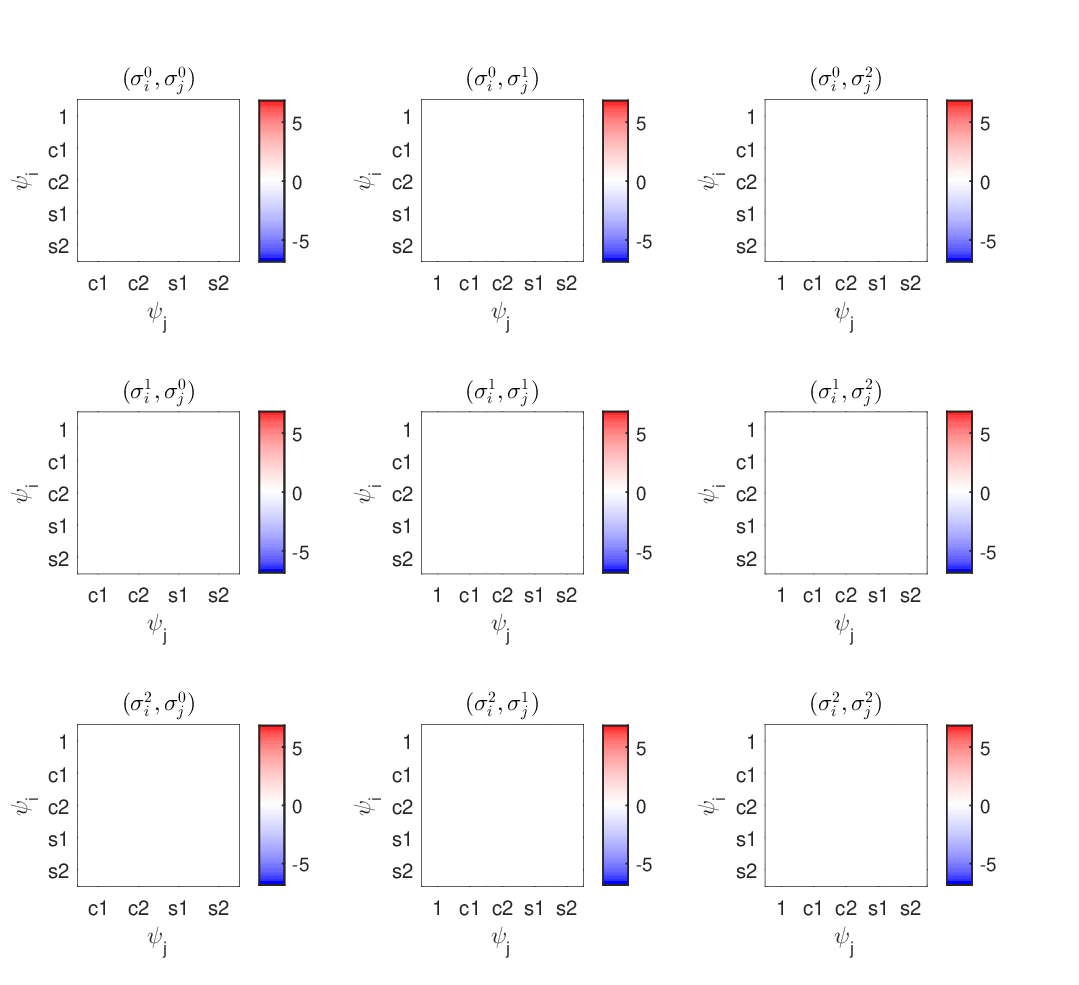}
    \includegraphics[width=0.6\textwidth]{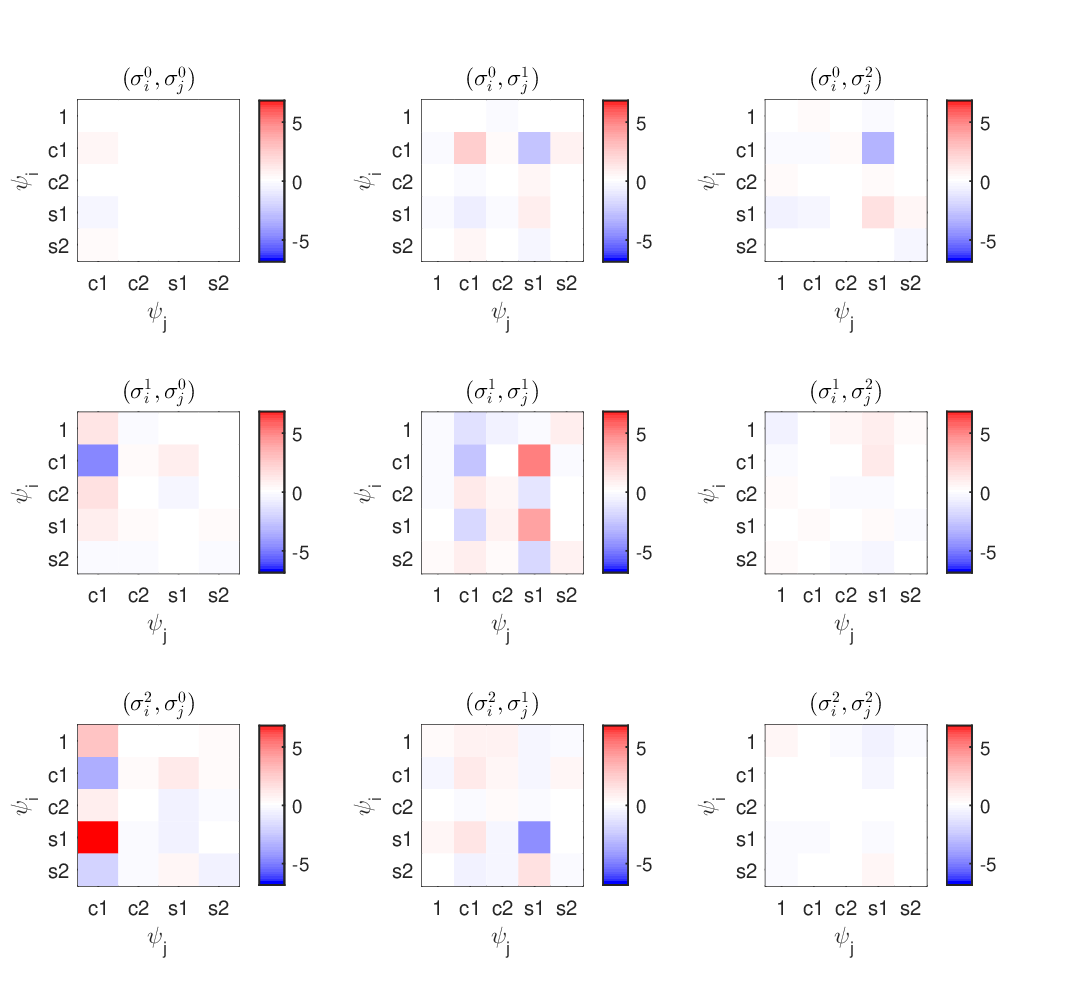}
	\caption{An example of the reconstructed coupling coefficients in the phase equation for two uni-directionally coupled Wilson-Cowan models. The upper panels correspond to the first oscillator without input, and the lower panels correspond to the second oscillator, that receives input from the first one. The coefficients are presented as a heat-map (see description in Fig.~\ref{fig:expcoefs} in the main text). \label{smfig:coupl_phase_wcmod}}
\end{figure}
\clearpage
\newpage
\section{References}